%% file: paper.tex
\documentclass[pdflatex,sn-mathphys-num]{sn-jnl}

\usepackage[utf8]{inputenc}
\usepackage{amsmath}
\usepackage{amsfonts}
\usepackage{amssymb}
\usepackage{amsthm}
\usepackage[dvipsnames,table]{xcolor}
\usepackage{graphicx}
\usepackage{appendix}
\usepackage{comment}

\usepackage{tikz,tikz-3dplot}
\usetikzlibrary{spy}
\usetikzlibrary{matrix}
\usetikzlibrary{arrows}
\usetikzlibrary{pgfplots.groupplots}
\usetikzlibrary{arrows.meta,decorations.markings} 
\usetikzlibrary{external}

\usepackage{booktabs}
\usepackage{multirow}
\usepackage{multicol}

\usepackage{pgfplots}
\usepgfplotslibrary{fillbetween}
\usepgfplotslibrary{groupplots}
\usepackage{subcaption}

\usepackage{bm}
\usepackage{dutchcal}
\usepackage{mathtools}
\usepackage[inline]{enumitem}
\usepackage{hyperref}

\usepackage{dirtree,array}

\usepackage{nicefrac}

\hypersetup{
    colorlinks = true,
}
\usepackage[capitalize]{cleveref}
\crefname{lstlisting}{listing}{listings}
\Crefname{lstlisting}{Listing}{Listings}

\def\ContinueLineNumber{\lstset{firstnumber=last}}
\def\ResetLineNumber{\lstset{firstnumber=1}}

\newcommand{\fullwidth}{0.90\linewidth}
\newcommand{\halfwidth}{0.45\linewidth}

\newcommand{\twothirdwidth}{0.60\linewidth}
\newcommand{\quartwidth}{0.225\linewidth}

\newcommand{\gismo}{\texttt{G+Smo}\:}
\newcommand{\gs}[1]{\lstinline[keepspaces=true,breaklines]|#1|}
\DeclareRobustCommand{\bottombar}{\text{\rotatebox[origin=t]{90}{\reflectbox{$\neg$}}}\:}
\usepackage{physics}
\usepackage{bm}

\newcommand{\R}{\mathbb{R}}

\newcommand{\VEC}[1]{\vb*{#1}}
\newcommand{\MAT}[1]{#1}
\newcommand{\TEN}[1]{\vb{#1}}

\newcommand{\uvec}{\vb*{u}}

\newcommand{\epsten}{\VEC{\varepsilon}}
\newcommand{\kapten}{\VEC{\kappa}}

\newcommand{\Nten}{\vb{N}}
\newcommand{\Mten}{\vb{M}}
\newcommand{\Sten}{\vb{S}}
\newcommand{\Cten}{\vb{C}}
\newcommand{\Ccten}{\VEC{\mathbcal{C}}}

\usepackage{listings}
\lstdefinestyle{all}{
    backgroundcolor=\color{white},
    commentstyle=\color{Green},
    keywordstyle=\color{Magenta},
    numberstyle=\tiny\color{Gray},
    stringstyle=\color{Purple},
    basicstyle=\ttfamily\footnotesize,
    breakatwhitespace=false,
    breaklines=true,
    captionpos=b,
    keepspaces=true,
    numbers=left,
    numbersep=5pt,
    showspaces=false,
    showstringspaces=false,
    showtabs=false,
    tabsize=2
}
\lstdefinestyle{cpp}{
    style=all,
    language=c++,
        escapechar=@,
    morekeywords={gsMultiPatch,gsMultiBasis,gsFunction,gsFunctionSet,gsBasis,gsGeometry,gsConstantFunction,gsOptionList,boundaryInterface,gsFileData,gsExprAssembler,gsBoundaryConditions,
    gsVector,gsMatrix,gsSparseMatrix,gsSparseSolver,
    gsMappedBasis,gsMappedSpline,gsAlmostC1,gsDPatch,gsDPatchBase,gsSurfMesh,gsC1SurfSpline,gsApproxC1Spline,gsMPBESSpline,gsSmoothInterfaces,
    gsMaterialMatrixBase,gsMaterialMatrixBaseDim,gsMaterialMatrixLinear,gsMaterialMatrixNonlinear,gsMaterialMatrixComposite,gsMaterialMatrixMuesli,gsMaterialMatrixTFT,getMaterialMatrix,gsThinShellAssembler,gsThinShellAssemblerDWR,
    gsStructuralAnalysisOps,gsStaticBase,gsStaticDR,gsStaticNewton,gsStaticComposite,gsEigenProblemBase,gsModalSolver,gsBucklingSolver,gsALMBase,gsALMLoadControl,gsALMRiks,gsALMCrisfield,
    Force_t,Matrix_t,Residual_t, ALResidual_t, Jacobian_t, dJacobian_t}
}
\lstdefinestyle{xml}{
    style=all,
    language=XML,
    morekeywords={MaterialMatrix,Thickness,Density,Parameters,Function}
}
\lstdefinestyle{python}{
    style=all,
    language=python
}

\newcommand{\code}[1]{\lstinline|#1|}

\definecolor{tudelft-cyan}{cmyk}{1,0,0,0}
\definecolor{tudelft-black}{cmyk}{0,0,0,1}
\definecolor{tudelft-white}{cmyk}{0,0,0,0}
\definecolor{tudelft-dark-blue}{cmyk}{1,0.8,0.08,0.7}
\definecolor{tudelft-turquoise}{cmyk}{0.72,0,0.24,0}
\definecolor{tudelft-blue}{cmyk}{0.98,0.4,0.0,0}
\definecolor{tudelft-purple}{cmyk}{0.65,1,0,0.12}
\definecolor{tudelft-pink}{cmyk}{0,0.7,0,0}
\definecolor{tudelft-raspberry}{cmyk}{0.05,1,0.48,0.3}
\definecolor{tudelft-red}{cmyk}{0,0.85,0.75,0}
\definecolor{tudelft-orange}{cmyk}{0,0.7,0.75,0}
\definecolor{tudelft-yellow}{cmyk}{0,0.31,0.98,0}
\definecolor{tudelft-light-green}{cmyk}{0.63,0,0.84,0}
\definecolor{tudelft-dark-green}{cmyk}{1,0,0.68,0.04}

\colorlet{TUcol1}{tudelft-cyan}
\colorlet{TUcol2}{tudelft-black}
\colorlet{TUcol3}{tudelft-white}
\colorlet{TUcol4}{tudelft-dark-blue}
\colorlet{TUcol5}{tudelft-turquoise}
\colorlet{TUcol6}{tudelft-blue}
\colorlet{TUcol7}{tudelft-purple}
\colorlet{TUcol8}{tudelft-pink}
\colorlet{TUcol9}{tudelft-raspberry}
\colorlet{TUcol10}{tudelft-red}
\colorlet{TUcol11}{tudelft-orange}
\colorlet{TUcol12}{tudelft-yellow}
\colorlet{TUcol13}{tudelft-light-green}
\colorlet{TUcol14}{tudelft-dark-green}

\colorlet{col0}{TUcol4}
\colorlet{col1}{TUcol9}
\colorlet{col2}{TUcol6}
\colorlet{col3}{TUcol14}
\colorlet{col4}{TUcol7}
\colorlet{col5}{TUcol12}

\colorlet{plotcol1}{TUcol10}
\colorlet{plotcol2}{TUcol12}
\colorlet{plotcol3}{TUcol13}
\colorlet{plotcol4}{TUcol14}
\colorlet{plotcol5}{TUcol6}
\colorlet{plotcol6}{TUcol7}

\pgfplotscreateplotcyclelist{mark list mod}{
    every mark/.append style={thin,solid,fill=\pgfplotsmarklistfill,line width = 0.1},mark size=1.5,mark=*,mark options=solid,fill opacity=0.5\\
    every mark/.append style={thin,solid,fill=\pgfplotsmarklistfill,line width = 0.1},mark size=1.5,mark=square*,mark options=solid,fill opacity=0.5\\
    every mark/.append style={thin,solid,fill=\pgfplotsmarklistfill,line width = 0.1},mark size=1.5,mark=triangle*,mark options=solid,fill opacity=0.5\\
    every mark/.append style={thin,solid,fill=\pgfplotsmarklistfill,line width = 0.1},mark size=1.5,mark=halfsquare*,mark options=solid,fill opacity=0.5\\
    every mark/.append style={thin,solid,fill=\pgfplotsmarklistfill,line width = 0.1},mark size=1.5,mark=pentagon*,mark options=solid,fill opacity=0.5\\
    every mark/.append style={thin,solid,fill=\pgfplotsmarklistfill,line width = 0.1},mark size=1.5,mark=halfcircle*,mark options=solid,fill opacity=0.5\\
    every mark/.append style={thin,solid,fill=\pgfplotsmarklistfill,rotate=180,line width = 0.1},mark size=1.5,mark=halfdiamond*,mark options=solid,fill opacity=0.5\\
    every mark/.append style={thin,solid,fill=\pgfplotsmarklistfill!40,line width = 0.1},mark size=1.5,mark=otimes*,mark options=solid,fill opacity=0.5\\
    every mark/.append style={thin,solid,fill=\pgfplotsmarklistfill,line width = 0.1},mark size=1.5,mark=diamond*,mark options=solid,fill opacity=0.5\\
    every mark/.append style={thin,solid,fill=\pgfplotsmarklistfill,line width = 0.1},mark size=1.5,mark=halfsquare right*,mark options=solid,fill opacity=0.5\\
    every mark/.append style={thin,solid,fill=\pgfplotsmarklistfill,line width = 0.1},mark size=1.5,mark=halfsquare left*,mark options=solid,fill opacity=0.5\\
}

\pgfplotscreateplotcyclelist{linestyles mod}{
    solid\\
    densely dashed\\
    densely dotted\\
    densely dashdotted\\
    densely dashdotdotted\\
    dashed\\
    dotted\\
    dashdotted\\
    dashdotdotted\\
}

\pgfplotscreateplotcyclelist{colorlist mod}{
    plotcol1\\
    plotcol2\\
    plotcol3\\
    plotcol4\\
    plotcol5\\
    plotcol6\\
}

\pgfplotsset{
    cycle multiindex* list={
        mark list mod\nextlist
        colorlist mod\nextlist
        linestyles mod\nextlist
    },
    width =\linewidth,
    height=0.25\textheight,
    grid=major,
    ylabel near ticks,
    xlabel near ticks,
    enlargelimits,
    grid=major,
    legend style={fill=white, fill opacity=0.6, draw opacity=1,draw=black,text opacity=1, font=\footnotesize},
    legend cell align={left},
    legend image with text/.style={
        legend image code/.code={%
            \node[anchor=center] at (0.3cm,0cm) {#1};
        }
    },
    tick label style={font=\footnotesize},
    label style={font=\footnotesize},
}

\colorlet{gray1}{black!20}
\colorlet{gray2}{black!40}
\colorlet{gray3}{black!60}
\colorlet{gray4}{black!80}

\definecolor{col1}{HTML}{0480B0}
\definecolor{col2}{HTML}{FFB600}
\definecolor{col3}{HTML}{FF2C00}
\definecolor{col4}{HTML}{AA6F39}

\pgfplotscreateplotcyclelist{MyCyclelist}{%
thin,solid,TUcol2,mark=o,mark size=1,every mark/.append style={solid,line width = 0.5}\\
thin,densely dashed,TUcol3,mark=square,mark size=1,every mark/.append style={solid,line width = 0.5}\\
thin,densely dotted,TUcol4,mark=triangle,mark size=1,every mark/.append style={solid,line width = 0.5}\\
thin,dotted,TUcol5,mark=diamond,mark size=1,every mark/.append style={solid,line width = 0.5}\\
thin,dashed,TUcol6,mark=x,mark size=1,every mark/.append style={solid,line width = 0.5}\\
thin,loosely dotted,TUcol7,mark=star,mark size=1,every mark/.append style={solid,line width = 0.5}\\
thin,loosely dashed, TUcol8,mark=+,mark size=1,every mark/.append style={solid,line width = 0.5}\\
thin,dashdotted,TUcol9,mark=asterisk,mark size=1,every mark/.append style={solid,line width = 0.5}\\
}

\pgfplotsset{
    layers/legend behind plots/.define layer set={
            axis background,axis grid,axis ticks,axis lines,axis tick labels,main,axis descriptions,axis foreground
    }{
        grid style= {/pgfplots/on layer=axis grid},
        tick style= {/pgfplots/on layer=axis ticks},
        axis line style= {/pgfplots/on layer=axis lines},
        label style= {/pgfplots/on layer=axis descriptions},
        legend style= {/pgfplots/on layer=axis tick labels}, 
        title style= {/pgfplots/on layer=axis descriptions},
        colorbar style= {/pgfplots/on layer=axis descriptions},
        ticklabel style= {/pgfplots/on layer=axis tick labels},
        axis background@ style={/pgfplots/on layer=axis background},
        3d box foreground style={/pgfplots/on layer=axis foreground},
    },
}

\tikzset{%
    style0/.style = {plotcol1,thin,solid,every mark/.append style={solid,fill=\pgfplotsmarklistfill,line width = 0.1},mark size=1.5,mark=*,mark options=solid,fill opacity=0.75,},
    style1/.style = {plotcol2,thin,densely dashed,every mark/.append style={solid,fill=\pgfplotsmarklistfill,line width = 0.1},mark size=1.5,mark=square*,mark options=solid,fill opacity=0.75},
    style2/.style = {plotcol3,thin,densely dotted,every mark/.append style={solid,fill=\pgfplotsmarklistfill,line width = 0.1},mark size=1.5,mark=triangle*,mark options=solid,fill opacity=0.75},
    style3/.style = {plotcol4,thin,densely dashdotted,every mark/.append style={solid,fill=\pgfplotsmarklistfill,line width = 0.1},mark size=1.5,mark=halfsquare*,mark options=solid,fill opacity=0.75},
    style4/.style = {plotcol5,thin,densely dashdotdotted,every mark/.append style={solid,fill=\pgfplotsmarklistfill,line width = 0.1},mark size=1.5,mark=pentagon*,mark options=solid,fill opacity=0.75},
}

\theoremstyle{thmstyleone}%
\theoremstyle{thmstyletwo}%
\newtheorem{example}{Example}%

\theoremstyle{thmstylethree}%

\begin{document}

\title[Isogeometric multi-patch shell analysis using the Geometry + Simulation Modules]{Isogeometric multi-patch shell analysis using the Geometry + Simulation Modules}


\author*[1]{\fnm{Hugo M.} \sur{Verhelst}}\email{hugomaarten.verhelst@unipv.it}
\author[2]{\fnm{Angelos} \sur{Mantzaflaris}}\email{angelos.mantzaflaris@inria.fr}
\author[3]{\fnm{Matthias} \sur{M\"oller}}\email{m.moller@tudelft.nl}
\affil*[1]{\orgdiv{Department of Civil Engineering and Architecture}, \orgname{University of Pavia}, \orgaddress{\street{Via A. Ferrata 3}, \city{Pavia}, \postcode{27100}, \country{Italy}}}
\affil[2]{\orgdiv{AROMATH}, \orgname{Inria Sophia Antipolis - M\'editerran\'ee, Universit\'e C\^ote d\'Azur}, \orgaddress{\street{2004 route des Lucioles - BP 93}, \city{Sophia Antipolis cedex}, \postcode{06902}, \country{France}}}
\affil[3]{\orgdiv{Delft Institute of Applied Mathematics}, \orgname{Delft University of Technology}, \orgaddress{\street{Mekelweg 4}, \city{Delft}, \postcode{2628 CD}, \country{The Netherlands}}}


\abstract{Isogeometric Analysis (IGA) bridges Computer-Aided Design (CAD) and Finite Element Analysis (FEA) by employing splines as a common basis for geometry and analysis. One of the advantages of IGA is in the realm of thin shell analysis: due to arbitrary continuity of the spline basis, Kirchhoff--Love shells can be modeled without the need to introduce unknowns for the mid-plane rotations, leading to a reduction in the number of unknowns. In this paper, we provide the background of an implementation of Isogeometric Kirchhoff--Love shells within the Geometry + Simulation Modules (G+Smo). This paper accompanies multiple previous publications and elaborates on the design of the software used in these papers, rather than the novelty of the methods presented therein. The presented implementation provides patch coupling via penalty methods and unstructured splines, goal-oriented error estimators, several algorithms for structural analysis and advanced algorithms for the modeling of wrinkling in hyperelastic membranes. These methods are all contained in three new modules in G+Smo: a module for Kirchhoff-Love shells, a module for structural analysis, and a module for unstructured spline constructions. As motivated in this paper, the modules are implemented with the aim of being compatible with future developments. For example, by providing base implementations of material laws, by using black-box functions for the structural analysis module, or by providing a standardized approach for the implementation of unstructured spline constructions. Through several examples with code snippets, simple routines are highlighted to illustrate how one can interact with the off-the-shelf routines provided within the modules. Overall, this paper demonstrates that the new modules contribute to a versatile ecosystem for the modeling of multi-patch shell problems through fast off-the-shelf solvers with a simple interface, designed to be extended in future research.}




\maketitle

\section{Introduction}\label{sec:Introduction}

With the advent of Isogeometric Analysis \cite{Hughes2005}, the fields of computer-aided design (CAD) and computer-aided engineering (CAE) slowly unify. By sharing the same mathematical foundation based on splines, highly smooth bases are introduced in classic Galerkin frameworks in CAE, and concepts like analysis suitability, water tightness, and de-featuring become important topics in the CAD community, supporting better analyses \cite{Cohen2010a,Bazilevs2010,Buffa2022a,Buffa2022b}. Since its introduction, many developments have occurred within the isogeometric analysis paradigm, providing a viable alternative to Finite Element Analysis (FEA) in the engineering discipline. With the aim of unifying the CAD and CAE pipelines, IGA substantially influences both communities. For example, new spline constructions with refinement properties suitable for CAE have been developed \cite{Giannelli2012,Giannelli2016,Deng2008,Dokken2013}. In addition, mathematical analysis of the properties of various existing spline constructions and their influence on CAE have been assessed \cite{Buffa2022}. On the other hand, the application of IGA to computational mechanics has increased the interest in problems with high continuity requirements, such as phase-field modeling \cite{Proserpio2020,Proserpio2021,Greco2024a,Greco2025a}, as well as shape-optimisation problems benefiting from the smooth and local geometric parametrisation of spline geometries \cite{Kiendl2014}, e.g., for patient-specific heart valve design \cite{Xu2018}.\\

With the increasing complexity of problems solved with isogeometric analysis, the demand for versatile software libraries increases. Since these software libraries operate in the CAD and CAE domains, advanced geometric and physical modeling capabilities are essential in performing state-of-the-art simulations. Since the advent of IGA in 2005, several software libraries have been developed, ranging from commercial closed-source libraries and free open-source libraries to in-house codes. In the following, an exhaustive yet incomplete overview of publicly available (closed-source and open-source) software libraries with IGA features is presented. On the one hand, IGA software libraries are based on existing FEA libraries that have been extended with IGA capabilities, such as the commercial software LS-DYNA \cite{Hartmann2016} with the ANSA pre-processor \cite{Leidinger2019c} or the open-source libraries Kratos multi-physics \cite{Ferrandiz2023,Dadvand2010} (C++), PetIGA (built on Petc in C++) \cite{Dalcin2016}, MFEM (C++)\cite{Anderson2021} and tIGAr (built on FEniCS in Python) \cite{Kamensky2019}. On the other hand, general-purpose IGA libraries have been presented, such as the open-source Nutils (Python) \cite{Zwieten2022}, GeoPDEs (MATLAB/Octave) \cite{Vazquez2016}, the Geometry + Simulation Modules (\gismo, C++ with Python bindings) \cite{Juttler2014,Mantzaflaris2020} and PSYDAC (Python) \cite{Guclu2022} for general problems, and Bembel (C++) \cite{Dolz2020} for BEM-IGA, YETI (Python) \cite{Duval2023} for structural optimisation and the closed-source software developed by Coreform LLC. In general, the FEA-based libraries are extensive libraries with large communities and heavily optimised routines (for FEA), whereas the IGA libraries have a much lower technology readiness level and emerge mainly from academic research projects.\\

The aim of this paper is to provide an overview of implementation aspects and reproducibility of previously published results \cite{Verhelst2020,Verhelst2021,Farahat2023a,Weinmuller2021,Weinmuller2022,Verhelst2023APALM,Verhelst2023Adaptive,Verhelst2023Coupling,Verhelst2025} through an open-source code for multi-patch shell structural analysis using IGA. The software is provided as modules within the library of the Geometry + Simulation Modules (\gismo), the latter providing basic routines for geometric modeling and isogeometric system assembly. In particular, this paper describes a module for Kirchhoff--Love shell analysis including hyperelasticity \cite{Verhelst2021}, a hyperelastic tension-field-based membrane model \cite{Verhelst2025}, as well as error estimation and mesh adaptivity \cite{Verhelst2023Adaptive}. This \gs{gsKLShell} module can be used as an off-the-shelf solver that does not require the user to implement PDEs manually. For unstructured spline-based multi-patch modeling, this paper presents the \gs{gsUnstructuredSpline} module, which has recently been used to provide a thorough comparison of unstructured spline constructions for isogeometric analysis \cite{Farahat2023a,Farahat2023,Weinmuller2021,Weinmuller2022,Verhelst2023Coupling}. Lastly, this paper elaborates on a module for structural analysis in \gismo named \gs{gsStructuralAnalysis}, which includes the novel Adaptive Parallel Arc-Length Method, presented in \cite{Verhelst2023APALM}. Most of the features presented in the present paper are off-the-shelf high-fidelity routines, aiming to be used in engineering applications. The goal of this paper is to provide a versatile software architecture for the IGA paradigm and demonstrate that demanding problems such as shell analysis can be solved with higher-order IGA methods, both efficiently and with superior quality in the numerical results.\\

In the remainder of this paper, the design of the modules for shell and structural analysis in \gismo are discussed. In \cref{sec:gismo}, an overview of \gismo is given to support the novel modules presented in this paper. \Cref{sec:gsKLShell} elaborates on the shell assembler in particular, describing how the isogeometric Kirchhoff--Love shell formulation is implemented into the module. \Cref{sec:gsStructuralAnalysis} describes the mathematical details behind different structural analysis routines and highlights the features of the implemented module. \Cref{sec:gsUnstructuredSplines} elaborates on different unstructured spline constructions for multi-patch analysis in \gismo. In \cref{sec:Results}, results for different benchmark problems are provided. Some of these results are adopted from previous publications, but in this paper, the aim is to elaborate on the model features behind them. Lastly, \cref{sec:Conclusion} provides a concluding summary.

\section{The Geometry + Simulation Modules}\label{sec:gismo}

Since this paper presents three novel modules for the Geometry + Simulation Modules library (\gismo), a brief overview of \gismo is provided here. The reader is referred to \cref{sec:installation} for download and installation instructions. \gismo is a header-only C++ library dedicated to isogeometric analysis. Being a collection of different modules, \gismo is a versatile library for both geometric modeling and isogeometric analysis, based purely on the mathematical foundations of spline modeling. Leaving an extensive overview of \gismo for another publication, this paper focusses on the relevant features of \gismo for the multi-patch shell analysis modules presented later (see \cref{sec:gsKLShell,sec:gsStructuralAnalysis,sec:gsUnstructuredSplines}).\\

\subsection{Geometric modeling}
Geometric modeling in \gismo can be performed on a wide range of spline bases and corresponding geometries. This flexibility in spline modeling is due to the fact that any basis in \gismo inherits from \gs{gsBasis}, for example tensor B-spline basis (\gs{gsTensorBSplineBasis}), their NURBS counterpart (\gs{gsTensorNurbsBasis}), or (Truncated) Hierarchical B-splines, or (T)HB splines \cite{Giannelli2012,Giannelli2016}, (\gs{gsTHBSplineBasis} and \gs{gsHBSplineBasis}). These abstract basis objects can be collected in a \gs{gsMultiBasis} to perform analysis on complex topologies. Similarly, the \gs{gsGeometry} is an abstract definition of the geometric counterparts of the implemented bases in \gismo, and the \gs{gsMultiPatch} provides a container for them. The \gs{gsBasis} and \gs{gsGeometry}, \gs{gsMultiBasis} and \gs{gsMultiPatch} are widely accepted classes throughout \gismo, meaning that multiple basis types can be used in all implemented solvers, e.g., for solving PDEs.\\

Geometric modeling in \gismo can be used for the pre-processing of geometries and bases for structural analysis with shells. Within \gismo, different spline geometry can be constructed and pre-processed in several ways. Firstly, CAD geometries can be imported through files, e.g., through the \gs{3dm} file format used in Rhinoceros, using the openNURBS plugin \cite{openNURBS}, through the \gs{IGES} format, or via the \gs{Parasolid x_t}. Secondly, geometries can be imported as a mesh, e.g., using the \gs{off} format, which can be transformed to a bi-linear \gs{gsMultiPatch}; see \cref{example:gismo}. Thirdly, geometries can be imported as point clouds and converted to splines using (adaptive) fitting routines \cite{Bracco2022}. Lastly, boundary-represented geometries can be meshed using Coon's patches or using optimisation or PDE-based parameterisation techniques \cite{Ji2021,Ji2022,Hinz2018,Hinz2021}.

\begin{example}[Geometry import using a mesh]
\label{example:gismo}
This example demonstrates the import of a mesh into \gismo, after which it is converted to bi-linear patches. Firstly, the \gs{gsFileData} class is used to read a \gs{gsSurfMesh} from an \gs{off} file; see \cref{fig:gismo_car_mesh} for the mesh:

\ResetLineNumber
\begin{lstlisting}
gsFileData<> fd("car.off");
gsSurfMesh HEmesh;
fd.getFirst<gsSurfMesh>(HEmesh);
\end{lstlisting}

The \gs{gsSurfMesh} uses the half-edge mesh structure. Using the algorithm presented in \cref{fig:coupling_Illustration_patches_from_mesh} \cite{Verhelst2023Coupling}, the half-edge mesh is converted into a smaller number of patches; see \cref{fig:gismo_car_interfaces}.

\begin{figure}
    \centering
    \begin{subfigure}[t]{0.3\linewidth}
        \centering
        \resizebox{\linewidth}{!}{\input{Figures/Illustration_patches_from_mesh1.tikz}}
        \caption{A simple mesh with boundary edges in black and interior edges in gray. The boundary extraordinary vertices (bEVs), i.e., the vertices on a boundary with valence $\nu\geq3$, are denoted by a black circle, and the interior extraordinary vertices (iEVs), i.e., interior vertices with valence $\nu\geq4, \nu\neq4$, are denoted by grey circles.}
        \label{fig:coupling_Illustration_patches_from_mesh1}
    \end{subfigure}
    \hfill
    \begin{subfigure}[t]{0.3\linewidth}
        \centering
        \resizebox{\linewidth}{!}{\input{Figures/Illustration_patches_from_mesh2.tikz}}
        \caption{Illustration of the interface tracing procedure. From each EV, all outgoing edges are traced as illustrated until another EV or a boundary is hit. }
        \label{fig:coupling_Illustration_patches_from_mesh2}
    \end{subfigure}
    \hfill
    \begin{subfigure}[t]{0.3\linewidth}
        \centering
        \resizebox{\linewidth}{!}{\input{Figures/Illustration_patches_from_mesh3.tikz}}
        \caption{Result of interface tracing from all the EVs. Every patch is now bounded by a set of boundary and traced interface curves. All patch corners are corners where a traced interface and/or a boundary edge form a corner. Along the hole, different patches are indicated with different shades of gray. In the part bottom-right of the hole, every face forms a patch since all traced curves denoted by colours intersect with other traced curves.}
        \label{fig:coupling_Illustration_patches_from_mesh3}
    \end{subfigure}
    \caption{Procedure to find a multi-patch segmentation from a given mesh. The original mesh in (a) has 46 vertices, 81 edges, and 45 faces, and the final multi-patch (c) has 20 patches. This figure is adopted from \cite{Verhelst2023Coupling}.}
    \label{fig:coupling_Illustration_patches_from_mesh}
\end{figure}

\ContinueLineNumber
\begin{lstlisting}
gsMultiPatch<> mp = HEmesh.linear_patches();
\end{lstlisting}
\begin{figure}[h]
\centering
\begin{subfigure}{\halfwidth}
\centering
\includegraphics[width=\linewidth]{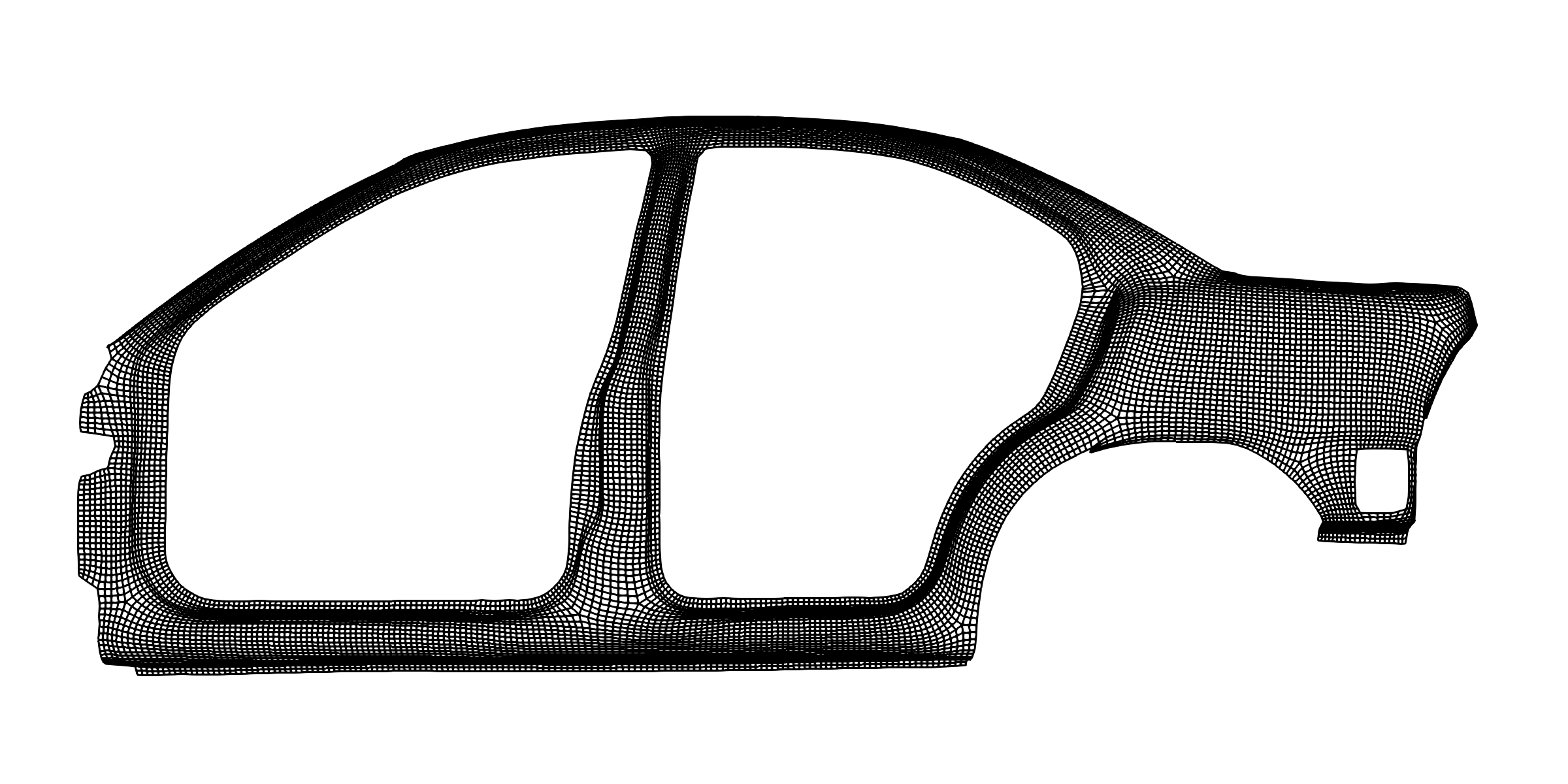}
\caption{Mesh of the geometry stored in \gs{neon_side.off}. The mesh consists of 31086 edges, 15895 vertices, and 62172 faces.}
\label{fig:gismo_car_mesh}
\end{subfigure}
\hfill
\begin{subfigure}{\halfwidth}
\centering
\includegraphics[width=\linewidth]{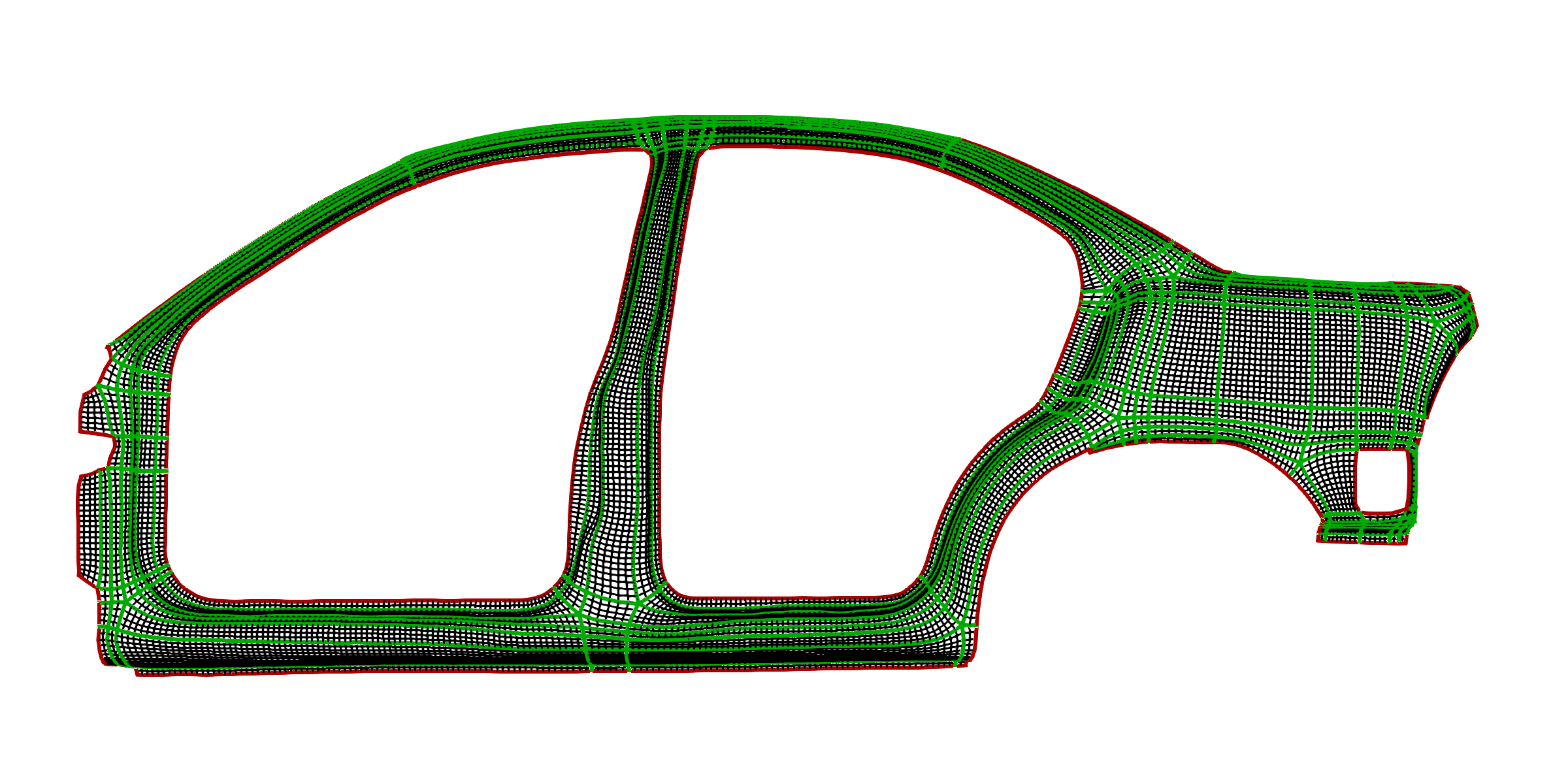}
\caption{Interfaces (green) and boundaries (red) of the bi-linear patches obtained from the mesh in (\subref{fig:gismo_car_mesh})}
\label{fig:gismo_car_interfaces}
\end{subfigure}
\caption{Mesh (\subref{fig:gismo_car_mesh}) and patch boundaries and interfaces (\subref{fig:gismo_car_interfaces}) of the side of a car, read from the file \gs{gismo/filedata/off/neon_side.off}. The plots are obtained via the Paraview \cite{Fabian2011} plotting functions in \gismo.}
\label{fig:coupling_car_side}
\end{figure}
\end{example}

\subsection{Assembly and Linear Algebra}
Besides the geometric modeling capabilities of \gismo, it features routines for the assembly and solving of linear systems for isogeometric analysis. Among the assembly routines available in \gismo is the novel \gs{gsExprAssembler}, used for assembly in the \gs{gsKLShell} module from \cref{sec:gsKLShell}. This assembler allows one to assemble linear systems based on the expressions of a weak formulation. It is shared-memory parallelised using OpenMP and supports exception handling. Furthermore, it uses abstract classes for mathematical functions, geometries (i.e., \gs{gsGeometry}), and bases (i.e., \gs{gsBasis}), such that it can assemble equations on different geometry definitions (e.g., a mathematical representation of a sphere or a THB-spline) using any available basis representation (e.g., \gs{gsBSplineBasis}, \gs{gsTHBSplineBasis}) and any implemented quadrature rule (e.g., the Gau\ss-Lobatto or Gau\ss-Legendre rules, or patch-wise quadrature \cite{Johannessen2017}).\\

For a given linear system, \gismo primarily uses \gs{Eigen} \cite{eigenweb} for linear algebra. For example, the \gs{gsVector}, \gs{gsMatrix}, and \gs{gsSparseMatrix} provide interfaces to vectors and dense and sparse matrices. Furthermore, \gs{gsSparseSolver} provides an interface to linear solvers in \gs{Eigen} and third-party solvers such as \gs{Pardiso} via the Intel$^\text{\textregistered}$ oneAPI Math Kernel Library. Eigenvalue problems can be solved using the eigenvalue solvers available in \gs{Eigen}, and for sparse systems \gismo provides an interface to \gs{Spectra} \cite{Qiu2023}.\\

\section{Kirchhoff-Love Shell Module}\label{sec:gsKLShell}

One of the novel modules presented in this chapter deals with isogeometric Kirchhoff-Love shells, as originally presented in \cite{Kiendl2009}. The module is referred to as the \gs{gsKLShell} module and is featured in recent works of the authors for (i) modeling stretch-based hyperelastic materials in the isogeometric Kirchhoff-Love shell framework \cite{Verhelst2021}, (ii) modeling wrinkling using isogeometric tension-field theory \cite{Verhelst2025}, (iii) the modeling of multi-patch shell problems using several unstructured spline constructions \cite{Farahat2023a,Verhelst2023Coupling}, and (iv) adaptive isogeometric Kirchhoff-Love shell analysis with THB splines \cite{Verhelst2023Adaptive}. As stated in the introduction of this paper (\cref{sec:Introduction}), the \gs{gsKLShell} module is developed within \gismo, hence using the functionalities described in \cref{sec:gismo}.\\

In this section, a brief mathematical background of the isogeometric Kirchhoff-Love shell model is provided; see \cref{subsec:gsKLShell_math}. This background is provided as a reference to show which equations are assembled into systems of equations. Thereafter, \cref{subsec:gsKLSHell_impl} elaborates on the design of the module. \Cref{tab:gsKLShell} in \cref{sec:overview} provides a list of the main classes in the \gs{gsKLShell} module.

\subsection{Mathematical Background}\label{subsec:gsKLShell_math}
The isogeometric Kirchhoff-Love shell was first presented in \cite{Kiendl2009}. More details on this shell model were later provided in \cite{Kiendl2011}, and various extensions of the model have been published by other authors, e.g. stretch-based hyperelasticity for Kirchhoff--Love shells \cite{Verhelst2021}, implicit wrinkling modeling of hyperelastic membranes using tension field theory \cite{Verhelst2025} or goal-oriented error estimators used for mesh adaptivity \cite{Verhelst2023Adaptive}. In this section, a brief derivation of the isogeometric Kirchoff-Love shell formulation is provided for the sake of reference, and the reader is referred to \cite{Verhelst2021,Kiendl2009,Kiendl2015} for a complete overview. Unless specified differently, we denote vectors by lower-case italic bold symbols (e.g., $\VEC{u}$), second-order tensors by upper-case bold symbols (e.g., $\TEN{N}$) and matrices by upper-case italic symbols (e.g., $\MAT{K}$).\\

The variational formulation of the isogeometric Kirchhoff-Love shell model is given by:
\begin{equation}\label{eq:gismo_gsKLShell_VirtualWork0}
\delta_{\uvec} \mathcal{W}^{\text{ext}}(\uvec,\VEC{v}) - \delta_{\uvec} \mathcal{W}^{\text{int}}(\uvec,\VEC{v})= 0,
\end{equation}
where $\delta_{\uvec} \mathcal{W}^{\text{ext}}(\uvec,\VEC{v})$ and $\delta_{\uvec} \mathcal{W}^{\text{int}}(\uvec,\VEC{v})$ are the variations of the external and internal virtual work, respectively, with respect to the displacement field $\uvec$ and a variation $\VEC{v}$, and are given in \cite{Kiendl2009,Kiendl2011}, among others. Using a discretization of $\uvec$ using B-splines, the $r^{\text{th}}$ component of the residual vector $\VEC{R}(\uvec)$ is given by:
\begin{equation}\label{eq:gismo_gsKLShell_Residual}
R_r(\uvec) = \int_\Omega \Nten(\uvec):\epsten_{,r}(\uvec)+\Mten(\uvec):\kapten_{,r}(\uvec)\dd{\Omega} - \int_{\Omega}\VEC{f}\cdot\uvec_{,r}\dd{\Omega} - \int_{\partial\Omega}\VEC{g}\cdot\uvec_{,r}\dd{\Gamma},
\end{equation}
where the product $:$ denotes the inner-product for second-order tensors, $\VEC{\varepsilon}$ and $\VEC{\kappa}$, $\TEN{N}(\uvec)$ and $\TEN{M}(\uvec)$ are the membrane strain, bending strain, force and bending moment tensors, respectively. Taking the second variation of \cref{eq:gismo_gsKLShell_VirtualWork0} with respect to the displacement, the Jacobian matrix $\MAT{K}(\uvec)$ is found, of which the entry with indices $r,s$ is given by:
\begin{equation}\label{eq:gismo_gsKLShell_Jacobian}
K_{rs} = \int_{\Omega} \Nten_{,s}(\uvec):\epsten_{,r}(\uvec) + \Nten(\uvec):\epsten_{,rs}(\uvec) + \Mten_{,s}(\uvec):\kapten_{,r}(\uvec) + \Mten(\uvec):\kapten_{,rs}(\uvec)\dd{\Omega}.
\end{equation}
Here, subscripts ${}_r$ and ${}_s$ denote first variations, and ${}_{rs}$ denotes second variations of the aforementioned tensors; see \cite{Verhelst2021,Kiendl2009,Kiendl2015}. For linear structural analysis, the external force vector $P$ is obtained by
\begin{equation}\label{eq:gismo_externalForce}
P_r = -R_r(\VEC{0}) = \int_{\Omega}\VEC{f}\cdot\uvec_{,r}\dd{\Omega} + \int_{\partial\Omega}\VEC{g}\cdot\uvec_{,r}\dd{\Gamma},
\end{equation}
and the linear stiffness matrix $\MAT{K}^L$ is obtained by
\begin{equation}\label{eq:gismo_gsKLShell_Linear}
K^L_{rs} = K_{rs}(\VEC{0}) = \int_{\Omega} \Nten_{,s}(\VEC{0}):\epsten_{,r}(\VEC{0}) + \Mten_{,s}(\VEC{0}):\kapten_{,r}(\VEC{0})\dd{\Omega}.
\end{equation}
For quasi-static and dynamic analyses, the external force vector and residual contain an extra parameter. In particular, quasi-static analysis involves the analysis of the problem with respect to a load scaled with load magnification factor $\lambda$, i.e., $\VEC{f}\equiv\lambda\VEC{f}$ or $\VEC{g}\equiv\lambda\VEC{g}$, yielding $\VEC{R}(\uvec,\lambda)$. For dynamic loading, the loads can be time-dependent, i.e., $\VEC{f}\equiv\VEC{f}(t)$ or $\VEC{g}\equiv\VEC{g}(t)$, yielding $\VEC{R}(\uvec,t)$. Furthermore, neglecting rotational inertia, the mass matrix $\MAT{M}$ for a Kirchhoff-Love shell with constant density $\rho$ and thickness $t$ is:
\begin{equation}\label{eq:gismo_massMatrix}
M_{rs} = \rho t \int_{\Omega} u_{,r}u_{,s}\dd{\Omega}.
\end{equation}

\subsection{Implementation}\label{subsec:gsKLSHell_impl}
The implementation of the isogeometric Kirchhoff-Love shell model is done in the \gs{gsKLShell} module in \gismo. The module consists of several classes, of which the most important are highlighted in \cref{tab:gsKLShell}; for a complete overview, the reader is referred to the documentation of \gismo. In general, the module is designed so that it is fast, versatile, and forward-compatible, particularly because of the aspects mentioned below.\\

Firstly, the \gs{gsThinShellAssembler} is designed in such a way that the constitutive law is decoupled from the \gs{gsThinShellAssembler}. The constitutive law, provided by the \gs{gsMaterialMatrixBase} and its derived classes, operates as a black box within the \gs{gsThinShellAssembler}. As a consequence, constitutive models can be developed as a family of models (i.e., like the \gs{gsBasis}-family), providing forward compatibility with user-defined or new developments in constitutive models.\\

Secondly, the \gs{gsThinShellAssembler} uses the \gs{gsExprAssembler} for matrix assembly, providing flexibility with respect to the basis and geometry types, as well as the quadrature rules. This provides versatility in the basis and geometry types that can be used to assemble the KL shell (e.g., THB splines); it uses optimised parallel routines and patch rules from the \gs{gsExprAssembler}; and it provides forward compatibility with respect to any future development that will be made to improve the \gs{gsExprAssembler}. In addition, since all operations in the \gs{gsThinShellAssembler} rely on generic \gs{gsBasis} and \gs{gsGeometry} classes, it is straightforward to apply the \gs{gsThinShellAssembler} on geometries from other file types such as \texttt{IGES} or \texttt{3dm}.\\

Lastly, the module uses \gismo's nightly build checks and unit-testing framework. The unit-tests include the uniaxial tension test and the inflated balloon test from \cite{Verhelst2021} for non-linear materials and material models based on tension-field theory \cite{Verhelst2025}. Furthermore, multi-patch capabilities are tested together with the \gs{gsUnstructuredSplines} module presented in \cref{sec:gsUnstructuredSplines}.

\begin{example}[Non-linear analysis using Kirchhoff-Love shells]
\label{example:gsKLShell}
To show the use of the \gs{gsThinShellAssembler} and the \gs{gsMaterialMatrix}, an example code is provided. It is assumed that the geometry (\gs{mp}), the basis (\gs{basis}), the boundary conditions (\gs{bcs}), the body force (\gs{force}), the optional point loads (\gs{pointLoads}), and the material parameters (\gs{E}, \gs{nu}, and \gs{rho}) are provided. Then, the \gs{gsThinShellAssembler} and the \gs{gsMaterialMatrix} are simply constructed as:

\ResetLineNumber
\begin{lstlisting}
// Provided mp, basis, bcs, force, pointLoads, E, nu, rho
gsMaterialMatrixLinear<3,real_t> materialMatrix(mp,thickness,E,nu,rho)
// Create a shell assemblers, with template parameters <target dimension, double-type, assemble bending stiffness>
gsThinShellAssembler<3,real_t,true> assembler(mp,basis,bcs,force,materialMatrix);
assembler.addPointLoads(pointLoads);
\end{lstlisting}

Using the assembler, a linear system can be assembled using \gs{assemble()}, and it can be solved using sparse linear solvers to obtain the linear solution vector, \gs{solution}. Then, the Jacobian matrix and the residual vector are obtained by using \gs{assembleMatrix(solution)} and \gs{assembleVector(solution)}, respectively, which can be used in Newton-Raphson iterations (see \cref{eq:gismo_SA_Newton}) as shown below.

\ContinueLineNumber
\begin{lstlisting}
assembler.assemble(); // Assembles @\cref{eq:gismo_gsKLShell_Linear}@ and @\cref{eq:gismo_externalForce}@

gsSparseSolver<real_t>::LU solver(assembler.matrix());

// Solve the linear static problem (@\cref{eq:gismo_SA_Newton}@)
gsVector<real_t> solution = solver.solve(assembler.rhs());

// Solve the non-linear static problem (@\cref{eq:gismo_SA_Newton}@)
for (index_t it = 0; it!=10; it++)
{
    assembler.assembleMatrix(solution);        // Assembles @\cref{eq:gismo_gsKLShell_Jacobian}@
    assembler.assembleVector(solution);        // Assembles @\cref{eq:gismo_gsKLShell_Residual}@
    solver.compute(assembler.matrix());        // Factorises the matrix
    solution += solver.solve(assembler.rhs()); // Computes the update
}
\end{lstlisting}

Using \gs{solution}, a \gs{gsMultiPatch} of the deformed shell can be constructed for post-processing:

\ContinueLineNumber
\begin{lstlisting}
gsMultiPatch<> mp_def = assembler.constructSolution(solution);
\end{lstlisting}
\end{example}

\section{Structural Analysis Module}\label{sec:gsStructuralAnalysis}

The \gs{gsThinShellAssembler} from the \gs{gsKLShell} module, see \cref{sec:gsKLShell}, can provide matrices and vectors resolving Kirchhoff-Love shell mechanics in isogeometric analysis. For different structural analysis applications, these matrices and vectors are used in a different way, depending on the type of analysis that is solved. For example, the analysis of eigenfrequencies of a structure requires modal analysis, and the analysis of structural stability requires linear buckling or quasi-static analysis. The \gs{gsStructuralAnalysis} module implements routines for structural analysis in \gismo. The main purpose of the module is to provide a broad range of structural analysis routines that can be reused in several applications. Moreover, the module is also home to the novel Adaptive Paralllel Arc-Length Method (APALM), as published in \cite{Verhelst2023APALM}.\\

In this section, a brief mathematical background of different structural analysis routines is provided in \cref{subsec:SA_math}, for reference in later presented functions. Thereafter, \cref{subsec:SA_impl} elaborates on the design of the \gs{gsStructuralAnalysis} module in \gismo. In general, different structural analyses require different types of operators, for example, the linear stiffness matrix $\MAT{K}$ or the force residual vector $\VEC{R}(\uvec)$, depending on the discrete displacements $\uvec$. \Cref{tab:operators} lists the discrete operators relevant in the structural analysis module, together with their type and a description. As will be explained in \cref{subsec:SA_impl}, the code is designed in such a way that the operators are independent of the assemblers or discretisation method, such that the \gs{gsStructuralAnalysis} module can be used with the \gs{gsKLShell} module for Kirchhoff-Love shells, but also with the \gs{gsElasticity} module for solids or with other modules providing discrete operators. \Cref{tab:gsStructuralAnalysis} in \cref{sec:overview} provides an overview of the most important classes in the \gs{gsStructuralAnalysis} module.

\begin{table}[h]
\centering
\caption{Operators for structural analysis, with their corresponding types in \gismo and a description linking the operators to the mathematical theory for shells, see \cref{sec:gsKLShell}.}
\label{tab:operators}
\footnotesize
\begin{tabular}{p{0.1\linewidth}p{0.15\linewidth}p{0.65\linewidth}}
\toprule
Operator & Type & Description\\ \midrule
$\VEC{P}$& \lstinline|Force_t|& The vector of external forces, see \cref{eq:gismo_externalForce}.\\ \midrule
$\VEC{R}(\uvec)$\newline$\VEC{R}(\uvec,\lambda)$ & \lstinline|Residual_t|\newline\lstinline|ALResidual_t| & The vector of residual forces, depending on the discrete displacement vector $\uvec$ and optionally depending on the load factor $\lambda$, see \cref{eq:gismo_gsKLShell_Residual}.\\ \midrule
$\MAT{K}_L$ & \lstinline|Matrix_t| & The linear stiffness matrix, see \cref{eq:gismo_gsKLShell_Linear}.\\ \midrule
$\MAT{M}$ & \lstinline|Matrix_t| & The mass matrix, see \cref{eq:gismo_massMatrix}.\\ \midrule
$\MAT{K}(\uvec)$\newline$\MAT{K}(\uvec,\Delta\uvec)$ & \lstinline|Jacobian_t|\newline\lstinline|dJacobian_t| & The Jacobian matrix, depending on the discrete displacement vector $\uvec$, see \cref{eq:gismo_gsKLShell_Jacobian}, and optionally on the displacement vector increment $\Delta\uvec$ for the Mixed Integration Point (MIP) method \cite{Leonetti2018}.\\
\bottomrule
\end{tabular}
\end{table}


\subsection{Mathematical Background}\label{subsec:SA_math}
A structural analysis problem typically involves a number of different operators providing matrices or vectors to construct a system to solve. As shown in \cref{subsec:gsKLShell_math}, the residual vector (\cref{eq:gismo_gsKLShell_Residual}), the Jacobian (\cref{eq:gismo_gsKLShell_Jacobian}), and the linearized Jacobian and force vector (\cref{eq:gismo_gsKLShell_Linear}) can, for example, be derived. More precisely, the different linear operators that are required for performing the structural analyses described in the remainder of this section are given in \cref{tab:operators}. This section provides a concise overview of the mathematics behind the structural analysis solvers implemented in the \gs{gsStructuralAnalysis} module.

\subsubsection{Static Analysis}
For (non-linear) static analysis, a time-independent boundary value problem with constant load magnitudes is solved. The most common approach is to solve the residual equation, i.e., the variational formulation from \cref{eq:gismo_gsKLShell_VirtualWork0}, using the Newton-Raphson method. This requires the discrete residual and the Jacobian of the system, such that the displacements $\uvec$ are found incrementally by:
\begin{equation}\label{eq:gismo_SA_Newton}
K(\uvec^i)\Delta\uvec = -\VEC{R}(\uvec^i),\quad\uvec^{i+1} = \uvec^i + \Delta\uvec,\quad i = 0,1,...
\end{equation}
Here, $\uvec^{i}$ is the solution in iteration $i$, and $\Delta\uvec$ is the solution increment. The initial solution $\uvec^0$ can be initialised by solving linear static analysis:
\begin{equation}\label{eq:gismo_SA_Linear}
K_L \uvec = \VEC{P},
\end{equation}
where $\MAT{K}_L=\MAT{K}(\VEC{0})$ and $\VEC{P}$ is the external load vector. The Newton-Raphson method is a well-known procedure to solve non-linear systems of equations. Although it requires the assembly of the Jacobian matrix and solving a linear system of equations in every iteration, it converges quadratically.\\

Alternative to the Newton-Raphson method, the Dynamic Relaxation method \cite{Otter1960,Otter1965,Otter1966} provides an explicit method for the computation of static problems with slow convergence but low costs per iteration. It solves the structural dynamics equation by expanding the accelerations using finite differences and mass-proportional damping:
\begin{equation}\label{eq:preliminaries_DR_equation}
    M\ddot{\uvec}(t)+C\dot{\uvec}(t)-\VEC{R}(\uvec(t))=\VEC{0},
\end{equation}
where $M$ is the mass matrix, $\ddot{\uvec}$ is the vector of discrete accelerations, $C$ is the damping matrix, $\dot{\uvec}$ is the vector of discrete velocities, and $\VEC{R}(\uvec)$ is the residual vector. Defining the velocity vector and a time step by $\dot{\uvec}$ and $\Delta t$, respectively, the acceleration vector $\ddot{\uvec}$ is obtained by central finite differences:
\begin{equation}\label{eq:preliminaries_DR_velocity}
    \ddot{\uvec}_{t} = \frac{\dot{\uvec}_{t+\Delta t / 2} - \dot{\uvec}_{t-\Delta t / 2}}{\Delta t},
\end{equation}
using $u_{t}=u(t)$ for the sake of clarity. A common assumption in dynamic relaxation methods is to define the damping matrix proportional to the mass matrix, i.e., $C = cM$ or by other appropriate scaling techniques \cite{Rodriguez2011,Rezaiee-Pajand2017}. However, to avoid the introduction of this extra damping parameter, the \emph{kinetic damping} approach \cite{Cundall1976} traces the kinetic energy in the system and sets the nodal velocity to zero when a peak in kinetic energy is detected yielding a robust method without damping parameter \cite{Barnes1988,Barnes1999,Shugar1990}. Let the kinetic energy in the system be defined by:
\begin{equation}
    E^K_{t} = \frac{1}{2}\dot{\uvec}^\top M\dot{\uvec}.
\end{equation}
Then, a peak is detected if $E^K_{t}> E^K_{K,t+\Delta t}$ for $\Delta t > 0$. Assuming the peak occurs in the middle of the interval $[t-\Delta t,t]$, hence at $t-\Delta t/2$. At this time instance, it holds that $E^K_{t-3\Delta t/2}<E^K_{t-\Delta t /2 }$ and $E^K_{t-\Delta t /2} > E^K_{t+\Delta t/2}$ \cite{Topping1994}. Furthermore, since central finite differences are used, the displacement vector $\uvec_{t+\Delta t}$ and velocity vector $\dot{\uvec}_{t+\Delta t /2}$ are known and the displacements at the peak can be computed by: 
\begin{equation}\label{eq:preliminaries_DR_displacement}
    \uvec_{t^\star} = \uvec_{t + \Delta t} - \frac{3}{2}\dot{\uvec}_{t+\Delta t/2} + \frac{\Delta t}{2}M^{-1}\VEC{R}(\uvec_{t}),
\end{equation}
Where the peak time is denoted by $t^\star$. Using the displacement vector $\uvec_{t^\star}$, the method is re-initiated using $\uvec_{t^\star}$.
Since the velocities are fully damped after a kinetic energy peak, they are set to zero upon re-initialization. Hence, to compute the next step after the restart at a peak on $\uvec_{t^\star}$, the velocity vector for $\uvec_{t^\star+\Delta t/2}$ becomes:
\begin{equation}\label{eq:preliminaries_DR_peak_velo}
    \dot{\uvec}_{\uvec_{t^\star+\Delta t / 2}} = \frac{\Delta t}{2}M^{-1}\VEC{R}(\uvec_{t^\star})
\end{equation}
Using \cref{eq:preliminaries_DR_peak_velo}, the displacement vector $\uvec_{t^\star+\Delta t}$ can be found using \cref{eq:preliminaries_DR_displacement}. This kinematic damping procedure is successfully applied in \cite{Barnes1988,Barnes1999,Taylor2014,Rezaiee-Pajand2017}, among others, showing the robustness of the method while eliminating the need to determine the damping coefficient $c$. \\

The equations above can be solved if and only if the mass matrix $M$ is invertible. However, having a full mass matrix drastically increases the computational costs of the method, since it involves solving a linear system. Therefore, alternatives include to select $M=\rho\:\text{diag}(K)$ with linear stiffness matrix $K$ \cite{Papadrakakis1981} to use a diagonal lumped mass matrix \cite{Joldes2010,Joldes2011,Barnes1988,Barnes1999,Topping1994,Alic2016}
or to use a column-sum of the stiffness matrix \cite{Underwood1983}.

\subsubsection{Modal Analysis}
Modal analysis provides the eigenfrequencies and eigenmodes of a structure in free vibration. Starting from the structural dynamics equation, assuming no damping and a harmonic solution $\uvec(t)=\uvec_A e^{i\omega t}$, the structural dynamics equations simplify to a generalised eigenvalue problem:
\begin{equation}\label{eq:gismo_ModalAnalysisProblem}
\omega^2 M \VEC{\phi} = K \VEC{\phi}
\end{equation}
With $\omega$ and $\VEC{\phi}$ the eigenfrequency and modeshape, respectively. Given $\MAT{M}$ and $\MAT{K}_L$ as discrete linear operators with $N$ degrees of freedom, the solution to \cref{eq:gismo_ModalAnalysisProblem} consists of $N$ eigenpairs $(\omega_k,\VEC{\phi}_k)$. For large $N$, $\MAT{K}_L$ and $\MAT{M}$ are sparse systems, and methods like the Shifted Block Lanczos Algorithm \cite{Grimes1994} can be used to find eigenpairs from a sparse general eigenvalue problem within a range.

\subsubsection{Linear Buckling Analysis}\label{subsec:gsStructuralAnalysis_buckling}
Linear buckling analysis aims to find the critical load in a certain load configuration. Mathematically, this means that one needs to find the point where the stability of the structure changes, coinciding with the solution $\uvec$ for which $\det\qty( \MAT{K}(\uvec))=\VEC{0}$. In brief, linear buckling analysis requires a linear static solution $\uvec_L$ of the problem $\MAT{K}_L\uvec_L=\VEC{P}$ in order to solve the following eigenvalue problem:
\begin{equation}\label{eq:gismo_gsStructuralAnalysis_BucklingAnalysisProblem}
\MAT{K}_L\VEC{\phi}_k = \lambda_k \qty(\MAT{K}(\uvec_L)-\MAT{K}_L)\VEC{\phi}_k.
\end{equation}
Here, the sub-script ${\cdot}_L$ is used to emphasise the linear stiffness matrix $\MAT{K}_L$, and the eigenpair $(\lambda_k,\VEC{\phi}_k)$ consists of the load magnification factor $\lambda_k$ and the corresponding mode shape $\VEC{\phi}_k$ for the $k^{\text{th}}$ buckling mode, $k=1,\dots,N$ for a system with $N$ degrees of freedom. For more information regarding buckling analysis, the reader is referred to \cite{Brendel1982} among other references.

\subsubsection{Post-Buckling Analysis}\label{sec:gsStructuralAnalysis_ArcLength}
Post-buckling analysis involves the quasi-static analysis of a problem with time-independent but varying loads. The goal of quasi-static analysis is to solve $\VEC{R}(\uvec,\lambda)=\VEC{0}$, forming a so-called equilibrium path in the $\VEC{w}=(\lambda,\uvec)$-space. Instead of varying $\lambda$ (load-control) or parts of $\uvec$ (displacement-control), arc-length methods find the solution to $\VEC{R}(\uvec,\lambda)=\VEC{0}$ by varying $\uvec$ and $\lambda$ simultaneously, providing a method that allows to find paths with snap-through behaviour. Given $N$ degrees of freedom, arc-length methods solve a system of $N+1$ equations iteratively by adding the constraint equation $f(\Delta\uvec,\Delta\lambda)=0$ to the incremental solution $\Delta\VEC{w}=(\Delta\uvec,\Delta\lambda)$. For example, the constraint equation corresponding to Crisfield's \cite{Crisfield1981} method is:
\begin{equation}\label{eq:gismo_gsStructuralAnalysis_Crisfield}
f(\Delta\uvec,\Delta\lambda) = \Delta\uvec^\top\Delta\uvec+\Delta\lambda^2\VEC{P}^\top\VEC{P} - \Delta\ell^2 = 0
\end{equation}
Crisfield's method requires solving a quadratic equation, and therefore two solutions are typically found. In \cite{Verhelst2020}, a brief summary of techniques for robust arc-length analysis using Crisfield's method is provided, e.g., the handling of complex roots.\\
In case a singular point is detected, e.g., by monitoring the sign of the determinant of the Jacobian at every solution point, the singular point can be computed by solving the following system of equations \cite{Wriggers1988}:
\begin{equation}\label{eq:gismo_gsStructuralAnalysis_CrisfieldExtended}
\begin{bmatrix}
\VEC{R}(\uvec,\lambda)\\
K(\uvec)\VEC{\phi}\\
l(\VEC{\phi})
\end{bmatrix}
=\VEC{0}.
\end{equation}
Here, the first equation consists of the discrete residual of the system (see \cref{eq:gismo_gsKLShell_Residual}), and the second equation implies finding the singular point by using the equivalence $\det K(\uvec)=0\iff K(\uvec)\VEC{\phi}=0$ for eigenvector $\VEC{\phi}$, normalised using the condition $l(\VEC{\phi})=\Vert\VEC{\phi}\Vert-1=0$. Since the solution at the singular point is the buckling mode shape, this so-called \emph{extended arc-length method} can be used to compute post-buckling simulations without \emph{a priori} perturbation of the solution to switch branches described in \cite{Verhelst2020}.\\

\subsection{Implementation}\label{subsec:SA_impl}
As mentioned in the beginning of this \cref{sec:gsStructuralAnalysis}, the \gs{gsStructuralAnalysis} module is implemented in a way such that it is independent of the discretisation method. In other words, the \gs{gsStructuralAnalysis} module requires the user to provide small interfaces to their discretisation methods, which will be used inside the module for the required computations. These discrete operators, summarised in \cref{tab:operators}, can be used in the solvers provided in the \gs{gsStructuralAnalysis} module, as listed in \cref{tab:gsStructuralAnalysis}. Furthermore, all routines in the module support exception handling, such that assembly errors or other errors can be passed robustly.\\

The operators listed in \cref{tab:operators}, e.g., \gs{ALResidual_t}, are used in different solvers throughout the \gs{gsStructuralAnalysis} module and have different corresponding C++ types. The types \gs{Force_t} and \gs{Matrix_t} are straightforward, since they correspond to \gs{gsVector<T>} and \gs{gsMatrix<T>}. However, the other types require function arguments, e.g., the displacement vector $\uvec$ or the load magnification factor $\lambda$, and are therefore passed using functions. For example, the arc-length residual $\VEC{R}(\uvec,\lambda)$ or the Jacobian $\MAT{K}(\uvec)$ have the following types (defined in \gs{gsStructuralAnalysisOps.h}):

\begin{lstlisting}
typedef std::function < bool ( gsVector<T> const &, T,    gsVector<T> & ) > ALResidual_t;
typedef std::function < bool ( gsVector<T> const &, gsSparseMatrix<T> & ) > Jacobian_t;
\end{lstlisting}
These definitions show that \gs{ALResidual_t} is a function that returns a boolean (\gs{true} if the assembly is successful, \gs{bool}) and takes the vector of displacments ($\uvec$, \gs{gsVector<T> const &}), a real number for the load magnification factor ($\lambda$, \gs{T}), and in-place the residual vector ($\VEC{R}(\uvec,\lambda)$, \gs{gsVector<T> &}) as input arguments. Similarly, the type \gs{Jacobian_t} is a function returning a boolean for successful assembly that takes the vector of displacements ($\uvec$, \gs{gsVector<T> const &}) and returns the Jacobian matrix in-place ($\MAT{K}(\uvec)$, \gs{gsSparseMatrix<T> & }). For the other operators with function arguments in \cref{tab:operators}, the definitions work in a similar way. Eventually, the definition of an operator depends on the interface with the discretisation method. In \cref{example:gsStructuralAnalysis}, an example is provided for the operator definitions using the \gs{gsThinShellAssembler} to initialise the \gs{gsStaticNewton} solver.\\

\begin{example}[Linear operators for a non-linear static analysis using Kirchhoff-Love shells]
\label{example:gsStructuralAnalysis}
When using the \gs{gsThinShellAssembler} from the \gs{gsKLShell} module, see \cref{sec:gsKLShell}, the linear stiffness matrix $\MAT{K}$ and the external force vector $\VEC{P}$ are simply assembled by calling \gs{assemble()} and obtained by calling \gs{matrix()} and \gs{rhs()}, respectively. As shown in \cref{example:gsKLShell}, the Jacobian and residual are assembled using \gs{assembleMatrix(solution)} and \gs{assembleVector(solution)}, respectively. Using this syntax, the function \gs{ALResidual} of type \gs{ALResidual_t} can be defined as follows:

\ContinueLineNumber
\begin{lstlisting}
gsStructuralAnalysisOps<T>::ALResidual_t ALResidual;
ALResidual = [&assembler,&Force](gsVector<T> const &u, T lam, gsVector<T> & result)
{
    // Assemble eq. @\ref{eq:gismo_gsKLShell_Residual}@
    ThinShellAssemblerStatus status = assembler->assembleVector(u);
    result = Force - lam * Force - assembler->rhs();
    return status == ThinShellAssemblerStatus::Success;
};
\end{lstlisting}
Here, it can be seen that the \gs{ThinShellAssemblerStatus} returned by the assembly in the \gs{gsThinShellAssembler} is used to verify whether the assembly is successful. Furthermore, an external force vector \gs{Force} is used together with the residual $\VEC{R}(\uvec)=\VEC{F}_{\text{int}}-\VEC{F}_{\text{ext}}$ to define $\VEC{R}(\uvec,\lambda)=\VEC{R}(\uvec)+\VEC{F}_{\text{ext}}-\lambda \VEC{F}_{\text{ext}}$. Similarly, the Jacobian operator can be defined as:

\ContinueLineNumber
\begin{lstlisting}
gsStructuralAnalysisOps<real_t>::Jacobian_t Jacobian;
Jacobian = [&assembler](gsVector<real_t> const &u, gsSparseMatrix<real_t> & m)
{
    // Assemble eq. @\ref{eq:gismo_gsKLShell_Jacobian}@
    ThinShellAssemblerStatus status = assembler->assembleMatrix(u);
    m = assembler->matrix();
    return status == ThinShellAssemblerStatus::Success;
};
\end{lstlisting}
The \gs{gsThinShellAssembler} also provides an implementation of the Mixed Integration Method (MIP) \cite{Leonetti2018}, which is a method to accelerate Newton iterations by decoupling the stress update from the displacements. In this case, the Jacobian matrix is given by $\MAT{K}(\uvec,\Delta\uvec)$ and defined as follows:
\ContinueLineNumber
\begin{lstlisting}
gsStructuralAnalysisOps<real_t>::dJacobian_t dJacobian;
dJacobian = [&assembler](gsVector<real_t> const &u,
                         gsVector<real_t> const &Du,
                         gsSparseMatrix<real_t> & m)
{
    ThinShellAssemblerStatus status = assembler->assembleMatrix(u,Du);
    m = assembler->matrix();
    return status == ThinShellAssemblerStatus::Success;
};
\end{lstlisting}
The \gs{gsStaticNewton} and \gs{gsALMCrisfield} solvers of the \gs{gsStructuralAnalysis} module are compatible with the MIP method. Using the external load factor, the arc-length residual and Jacobian operators, $\VEC{P}$, $\VEC{R}(\uvec,\lambda)$, and $\MAT{K}(\uvec)$, respectively, an arc-length method can be setup. For example, a Crisfield arc-length solver is defined as follows:

\ContinueLineNumber
\begin{lstlisting}
// Provided Force, ALResidual, Jacobian, options, dL
// Define a Crisfield method, with eq. @\ref{eq:gismo_gsStructuralAnalysis_Crisfield}@
gsALMCrisfield<T> arcLength(Force,ALResidual,Jacobian);
// Set the length
arcLength.options().setReal("Length",dL);
\end{lstlisting}
\end{example}

\section{Unstructured Splines Module}\label{sec:gsUnstructuredSplines}


The unstructured splines module is intensively used in the previous work \cite{Verhelst2023Coupling}, where qualitative and quantitative comparisons are given of multiple unstructured spline routines. Furthermore, the work \cite{Marsala2022a} develops a $G^1$ construction constructed in third-party software but imported into \gismo as a combination of a sparse matrix and a local basis to perform different kinds of analyses. To the best of the author's knowledge, \gismo has a unique position in unstructured splines modeling, being the only open-source code for isogeometric analysis and providing an interface for the development of unstructured splines together with advanced implementations of existing schemes.\\

\subsection{Mathematical Background}\label{subsec:gsUnstructuredSplines_mathematics}
The mathematical background behind the unstructured spline module is based on the concept of mapping splines defined on local patches into a global smooth space. In other words, a globally smooth basis function $\varphi_k\in\mathbb{S}^{\text{global}}$ can be represented by a weighted sum of locally defined basis functions $\psi_l\in\mathbb{S}^{\text{local}}$ through the coefficients $A_{kl}$. Therefore, provided a vector of evaluations $\psi_l$ on a point $\xi$ located in the parametric domain of one of the patches, $\VEC{\psi}(\xi)$, the globally smooth basis at point $\xi$ is:
\begin{equation}\label{eq:gismo_map_vec}
\VEC{\varphi}(\xi) = A\VEC{\psi}(\xi).
\end{equation}
For example, for a simple one-dimensional domain with two bases separated by a $C^0$ discontinuity, a $C^1$-smooth basis over the $C^0$ interface is constructed as shown in \cref{fig:USplines}. In higher dimensions, interface smoothing, as illustrated in \cref{fig:USplines}, can be performed to construct interface basis functions, but the construction of $C^1$ functions over vertices can be challenging. Examples of such constructions include the D-Patch \cite{Reif1997,Toshniwal2017}, the Almost-$C^1$ construction \cite{Takacs2023}, the Approximate $C^1$ basis \cite{Weinmuller2021,Weinmuller2022}, the Analysis-Suitable $G^1$ construction \cite{Collin2016,Farahat2023,Farahat2023a}, Polar spline constructions \cite{Toshniwal2017a} and constructions based on subdivision surfaces \cite{Marsala2022a,Barendrecht2013}.

\begin{example}[1-Dimensional interface smoothing]
    Consider two bases of degree 2 with unique knot vectors $\Xi^{(1)} = \{0,1/8,2/8,...,7/8,1\}$ and $\Xi^{(2)} = \{1,9/8,10/8,...,15/8,2\}$. The basis functions with non-zero value on the interface are denoted by $\psi^{(1)}_{10}$ and $\psi^{(2)}_{11}$, and the basis function with non-zero derivative on the interface are denoted by, $\psi^{(1)}_9$ and $\psi^{(2)}_{12}$. When scaling $\psi^{(1)}_{10}$ and $\psi^{(2)}_1$ by a factor of $1/2$ and all the other functions with a factor of 1 in their support, $C^1$ smoothing over the interface is achieved, illustrated in \Cref{fig:USplines}. The local bases in this example have $10$ basis functions each; hence, $20$ in total. The global basis consists of $18$ functions; hence, the matrix $\MAT{A}$ is a $18\times20$ matrix. For interface basis function $\varphi_9$, the non-zero coefficients in matrix $A$ are $A_{9,9}=1$, $A_{9,10}=A_{9,11}=1/2$ and for interface basis function $\varphi_{10}$, the non-zero coefficients are $A_{10,12}=1$, $A_{10,10}=A_{9,11}=1/2$.

    \begin{figure}
        \centering
        \begin{tikzpicture}
        \begin{groupplot}
        [
        height=0.2\textheight,
        width=\fullwidth,
        ytick = {0,1},
        legend pos = outer north east,
        group style={
                        group name=my plots,
                        group size=1 by 3,
                        xlabels at=edge bottom,
                        xticklabels at=edge bottom,
                        ylabels at=edge left,
                        yticklabels at=edge left,
                        vertical sep=5pt,
                        horizontal sep=5pt},
        xlabel={$\xi$},
        ]
        \nextgroupplot[ylabel={$\psi_j$},no markers]
        \foreach \k in {0,...,7}
        {
        \pgfmathtruncatemacro{\y}{2*\k+1};
        \pgfmathtruncatemacro{\x}{2*\k};
        \addplot+[black!50,thin,solid] table[col sep = comma, header=false, y index = {\y}, x expr = \thisrowno{\x}] {Data/Refinement_illustration/init.csv};
        }

        \addplot+[black,thick,densely dashed] table[col sep = comma, header=false, y index = {17}, x expr = \thisrowno{16}] {Data/Refinement_illustration/init.csv};

        \addplot+[black,thick,densely dotted] table[col sep = comma, header=false, y index = {19}, x expr = \thisrowno{18}] {Data/Refinement_illustration/init.csv};

        \addplot+[black,thick,densely dotted] table[col sep = comma, header=false, y index = {1}, x expr = \thisrowno{0} + 1] {Data/Refinement_illustration/init.csv};

        \addplot+[black,thick,densely dashed] table[col sep = comma, header=false, y index = {3}, x expr = \thisrowno{2} + 1] {Data/Refinement_illustration/init.csv};

        \foreach \k in {2,...,9}
        {
            \pgfmathtruncatemacro{\y}{2*\k+1};
            \pgfmathtruncatemacro{\x}{2*\k};
            \addplot+[black!50,thin,solid] table[col sep = comma, header=false, y index = {\y}, x expr = \thisrowno{\x} + 1] {Data/Refinement_illustration/init.csv};
        }

        \node[left] (a) at (axis cs: 1,0.9){$\psi^{(1)}_{10}$};
        \node[right] (b) at (axis cs: 1,0.9){$\psi^{(2)}_{11}$};
        \node[] (a) at (axis cs: 0.75,0.6){$\psi^{(1)}_{9}$};
        \node[] (b) at (axis cs: 1.25,0.6){$\psi^{(2)}_{12}$};

        \nextgroupplot[ylabel={$A_{ij}\psi_j$},no markers]
        \foreach \k in {0,...,7}
        {
            \pgfmathtruncatemacro{\y}{2*\k+1};
            \pgfmathtruncatemacro{\x}{2*\k};
            \addplot+[black!50,thin,solid] table[col sep = comma, header=false, y expr = 1.0*\thisrowno{\y}, x expr = \thisrowno{\x}] {Data/Refinement_illustration/init.csv};
        }

        \foreach \k in {2,...,9}
        {
            \pgfmathtruncatemacro{\y}{2*\k+1};
            \pgfmathtruncatemacro{\x}{2*\k};
            \addplot+[black!50,thin,solid] table[col sep = comma, header=false, y expr = 1.0*\thisrowno{\y}, x expr = \thisrowno{\x} + 1] {Data/Refinement_illustration/init.csv};
        }

        \addplot+[black,thick,densely dotted] table[col sep = comma, header=false, y expr = 0.5*\thisrowno{1}, x expr = \thisrowno{0} + 1] {Data/Refinement_illustration/init.csv};
        \addplot+[black,thick,densely dashed] table[col sep = comma, header=false, y expr = \thisrowno{3}, x expr = \thisrowno{2} + 1] {Data/Refinement_illustration/init.csv};

        \addplot+[black,thick,densely dotted] table[col sep = comma, header=false, y expr = 0.5*\thisrowno{19}, x expr = \thisrowno{18}] {Data/Refinement_illustration/init.csv};
        \addplot+[black,thick,densely dashed] table[col sep = comma, header=false, y expr = \thisrowno{17}, x expr = \thisrowno{16}] {Data/Refinement_illustration/init.csv};

        \nextgroupplot[ylabel={$\varphi_i$},no markers]
        \foreach \k in {0,...,7}
        {
            \pgfmathtruncatemacro{\y}{2*\k+1};
            \pgfmathtruncatemacro{\x}{2*\k};
            \addplot+[black!50,thin,solid] table[col sep = comma, header=false, y index = {\y}, x expr = 2*\thisrowno{\x}] {Data/Refinement_illustration/BSpline.csv};
        }

        \addplot+[col2,thick,densely dashdotted] table[col sep = comma, header=false, y expr = \thisrowno{1}, x expr = \thisrowno{0} + 1] {Data/Refinement_illustration/USBspline.csv};
        \addplot+[col2,thick,densely dotted] table[col sep = comma, header=false, y expr = 0.5*\thisrowno{19}, x expr = \thisrowno{18}] {Data/Refinement_illustration/init.csv};

        \addplot+[col1,thick,densely dotted] table[col sep = comma, header=false, y expr = 0.5*\thisrowno{1}, x expr = \thisrowno{0} + 1] {Data/Refinement_illustration/init.csv};
        \addplot+[col1,thick,densely dashdotted] table[col sep = comma, header=false, y expr = \thisrowno{3}, x expr = \thisrowno{2}] {Data/Refinement_illustration/USBspline.csv};

        %

        \foreach \k in {10,...,17}
        {
            \pgfmathtruncatemacro{\y}{2*\k+1};
            \pgfmathtruncatemacro{\x}{2*\k};
            \addplot+[black!50,thin,solid] table[col sep = comma, header=false, y index = {\y}, x expr = 2*\thisrowno{\x}] {Data/Refinement_illustration/BSpline.csv};
        }

        \node[below left] (a) at (axis cs: 1,1){\textcolor{col1}{$\varphi_9$}};
        \node[below right] (b) at (axis cs: 1,1){\textcolor{col2}{$\varphi_{10}$}};
        \end{groupplot}
        \end{tikzpicture}

        \caption{The concept of interface smoothing between two bases of degree 2 with unique knot vector $\Xi^{(1)} = \{0,1/8,2/8,...,7/8,1\}$ (left) and $\Xi^{(2)} = \{1,9/8,10/8,...,15/8,2\}$ (right). The top figure plots the basis functions of the two bases, where the dashed functions have zero value but non-zero derivative on the interface. The dotted functions are non-zero on the boundary. The middle row presents scaled basis functions $A_{ij}\psi_j$. Here, all functions are scaled by a factor of 1, except the dotted functions, which are scaled by a factor of $1/2$. The bottom row presents the basis $\varphi_i=A_{ij}\psi_j$, where the sum is evaluated over the repeated index $j$. The blue and yellow functions, respectively $\varphi_9$ and $\varphi_{10}$, are constructed by taking the sum of the dashed and dotted functions in their support on the one side (resulting in the dash-dotted line) and taking the dotted line on the other side.}
        \label{fig:USplines}
    \end{figure}
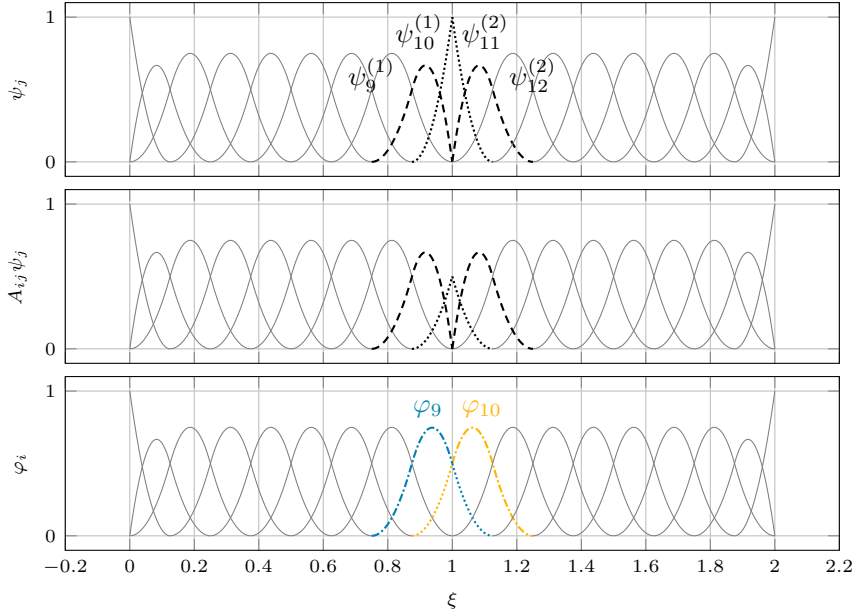
\end{example}

\subsection{Implementation}\label{subsec:gsUnstructuredSplines_implementation}
In \gismo, unstructured splines can be modelled using the \gs{gsMappedBasis} and \gs{gsMappedSpline} classes, both providing functions to evaluate a basis or a geometry that is mapped as a linear combination of underlying bases; see \cref{eq:gismo_map_vec}. The \gs{gsMappedBasis} and \gs{gsMappedSpline} classes are part of \gismo's core and are respectively both derived from \gs{gsFunctionSet}. This implies that both classes can be used as a basis or geometry, respectively, in assembly routines, for example.\\

The \gs{gsUnstructuredSplines} provide the implementation of the D-Patch \cite{Reif1997,Toshniwal2017}, the Almost-$C^1$ construction \cite{Takacs2023}, the Approximate $C^1$ basis \cite{Weinmuller2021,Weinmuller2022}, and the Analysis-Suitable $G^1$ construction \cite{Collin2016,Farahat2023,Farahat2023a} constructions as ready-to-use classes; see \cref{tab:gsUnstructuredSplines} for an overview. In other words, provided a multi-patch geometry, these classes provide a local basis as a \gs{gsMultiBasis} and a mapping matrix \gs{gsSparseMatrix}, or they directly construct a \gs{gsMappedBasis} or \gs{gsMappedSpline}. \Cref{example:gsUnstructuredSpline} provides an example of the construction of the Almost-$C^1$ basis and geometry.

\begin{example}[Unstructured spline construction]
\label{example:gsUnstructuredSpline}
Given the bi-linear multi-patch object \gs{mp} from \cref{example:gismo}, knot refinement and degree elevation steps can be performed to make the geometry suitable for the Almost-$C^1$ construction, after which this construction can be computed:

\ResetLineNumber
\begin{lstlisting}
mp.degreeIncrease(); // Increases the degree by 1 (default)
mp.uniformRefine();  // Refines the mesh by inserting 1 knot in each interval (default)

gsAlmostC1<2,real_t> almostC1(mp);
almostC1.options().setSwitch("SharpCorners",true); // Maintain C0 corners
almostC1.compute();
almostC1.matrix_into(global2local);

gsMappedBasis<2,real_t>  mbasis;
gsMappedSpline<2,real_t> mspline;
almostC1.globalBasis_into(mbasis);
almostC1.globalGeometry_into(mspline);
\end{lstlisting}
The resulting geometry, using the mesh from \cref{fig:gismo_car_mesh} in \cref{example:gismo}, is given in \cref{fig:gismo_car_patches}, with the contours of one Almost-$C^1$ basis around a valence $\nu=5$ extraordinary vertex highlighted in \cref{fig:gismo_car_function}.
\begin{figure}
\centering
\begin{subfigure}{0.59\linewidth}
\centering
\includegraphics[width=\linewidth]{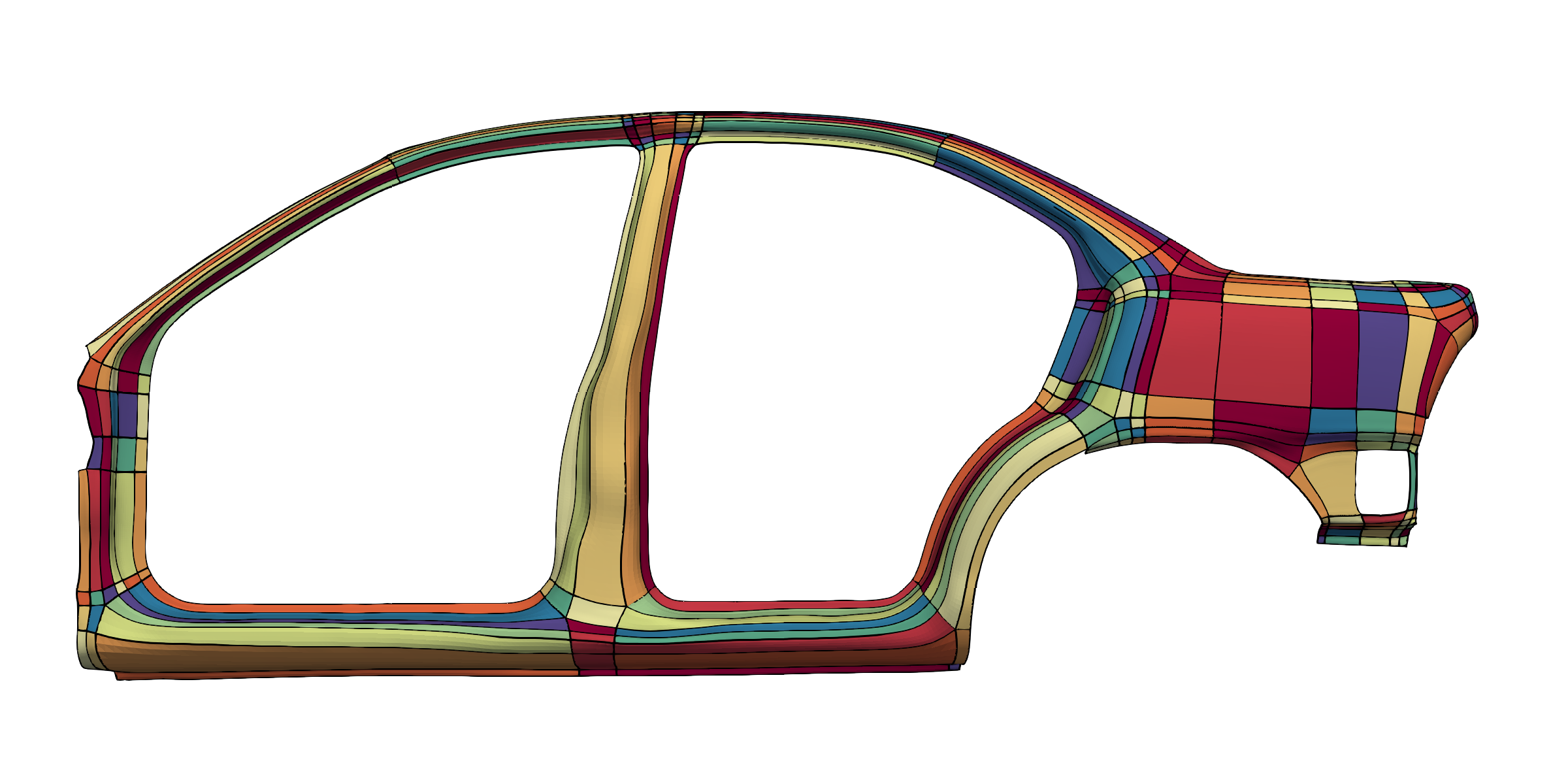}
\caption{Geometry of the side of the car.}
\label{fig:gismo_car_patches}
\end{subfigure}
\hfill
\begin{subfigure}{0.39\linewidth}
\centering
\includegraphics[width=\linewidth]{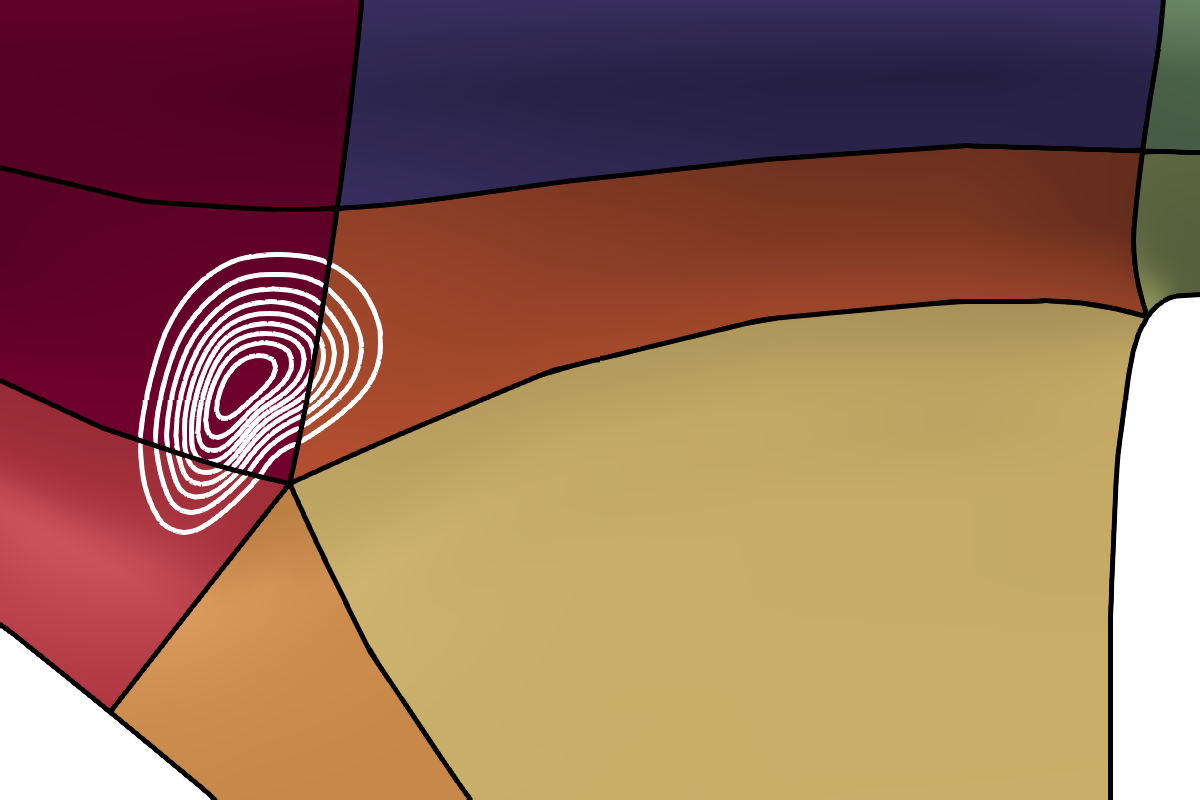}
\caption{Contour lines of an Almost-$C^1$ vertex function around an extraordinary vertex (bottom-left of \subref{fig:gismo_car_patches}).}
\label{fig:gismo_car_function}
\end{subfigure}
\caption{The multi-patch geometry constructed using the Almost-$C^1$ method \cite{Takacs2023} (\subref{fig:gismo_car_patches}) with a basis function highlighted using contour lines (\subref{fig:gismo_car_patches}). The patches are coloured randomly for the sake of referencing.}
\end{figure}

\end{example}

\section{Usage examples}\label{sec:Results}
To demonstrate the capabilities of the software presented in this paper, this section presents a number of usage examples. Most of the use cases are based on previous works, and for physical or mathematical interpretation, the reader is referred to these works. Instead, the usage examples are presented here with the aim of providing implementation background. Firstly, \cref{subsec:examples_materials} demonstrates different ways of defining materials in \gismo, and is based on \cite{Verhelst2021}. Secondly, \cref{subsec:examples_coupling} provides examples of multi-patch shell analysis using the penalty method \cite{Herrema2019} or unstructured splines \cite{Verhelst2023Coupling}. Thirdly, \cref{subsec:examples_ALM} demonstrates the use of arc-length methods and the Adaptive Parallel Arc-Length Method (APALM) \cite{Verhelst2023APALM} for quasi-static analysis. Fourth, \cref{subsec:examples_adaptivity} demonstrates adaptive meshing using goal-oriented error estimators. 

Throughout the usage examples, code snippets are provided to illustrate the use of the feature presented in the example. For the sake of brevity, the following variables are defined here and not repeated in the examples below:

\begin{table}
\caption{A list of common objects used in the examples provided in \cref{sec:Results}.}
\begin{tabular}{lll}
    \toprule
    \textbf{Symbol} & \textbf{Type} & \textbf{Description} \\
    \midrule
    \lstinline|mp|, \lstinline|mp_def| & \lstinline|gsMultiPatch| & a (deformed) multi-patch geometry. \\
    \lstinline|basis| & \lstinline|gsMultiBasis| & a basis corresponding to a multi-patch.\\
    \lstinline|bcs| & \lstinline|gsBoundaryConditions| & a set of boundary conditions.\\
    \lstinline|force| & \lstinline|gsFunctionSet| & a body force function, possibly defined on multiple patches.\\
    \lstinline|rho|, \lstinline|t| & \lstinline|gsFunctionSet| & material density and shell thickness.\\
    \lstinline|materialMatrix| & \lstinline|gsMaterialMatrixBase| & A constitutive model.\\
    \lstinline|assembler| & \lstinline|gsThinShellAssembler| & The shell assembler object.\\
    \lstinline|Residual| & \lstinline|Residual_t| & The residual operator, using \cref{eq:gismo_gsKLShell_Residual}.\\
    \lstinline|ALResidual| & \lstinline|ALResidual_t| & The residual operator for arc-length methods.\\
    \lstinline|Jacobian| & \lstinline|Jacobian_t| & The Jacobian matrix, using \cref{eq:gismo_gsKLShell_Jacobian}.\\
    \bottomrule
\end{tabular}
\end{table}

For details regarding the installation of \gismo, the \gs{gsKLShell}, \gs{gsStructuralAnalysis} and \gs{gsUnstructuredSplines} modules, the reader is referred to \cref{sec:installation}. In addition \cref{sec:reproduction} provides detailed instructions to reproduce every example provided in the following sub-sections.


\subsection{Material modeling}\label{subsec:examples_materials}
One of the features of the \gs{gsKLShell} module is its versatility with respect to constitutive modeling. By decoupling the constitutive models \gs{gsMaterialMatrixBase} from the kinematics and assembly \gs{gsThinShellAssembler}, see \cref{tab:gsKLShell}, different material models can be used with the generic shell assembler with a minimal interface. In this section, an example of material selection is provided, after which the novel hyperelastic tension field theory-based material model is demonstrated.

\subsubsection{Material Specification}
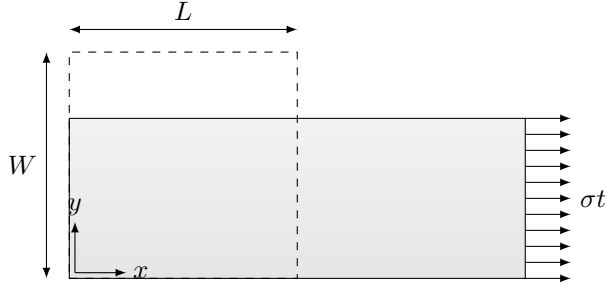
\begin{figure}
    \centering
    \begin{tikzpicture}[scale=0.75]
        \def\B{4}
        \def\L{4}
        \def\Bdef{2.828427}
        \def\Ldef{8}

        \filldraw[opacity=1, bottom color=black!10, top color=black!05, draw=black] (0,0) -- (\Ldef,0) -- (\Ldef,\Bdef) -- (0,\Bdef) --cycle;

        \draw[draw=black, dashed] (0,0) -- (\L,0) -- (\L,\B) -- (0,\B) --cycle;

        \draw[latex-latex] (0,\B+0.1*\B) -- node[midway,above] {$L$} (\L,\B+0.1*\B);
        \draw[latex-latex] (-0.1*\B,0) -- node[midway,left] {$W$} (-0.1*\B,\B);

        \draw[-latex] (0.1,0.1) -- node[right,inner sep=12pt] {$x$}(1,0.1);
        \draw[-latex] (0.1,0.1) -- node[above,inner sep=12pt] {$y$}(0.1,1);

        \foreach \k in {0,1,2,3,4,5,6,7,8,9,10}
        {
            \pgfmathsetmacro\y{\k*\Bdef / 10}
            \draw[-latex] (\Ldef,\y) -- (\Ldef + \Ldef/10,\y);
        }
        \node[right] at (\Ldef + \Ldef/10,\Bdef / 2) {$\sigma t$};

    \end{tikzpicture}
    \caption{Geometry for the uniaxial tension case. The filled geometry represents the deformed configuration, and the dashed line indicates the undeformed geometry. The bottom side of the undeformed sheet is fixed in $y$-direction and the left side of the sheet is fixed in $x$-direction. The applied load is $\sigma t$, where $\sigma$ is the actual Cauchy stress and $t$ is the thickness of the sheet.}
    \label{fig:axialGeometry}
\end{figure}

Consider the uniaxial tension case depicted in \cref{fig:axialGeometry} for a membrane with geometric dimensions $L\times W \times t = 1\times1\times0.001\:[\text{m}^3]$ and material parameters defined via Lam\'e parameter $\mu=E/(2(1+\nu))=1.5\cdot10^6\:[\text{N}/\text{m}^2]$ using a Poisson's ratio $\nu=0.45\:[-]$ such that the Young's modulus $E$ can be inferred. For a Mooney-Rivlin material, $\mu=\mu_1+\mu_2$, with $\mu_1/\mu_2=7$. To investigate the uniaxial tension behaviour of a material using different material laws, material model specification of different models is necessary. In the following, a demonstration of defining different material models in different ways using the \gs{gsKLShell} module is given. The uniaxial tension example is adopted from \cite{Verhelst2021}, and the model parameters are given in \cref{fig:axialGeometry}. Here, the example material definitions are given for compressible Neo-Hookean and Mooney-Rivlin material models.\\

To model the material laws for the uniaxial tension case, the material classes from \cref{tab:gsKLShell} need to be defined. As can be seen in \cref{example:gsKLShell}, a pointer to a material matrix, deriving from \gs{gsMaterialMatrixBase}, needs to be passed to the \gs{gsThinShellAssembler}. This can be done by explicitly constructing a material matrix and passing it into the assembler, as done in \cref{example:gsKLShell}. Alternatively, one can use the \gs{getMaterialMatrix} helper function, for example:

\begin{lstlisting}
// Provided mp, t, rho, ratio
// Define a MR material model with analytically implemented S and C tensors.
gsOptionList options;
options.addInt("Material","3 = Mooney-Rivlin",3);
options.addInt("Implementation","1 = Analytical",1);
gsConstantFunction<real_t> E(1.5e6*2*(1+0.45),3); // Young's modulus: E=mu*2*(1+nu)
gsConstantFunction<real_t> nu(0.45,3);            // Poisson's ratio
gsConstantFunction<real_t> Ratio(7.0,3);          // The material lives in the 3D domain
std::vector<gsFunctionSet<real_t> *> parameters{&E,&nu,&Ratio};
gsMaterialMatrixBase<real_t> * materialMatrix = getMaterialMatrix(mp,t,parameters,rho,options);
\end{lstlisting}
Here, the \gs{parameters} object contains a list of functions defining the material properties. These functions, as well as the thickness \gs{t} and the density \gs{rho}, are passed as \gs{gsFunctionSet} pointers, meaning that they can be defined using any type of function definition in \gismo. This allows us to define material properties as splines, allowing optimisation of material properties in an isogeometric setting. In the example above, the material parameters are simply passed as \gs{gsConstantFunction} objects, with the function value as the first argument and the geometric domain dimension as the second argument.\\

\begin{figure}[tb!]
    \centering
    \begin{tikzpicture}
        \footnotesize
        \begin{groupplot}
            [
            group style={
                group size=2 by 3,
                ylabels at=edge left,
                xlabels at=edge bottom,
                x descriptions at=edge bottom,
                vertical sep=0.05\linewidth,
                horizontal sep=0.05\textheight,
            },
            legend pos = north east,
            width=\halfwidth-0.05\linewidth,
            height=0.25\textheight-0.05\textheight,
            ]
            \nextgroupplot
            [
            title={Compressible},
            xlabel={\vphantom{x}},
            ylabel=Through-thickness stretch $\lambda_3$,
            xmin = 0,
            xmax = 12.5,
            ymin = 0.2,
            ymax = 1,
            ]

            \addplot+[style=style0,only marks] table[header=true,x index = {0},y index = {1}, col sep = comma]{Data/UAT/NH_compressible.csv};
            \addplot+[style=style1,only marks] table[header=true,x index = {0},y index = {1}, col sep = comma]{Data/UAT/MR_compressible.csv};

            \addplot+[style=style0,mark=none] table[header=true,x index = {0},y expr = sqrt(\thisrowno{1}/x), col sep = comma]{Data/UAT/Compressible_analytical.csv};
            \addplot+[style=style1,mark=none] table[header=true,x index = {0},y expr = sqrt(\thisrowno{2}/x), col sep = comma]{Data/UAT/Compressible_analytical.csv};

            \nextgroupplot
            [
            title={Incompressible},
            xlabel={\vphantom{x}},
            xmin = 0,
            xmax = 12.5,
            ymin = 0.2,
            ymax = 1,
            ]
            \addlegendimage{line legend,black,mark=none}
            \addlegendentry{$1/\sqrt{\lambda}$}

            \addplot+[style=style0,only marks,fill opacity=0] table[header=true,x index = {0},y index = {1}, col sep = comma]{Data/UAT/NH_incompressible.csv};
            \addplot+[style=style1,only marks,fill opacity=0] table[header=true,x index = {0},y index = {1}, col sep = comma]{Data/UAT/MR_incompressible.csv};
            \addplot+[style=style2,only marks,fill opacity=0] table[header=true,x index = {0},y index = {1}, col sep = comma]{Data/UAT/OG_incompressible.csv};

            \addplot[thin,mark=none,black,domain=0.01:12.5,samples = 100] {1/sqrt(x)};

            \nextgroupplot
            [
            xlabel={\vphantom{x}},
            ylabel=Axial Cauchy-stress $\sigma$,
            xmin = 0,
            xmax = 12.5,
            ]
            \addplot+[style=style0,only marks] table[header=true,x index = {0},y index = {2}, col sep = comma]{Data/UAT/NH_compressible.csv};
            \addplot+[style=style1,only marks] table[header=true,x index = {0},y index = {2}, col sep = comma]{Data/UAT/MR_compressible.csv};
            \addplot+[style=style2,only marks] table[header=true,x index = {0},y index = {2}, col sep = comma]{Data/UAT/OG_compressible.csv};

            \addplot+[style=style0,mark=none] table[header=true,x index = {0},y index = {4}, col sep = comma]{Data/UAT/Compressible_analytical.csv};
            \addplot+[style=style1,mark=none] table[header=true,x index = {0},y index = {5}, col sep = comma]{Data/UAT/Compressible_analytical.csv};
            \addplot+[style=style2,mark=none] table[header=true,x index = {0},y index = {6}, col sep = comma]{Data/UAT/Compressible_analytical.csv};

            \nextgroupplot
            [
            xlabel=Stretch $\lambda$,
            xmin = 0,
            xmax = 12.5,
            ]
            \addplot+[style=style0,only marks,fill opacity=0] table[header=true,x index = {0},y index = {2}, col sep = comma]{Data/UAT/NH_incompressible.csv};
            \addplot+[style=style1,only marks,fill opacity=0] table[header=true,x index = {0},y index = {2}, col sep = comma]{Data/UAT/MR_incompressible.csv};
            \addplot+[style=style2,only marks,fill opacity=0] table[header=true,x index = {0},y index = {2}, col sep = comma]{Data/UAT/OG_incompressible.csv};

            \addplot+[style=style0,mark=none] table[header=true,x index = {0},y index = {1}, col sep = comma]{Data/UAT/Incompressible_analytical.csv};
            \addplot+[style=style1,mark=none] table[header=true,x index = {0},y index = {2}, col sep = comma]{Data/UAT/Incompressible_analytical.csv};
            \addplot+[style=style2,mark=none] table[header=true,x index = {0},y index = {3}, col sep = comma]{Data/UAT/Incompressible_analytical.csv};

            \nextgroupplot
            [
            xlabel=Stretch $\lambda$,
            ylabel=Jacobian determinant $J$,
            xmin = 0,
            xmax = 12.5,
            ]
            \addplot+[style=style0,only marks] table[header=true,x index = {0},y index = {3}, col sep = comma]{Data/UAT/NH_compressible.csv};
            \addplot+[style=style1,only marks] table[header=true,x index = {0},y index = {3}, col sep = comma]{Data/UAT/MR_compressible.csv};
            \addplot+[style=style2,only marks] table[header=true,x index = {0},y index = {3}, col sep = comma]{Data/UAT/OG_compressible.csv};

            \addplot+[style=style0,mark=none] table[header=true,x index = {0},y index = {1}, col sep = comma]{Data/UAT/Compressible_analytical.csv};
            \addplot+[style=style1,mark=none] table[header=true,x index = {0},y index = {2}, col sep = comma]{Data/UAT/Compressible_analytical.csv};
            \addplot+[style=style2,mark=none] table[header=true,x index = {0},y index = {3}, col sep = comma]{Data/UAT/Compressible_analytical.csv};

            \nextgroupplot
            [%
            hide axis,
            xmin=0,
            xmax=0.1,
            ymin=0,
            ymax=0.1,
            transpose legend,
            legend columns=4,
            legend style={
                at={(0.5,0.5)},
                anchor=center}
            ]

            \addlegendimage{legend image with text= }\addlegendentry{};
            \addlegendimage{legend image with text=NH}\addlegendentry{};
            \addlegendimage{legend image with text=MR}\addlegendentry{};
            \addlegendimage{legend image with text=OG}\addlegendentry{};

            \addlegendimage{legend image with text=C}\addlegendentry{};
            \addlegendimage{style=style0,only marks};\addlegendentry{};
            \addlegendimage{style=style1,only marks};\addlegendentry{};
            \addlegendimage{style=style2,only marks};\addlegendentry{};

            \addlegendimage{legend image with text=I}\addlegendentry{};
            \addlegendimage{style=style0,fill opacity=0,only marks};\addlegendentry{};
            \addlegendimage{style=style1,fill opacity=0,only marks};\addlegendentry{};
            \addlegendimage{style=style2,fill opacity=0,only marks};\addlegendentry{};

            \addlegendimage{legend image with text=A}\addlegendentry{};
            \addlegendimage{style=style0,no markers};\addlegendentry{};
            \addlegendimage{style=style1,no markers};\addlegendentry{};
            \addlegendimage{style=style2,no markers};\addlegendentry{};
        \end{groupplot}
    \end{tikzpicture}
    \caption{Results for uniaxial tension for compressible (C, left column) and incompressible materials (I, right column), where the first row presents the thickness decrease $\lambda_3$, the second row the axial Cauchy stress or true axial stress $\sigma$, and the last row the Jacobian determinant $J$ for compressible materials, all against the stretch $\lambda$. The material models that are used are the Neo-Hookean (NH), the Mooney-Rivlin (MR), and the Ogden (OG) material models, and comparisons are made to analytical (A) solutions from \cite[ex. 1]{Holzapfel2000}.}
    \label{fig:hyperelasticity_uniaxial}
\end{figure}
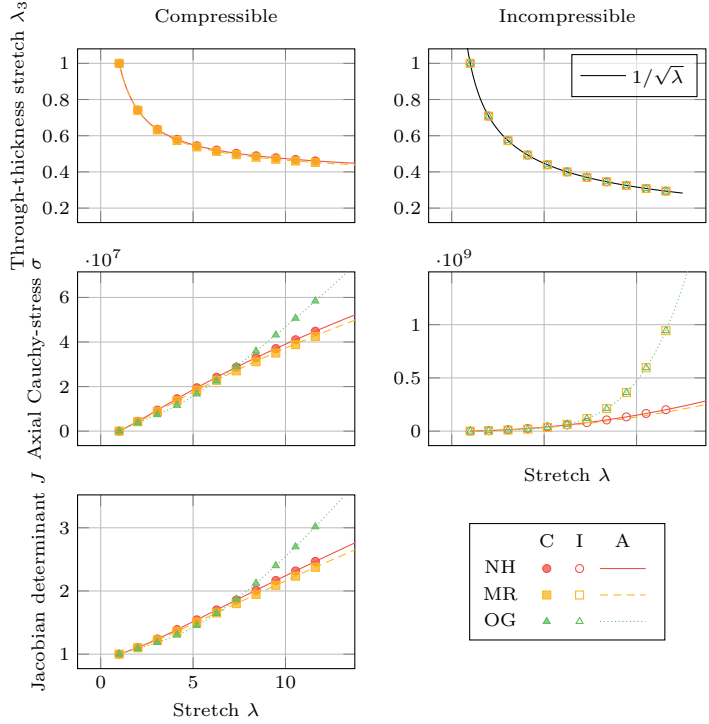

Alternative to the \gs{getMaterialMatrix} routine, one can define a material matrix using \gismo's XML format, providing an XML item, for example, for a compressible Neo-Hookean material:

\ResetLineNumber
\begin{lstlisting}[style=XML]
<MaterialMatrix type="IncompressibleNH3" id="0" TFT="false">
    <Thickness>
        <Function type="FunctionExpr" dim="3" index="0">0.0001</Function>
    </Thickness>
    <Density>
        <Function type="FunctionExpr" dim="3" index="0">0</Function>
    </Density>
    <Parameters>
        <!-- Youngs Modulus -->
        <Function type="FunctionExpr" dim="3" index="0">1.5e6*2*(1+0.45)</Function>
        <-- Poisson Ratio     -->
        <Function type="FunctionExpr" dim="3" index="1">0.45</Function>
    </Parameters>
</MaterialMatrix>
\end{lstlisting}
Using the material matrices defined above, either defined using \gs{getMaterialMatrix}, constructed directly, or read via XML, the \gs{gsThinShellAssembler} can be constructed as in \cref{example:gsKLShell}, and the \gs{Residual_t} and \gs{Jacobian_t} operators can be defined to perform non-linear arc-length analysis with an increasing load. The reader is referred to \cref{example:gsStructuralAnalysis} for the definition of the required solver. Finally, the results presented in \cref{fig:hyperelasticity_uniaxial} are obtained by using the material specifications from above.

\subsubsection{Tension Field Theory}
Consider the cylinder depicted in \cref{fig:cylinder_setup}, with a radius of $R=250\:[\text{mm}]$, length $H=1.0\:[\text{mm}]$, thickness $t=0.05\:[\text{mm}]$ and material parameters $E=1.0\:[\text{GPa}]$ and $\nu=0.5\:[-]$ subject to a boundary rotation $\theta=\pi/2\:[\text{rad}]$ and translation $u_x=H$. When modeling this example with a shell model, the deformation pattern is composed of diagonal wrinkles over the length of the cylinder; see \cref{fig:cylinder}. Since modeling the membrane wrinkles could require a large number of elements, an alternative way of modeling this benchmark is to model only the mid-plane of the wrinkled geometry using a modified membrane wrinkling model. Here, the use of the novel hyperelastic tension field theory-based model for general hyperelastic material models \cite{Verhelst2025} is demonstrated, within the \gs{gsKLShell} module.\\

\begin{figure}[tb!]
    \centering
    \begin{minipage}{\twothirdwidth}
        \centering
        \begin{tikzpicture}[]
            \def\R{2}
            \begin{scope}[rotate=00]
                \begin{axis}[
                    axis equal image,
                    hide axis,
                    width=\linewidth,
                    view = {20}{10},
                    scale = 2,
                    xmin = 0.0,
                    xmax = 8.0,
                    ymin = -\R-0.5,
                    ymax = \R+0.5,
                    zmin = -\R-0.5,
                    zmax = \R+0.5,
                    line cap=round
                    ]\

                    \addplot3[domain=-0.5*pi:0.5*pi, samples=100, samples y=0, no marks, smooth, thick,gray](
                    {0},
                    {\R*cos(deg(\x)},
                    {\R*sin(deg(\x))}
                    );

                    \addplot3[
                    surf,
                    samples = 20,
                    fill opacity=0.5,
                    samples y = 2,
                    domain = 0.0*pi:2*pi,
                    domain y = 0:3,
                    draw=none,
                    z buffer = sort,
                    no marks,
                    mesh/interior colormap={blueblack}{color=(black) color=(white)},
                    colormap ={blueblack}{color=(black) color=(white)},
                    shader=interp,
                    opacity=0.5,
                    point meta=x+3*y*y-0.25*z,
                    ](
                    {2*\y},
                    {(\R)*cos(deg(\x)},
                    {(\R)*sin(deg(\x))}
                    );

                    \addplot3[domain=-0.5*pi:0.5*pi, samples=100, samples y=0, no marks, smooth, thick,black](
                    {0},
                    {-\R*cos(deg(\x)},
                    {-\R*sin(deg(\x))}
                    );

                    \addplot3[domain=0:2*pi, samples=100, samples y=0, no marks, smooth, thick,black](
                    {6},
                    {\R*cos(deg(\x)},
                    {\R*sin(deg(\x))}
                    );

                    \addplot3[domain=0.0*pi:0.5*pi, samples=100, samples y=0, no marks, smooth, thick,black,latex-](
                    {7.5},
                    {\R*cos(deg(-\x)},
                    {\R*sin(deg(-\x))}
                    );

                    \foreach \t in {0.0,0.1,...,2} {

                        \edef\x{\R*cos(deg(\t*3.1415} 
                        \edef\y{\R*sin(deg(\t*3.1415} 
                        \edef
                        \temp{
                            \noexpand
                            \draw [-latex] (axis cs:6,\x,\y) -- (axis cs:7,\x,\y);
                        }
                        \temp
                    }

                    \addplot3[domain=0.0*pi:2*pi, samples=100, samples y=0, no marks, smooth, thick](
                    {7},
                    {\R*cos(deg(\x)},
                    {\R*sin(deg(\x))}
                    );



                    \node[right] at (axis cs: 7,0.707*\R,0.707*\R) {$u_x$};
                    \node[right] at (axis cs: 7.5,0.707*\R,-0.707*\R) {$\varTheta$};

                    \draw[thick] (axis cs: 0,0,\R) -- (axis cs: 6,0,\R);
                    \draw[thick] (axis cs: 0,0,-\R) -- (axis cs: 6,0,-\R);

                    \draw[latex-latex] (axis cs: 0,0,0) -- (axis cs: 0,-0.707*\R,-0.707*\R) node[midway, right]{$R$};

                    \draw[latex-latex] (axis cs: 0,0,-2.2) -- (axis cs: 6,0,-2.2) node[midway,below]{$H$};
                    \node at (axis cs: 0,0,0) {$\times$};

                    \draw[-latex] (axis cs: 2,0,0) -- (axis cs: 2,0,1.0) node[above]{$z$};
                    \draw[-latex] (axis cs: 2,0,0) -- (axis cs: 2,2.0,0) node[above right]{$y$};
                    \draw[-latex] (axis cs: 2,0,0) -- (axis cs: 3,0,0) node[right]{$x$};

                \end{axis}
            \end{scope}
        \end{tikzpicture}
    \end{minipage}
    \hfill
    \begin{minipage}{0.3\linewidth}
        \centering
        \footnotesize
        \begin{tabular}{llr}
            \toprule
            \multicolumn{3}{c}{Geometry}\\
            \midrule
            $R$ & 250  & $[\text{mm}]$\\
            $L$ & 1.0 & $[\text{m}]$\\
            $t$ & 0.05 & $[\text{mm}]$\\
            \midrule
            \multicolumn{3}{c}{Material}\\
            \multicolumn{3}{c}{\textit{Incompressible Neo-Hookean}}\\
            \midrule
            $E$ & 1.0 & $[\text{GPa}]$\\
            $\nu$ & 0.5 & $[-]$\\
            \midrule
            \multicolumn{3}{c}{Boundary Conditions}\\
            \midrule
            $\varTheta$ & $\pi/2$ & $[\text{rad}]$\\
            $u_x$ & $1.0$ & $[\text{m}]$\\
            \bottomrule
        \end{tabular}
    \end{minipage}

    \caption{Problem definition for an cylinder with inner radius $R$ and length $L$ subject to an elongation $u_x$ and a rotation $\varTheta$ on the right boundary $\Gamma_r$ while being fixed on the left boundary $\Gamma_l$. The cylinder has a Neo-Hookean material model with the parameters provided in the table on the right.}
    \label{fig:cylinder_setup}
\end{figure}

Tension field theory-based membrane models are enabled by the \gs{gsMaterialMatrixTFT} class. This class uses tension field -- fields based on principal strains and stresses -- to modify the stress and material tensors $\TEN{S}$ and $\Ccten$ according to the modification scheme from \cite{Nakashino2005,Nakashino2020} for linear elastic materials and from \cite{Verhelst2025} for hyperelastic materials. For each quadrature point, the stress and material tensors are determined by the tension field: they are set to zero (or a small number) in the slack state, they are modified according to the schemes in \cite{Nakashino2005,Nakashino2020,Verhelst2025} in the wrinkled state, and they are unmodified in the taut state.\\

Since the \gs{gsMaterialMatrixTFT} class inherits from \gs{gsMaterialMatrixBase}, it can be used in the \gs{gsThinShellAssembler} like any other material model in \gs{gsKLShell}. The class \gs{gsMaterialMatrixTFT} can be initialised by setting the \gs{TFT} flag in the XML specification above to \gs{true}, or by constructing the object as:

\begin{lstlisting}
// Provided materialMatrix<d,real_t>
gsMaterialMatrixTFT<d,real_t,true> materialMatrixTFT(&materialMatrix);
materialMatrixTFT.options().setInt();
\end{lstlisting}
Here, the template parameter \gs{d} specifies the geometric dimension, either 2 or 3, and the last template parameter specifies whether a linear modification (\gs{true}) \cite{Nakashino2005,Nakashino2020} or a non-linear modification (\gs{false}) \cite{Verhelst2025} should be used. The results in \cref{fig:cylinder} are obtained by a two-stage static solver composed of a dynamic relaxation method followed by a Newton-Raphson method, defined as follows:
\begin{lstlisting}
// Provided F, M, K, Residual, Jacobian, alpha
gsStaticDR<real_t> DR(M,F,Residual);
DR.options().setReal("alpha",alpha);
DR.options().setReal("tolF",1e-4);
DR.initialize();
gsStaticNewton<real_t> NR(K,F,Jacobian,Residual);
NR.options().setReal("tolF",1e-6);
NR.initialize();
gsStaticComposite<real_t> solver({&DR,&NR});
\end{lstlisting}
Here, the parameter \gs{alpha} is the tuning parameter for the dynamic relaxation method. Furthermore, the tolerances specified are the relative residual tolerances.

Inside the displacement-step loop, this composite solver is accordinly called after updating the boundary conditions:
\begin{lstlisting}
dx.set_t(D); dy.set_t(D); dz.set_t(D);
assembler->updateBCs(BCs);
StaticSolver.setDisplacement(U);
StaticSolver.solve();
\end{lstlisting}
In the snippet above, \gs{dx}, \gs{dy} and \gs{dz} are displacement functions defined on the boundaries of the cylinder for $x$, $y$ and $z$ directions, respectively. These functions are defined in terms of a variable $t$ representing the rotation $\theta$, which is incrementally set by the variable \gs{D}. After updating the the displacement functions, the previous displacement vector \gs{U} is passed to the solver by means of an initial guess, after which the system is solved.\\

The final results are provided in \cref{fig:cylinder}, and the reader is referred to \cite{Verhelst2025} for a detailed discussion and numerical validation of the obtained results.

\begin{figure}
    \centering
    \begin{subfigure}[t]{0.45\linewidth}
        \centering
        \includegraphics[width=\linewidth]{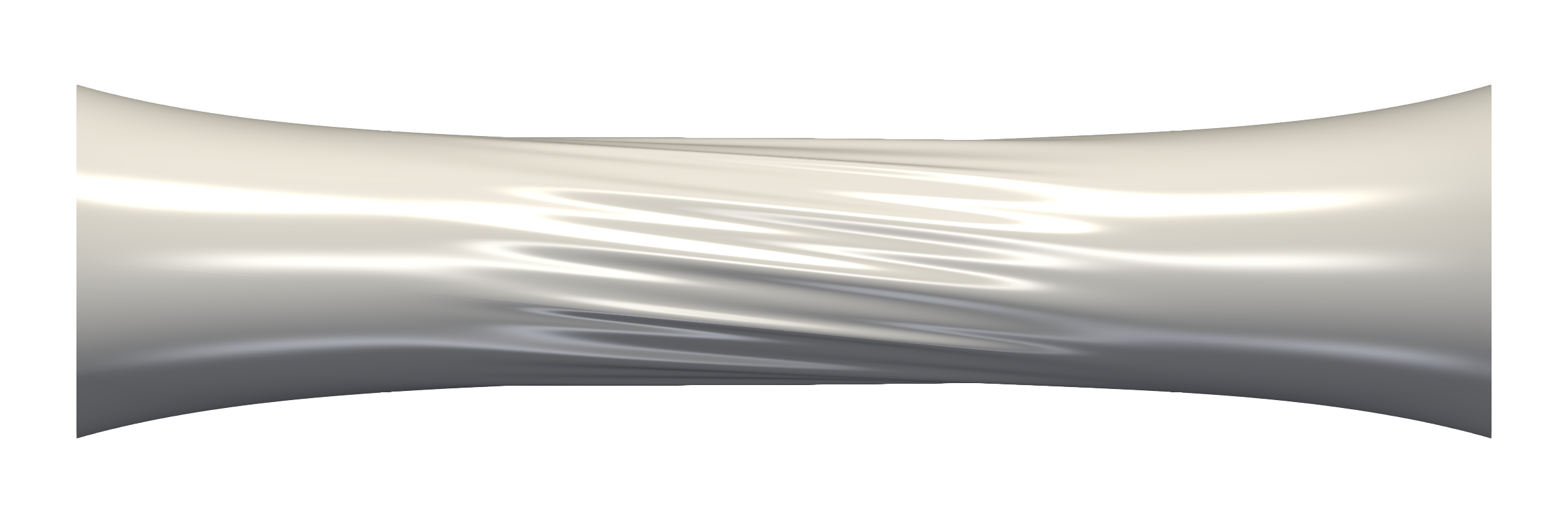}
        \includegraphics[width=\linewidth]{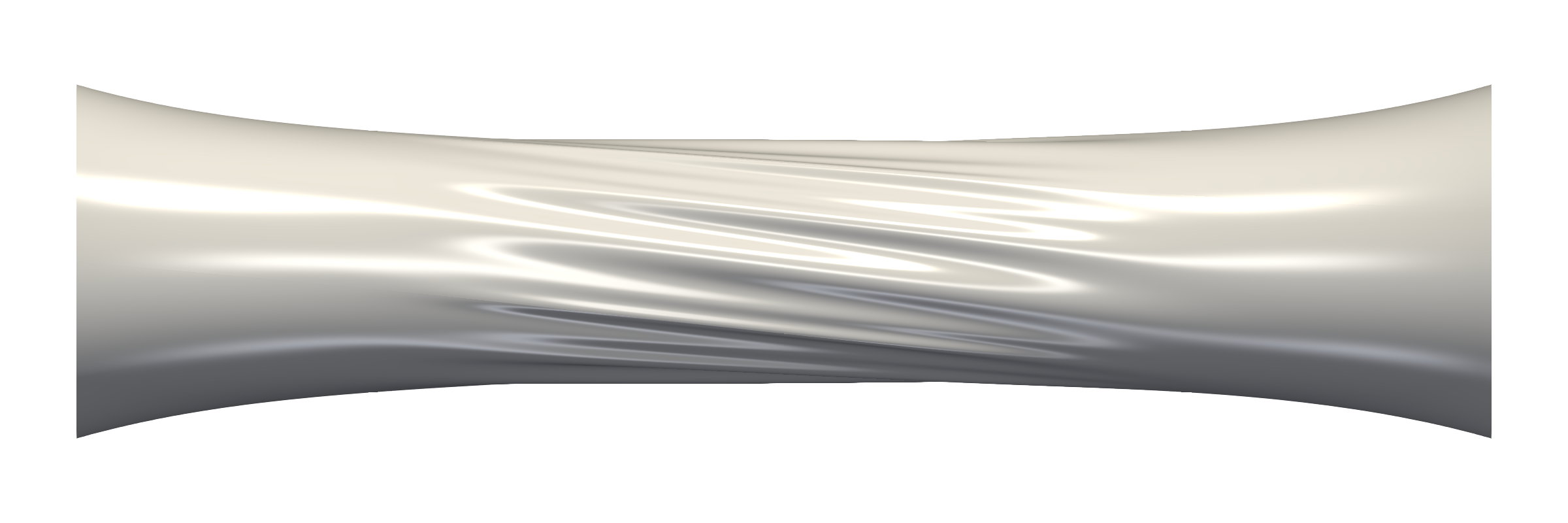}
        \includegraphics[width=\linewidth]{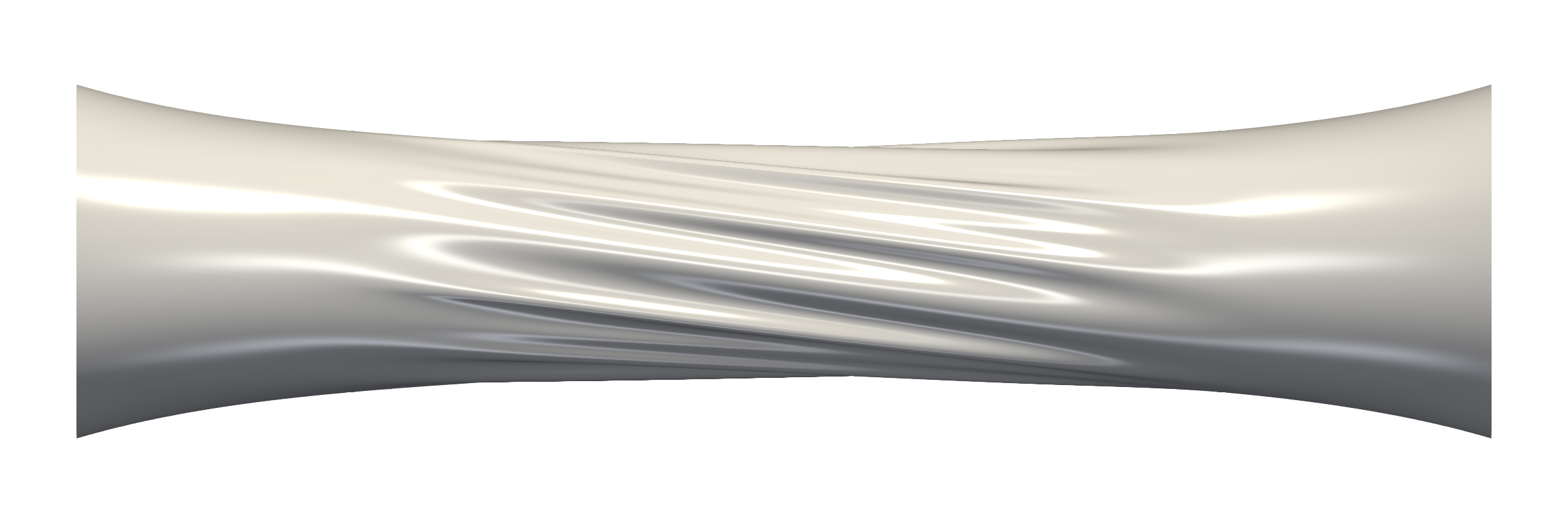}     \includegraphics[width=\linewidth]{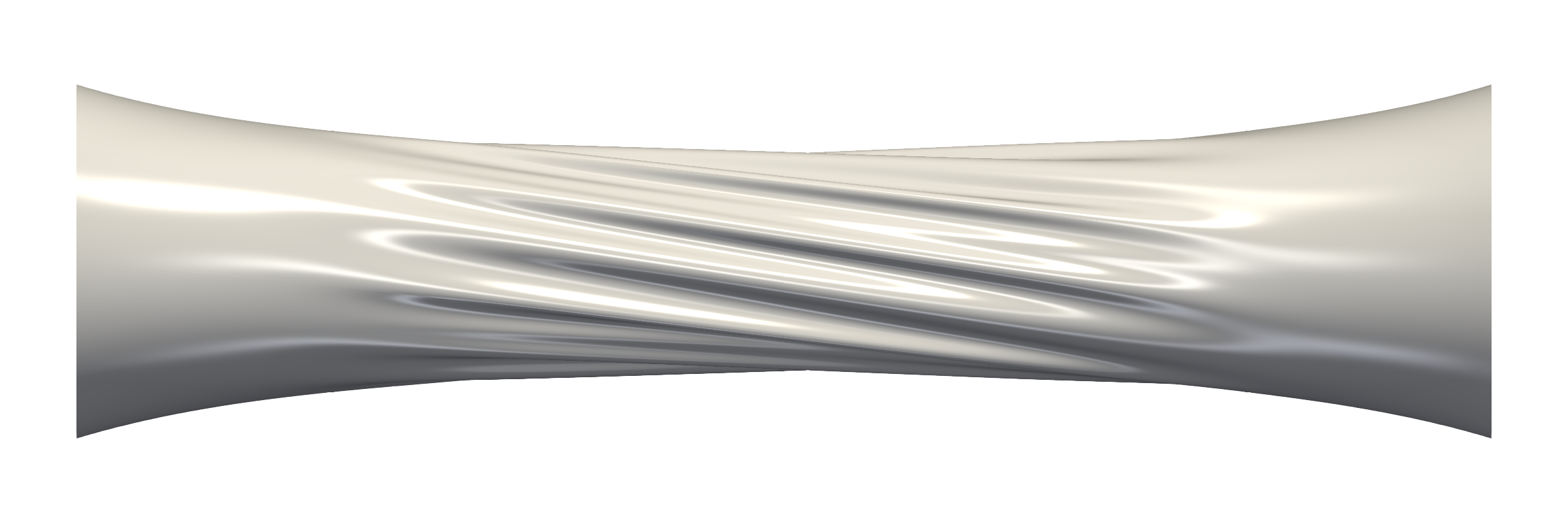}
        \caption{Side view of the deformed cylinder using the Kirchhoff--Love shell model with a rotation of $\{0.35,0.40,0.45,0.50\}\:[\pi\:\text{rad}]$ from top to bottom.}
        \label{fig:cylinder_side}
    \end{subfigure}
    \hfill
    \begin{subfigure}[t]{0.45\linewidth}
        \centering
        \includegraphics[width=\linewidth]{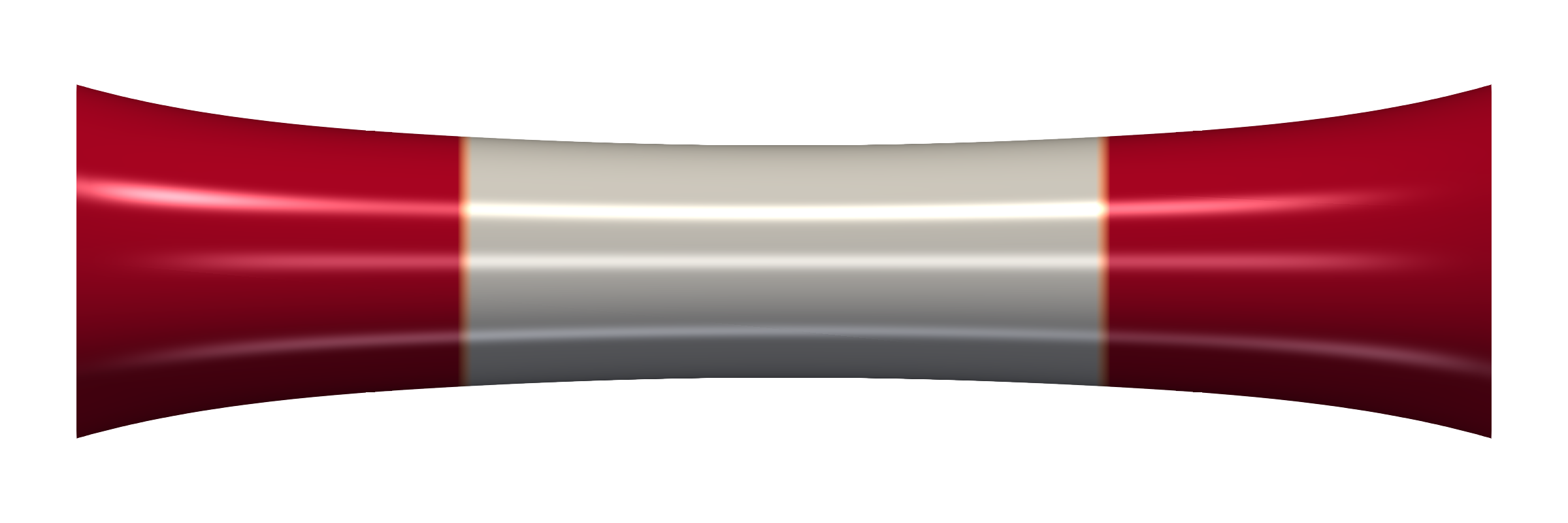}
        \includegraphics[width=\linewidth]{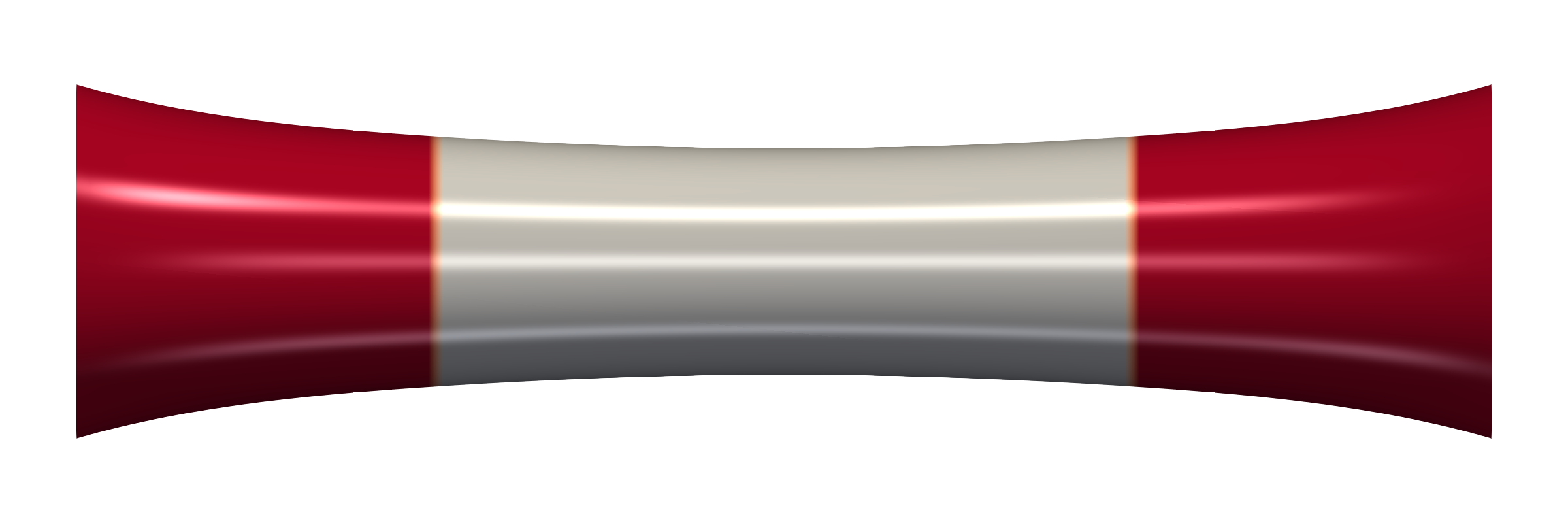}
        \includegraphics[width=\linewidth]{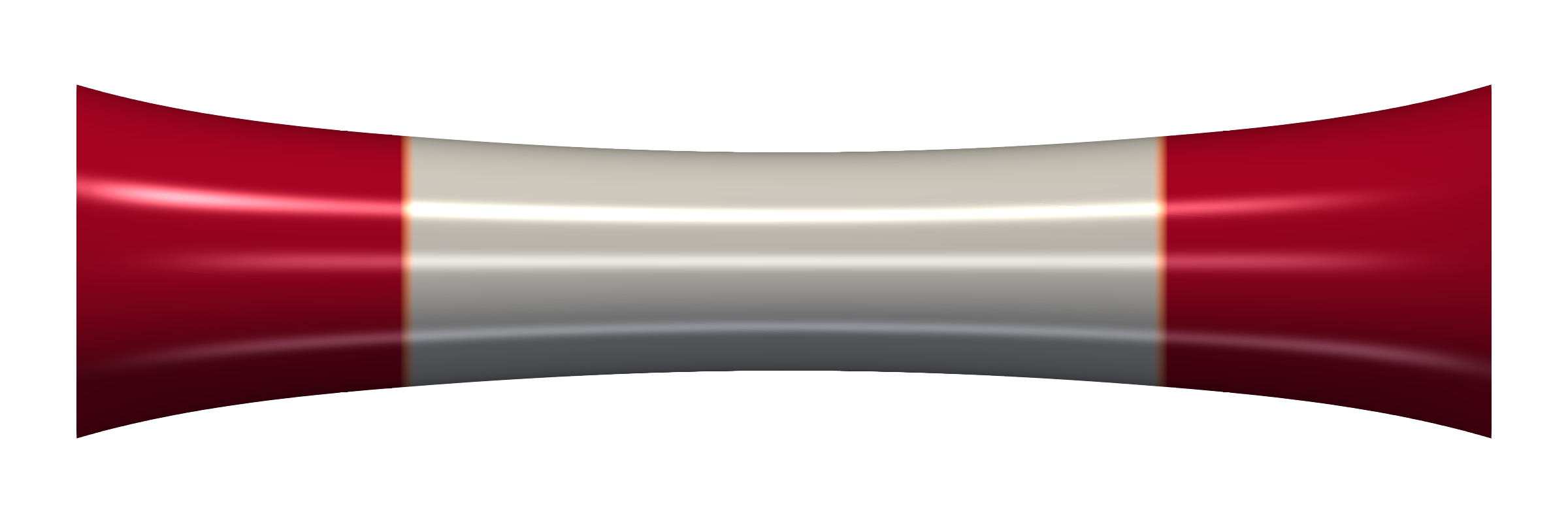}
        \includegraphics[width=\linewidth]{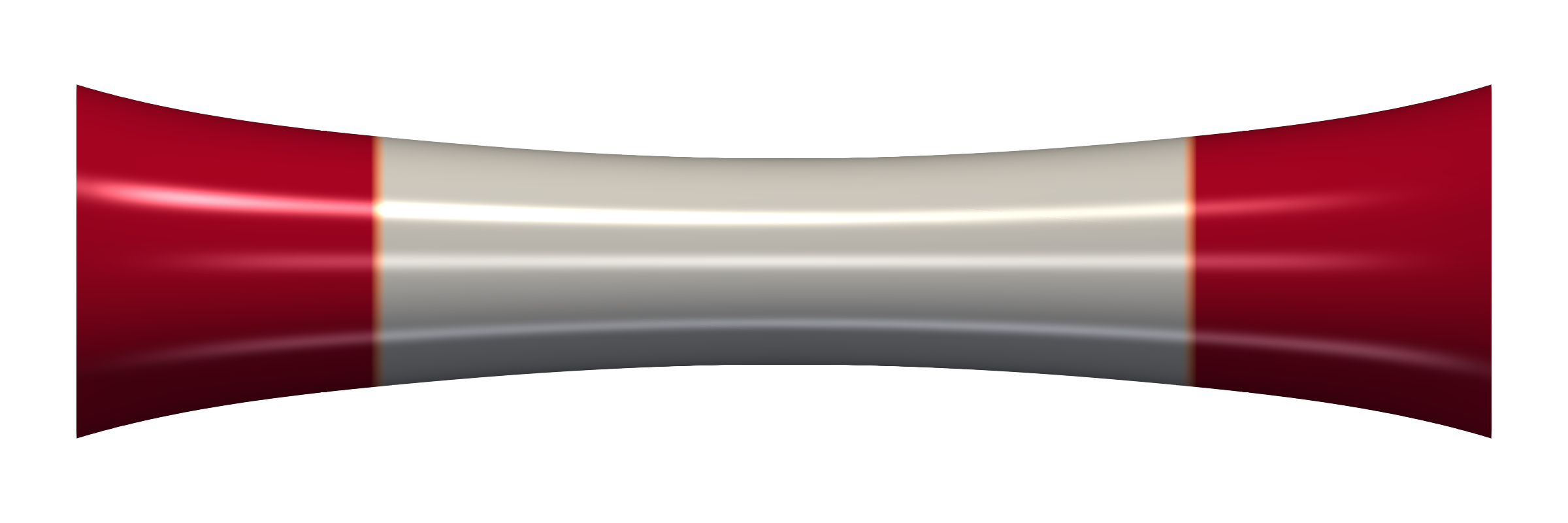}
        \caption{Side view of the deformed cylinder using the tension field theory membrane model with a rotation of $\{0.35,0.40,0.45,0.50\}\:[\pi\:\text{rad}]$ from top to bottom. The red region denotes a taut, and the gray region denotes a wrinkled region.}
        \label{fig:cylinder_side_tensionfield}
    \end{subfigure}
\caption{(Caption on next page)}
\end{figure}

\begin{figure}
\ContinuedFloat
    \centering
    \begin{subfigure}[t]{0.45\linewidth}
        \centering
        \includegraphics[width=\linewidth]{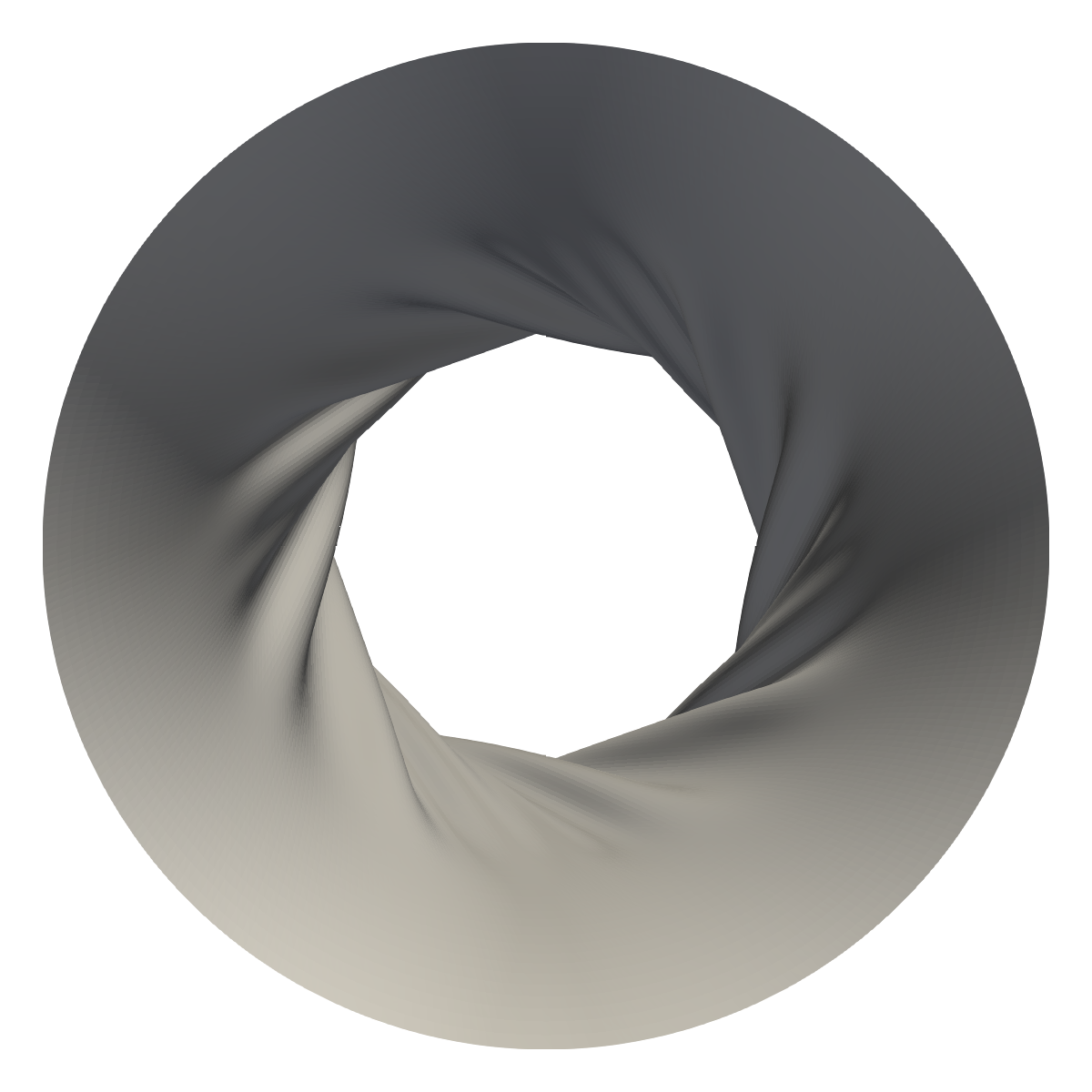}
        \caption{Top view of the deformed cylinder using the Kirchhoff--Love shell model.}
        \label{fig:cylinder_top}
    \end{subfigure}
    \hfill
    \begin{subfigure}[t]{0.45\linewidth}
        \centering
        \includegraphics[width=\linewidth]{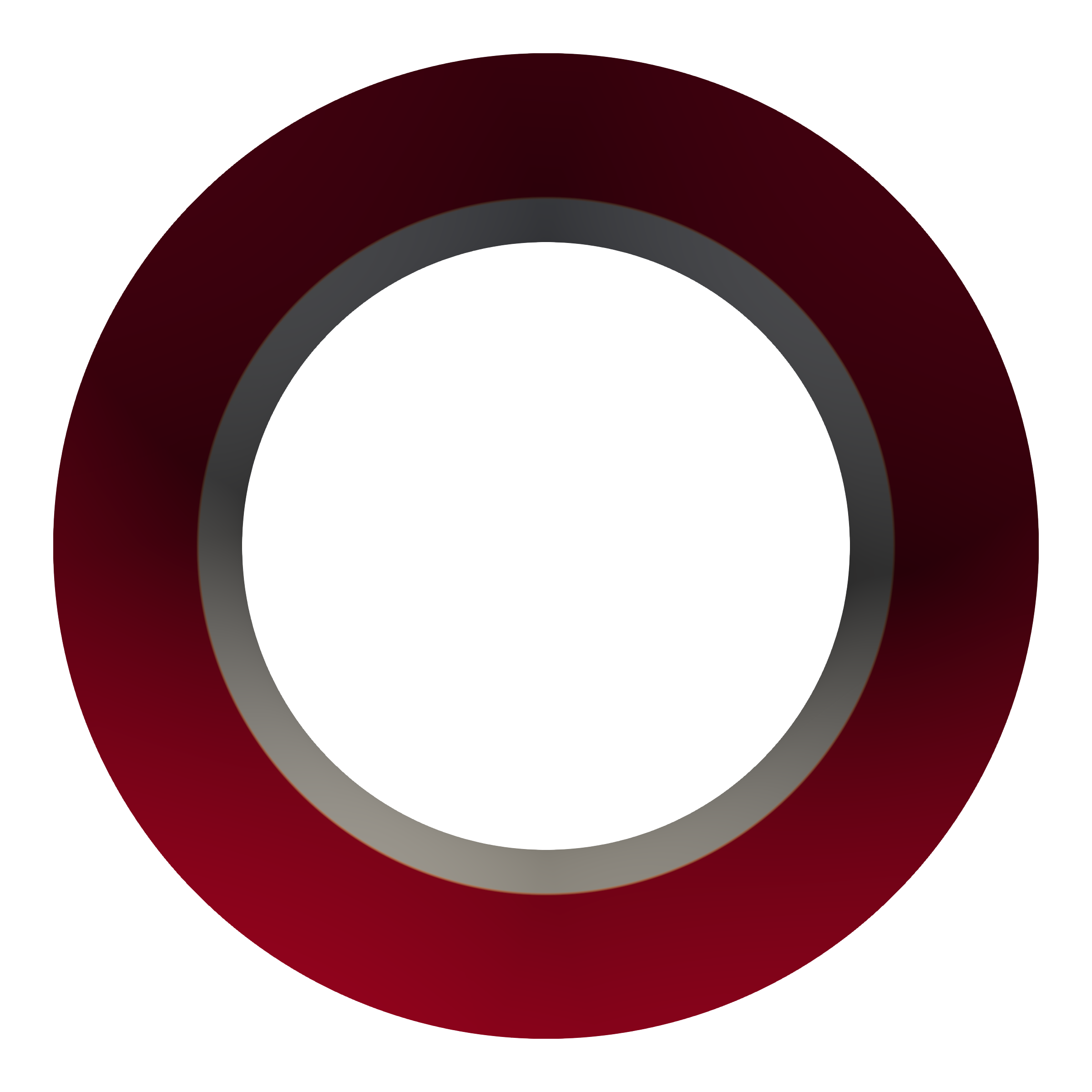}
        \caption{Top view of the deformed cylinder using the tension field theory membrane model. The red region denotes a taut, and the gray region denotes a wrinkled region.}
        \label{fig:cylinder_top_tensionfield}
    \end{subfigure}
    \caption{Results of the example with an cylinder with fixed bottom boundary and with a top boundary subject to a translation and a rotation, see \cref{fig:cylinder_setup}. A side view (\subref{fig:cylinder_side}) and a top view (\subref{fig:cylinder_top}) of the wrinkled membrane using the Kirchhoff--Love shell simulation are provided, as well as the deformed geometry from the tension field theory membrane simulation with the tension field for colouring (\subref{fig:cylinder_top_tensionfield}).}
    \label{fig:cylinder}
\end{figure}

\subsection{Multi-Patch Analysis}\label{subsec:examples_coupling}
In order to model Kirchhoff-Love shell problems on complex multi-patch geometries, $C^1$ continuity of the basis needs to be ensured on patch interfaces. To meet this continuity requirement, various techniques can be used, as reviewed in \cite{Verhelst2023Coupling}. Among these techniques are the penalty coupling method for shells \cite{Herrema2019} and unstructured spline constructions. In this section, two examples for multi-patch shell analysis are given: using penalty coupling (\cref{subsubsec:examples_coupling_weak}) and using unstructured spline constructions (\cref{subsubsec:examples_coupling_strong}).

\subsubsection{Weak Coupling}\label{subsubsec:examples_coupling_weak}
Consider the T-beam depicted in \cref{fig:gismo_examples_coupling_weak_beam_setup} fixed at its side, $\Gamma$. The beam has length $L=10.0\:[\text{m}]$, flange width $w=2.0\:[\text{m}]$, height $h=2.0\:[\text{m}]$ and thickness $t=0.1\:[\text{m}]$ and is exposed to a point load $F=10.0\:[\text{N}]$ and has material properties $E=10.0\cdot 10^7\:[\text{Pa}]$, $\nu=0.0\:[-]$, modelled through a linear Saint-Venant Kirchhoff material model. The beam can be modelled using a multi-patch geometry with three patches: one for the web and two for the flange, as shown in the work by \cite{Herrema2019}. On the interface of these patches, penalty coupling for the isogeometric Kirchhoff-Love shell can be applied to find the response of the beam. Contrary to unstructured splines, the penalty method provides versatility in interface coupling, for instance, coupling patches with non-matching interfaces or patches that are joined with $C^0$ continuity. However, the penalty method requires a penalty coefficient, which needs to be selected a priori.\\

\begin{figure}
\centering
\begin{minipage}[t]{\halfwidth}
\centering
\begin{tikzpicture}[]
\tdplotsetmaincoords{70}{60}
\begin{scope}[xshift=4cm,tdplot_main_coords]
\fill[draw=black,top color=black!60,bottom color=black!20](0,0,0) -- (10,0,0) node[](E){} -- (10,0,-2) coordinate (B){}-- (0,0,-2) coordinate (A){}--cycle;
\draw[thick] (0,-1,0) -- (0,1,0) coordinate[pos=0.75](midtop){};
\draw[thick] (0,0,-2) -- (0,0,0) coordinate[midway](midbelow){};
\draw[latex-](midtop) to[in=-90,out=110] (-0.5,0,0.5);
\draw[latex-](midbelow) to[in=-90,out=140] (-0.5,0,0.5) node[above]{$\Gamma$};

\fill[draw=black,top color=black!60,bottom color=black!20,fill opacity=0.7](0,1,0) -- (10,1,0) coordinate[midway](C)-- (10,-1,0) -- (0,-1,0) coordinate[midway](D){}--cycle;
\draw[latex-,thick](10,-1,0)--(10,-1,1) node[above left]{$F$};
\draw[latex-latex] ($(A)+(0,0,-0.2)$)--($(B)+(0,0,-0.2)$) node[midway,below left]{$L$};
\draw[latex-latex] (C)--(D) node[midway,below]{$w$};
\draw[latex-latex] ($(B)+(0.2,0,0)$)--($(E)+(0.2,0,0)$) node[midway,right]{$h$};
\end{scope}
\end{tikzpicture}
\caption{Problem setup for the T-beam example, inspired by \cite{Herrema2019}.}
\label{fig:gismo_examples_coupling_weak_beam_setup}
\end{minipage}
\hfill
\begin{minipage}[t]{\halfwidth}
\centering
\begin{tikzpicture}
\begin{axis}
[
legend pos = north east,
width=\linewidth,
height=0.3\textheight,
grid=major,
legend style={fill=white, fill opacity=0.6, draw opacity=1,draw=none,text opacity=1},
xlabel={Penalty parameter $\alpha\:[-]$},
ylabel={Free-end rotation $\theta\:[{}^\circ]$},
xmode=log,
ymode=normal,
restrict y to domain=0:5,
]
\addplot+[solid,no markers] table[col sep=comma,header=true,x index={1},y index={9}] {Data/Beam_alpha.csv};
\end{axis}
\end{tikzpicture}
\caption{The rotation of the end-point of the beam with respect to the penalty parameter $\alpha$.}
\label{fig:gismo_examples_coupling_weak_beam_plot}
\end{minipage}
\end{figure}

Penalty coupling is established by adding extra energy contributions to the weak formulation of the isogeometric Kirchhoff-Love shell (see \cite{Herrema2019}). The penalty method is implemented using a single parameter for both rotations and displacements, according to the work by \cite{Herrema2019}. To activate penalty coupling on selected interfaces, stored in \gs{interfaces}, of a multi-patch geometry, only the following options need to be set to the \gs{gsThinShellAssembler}:

\begin{lstlisting}
// Provided mp, assembler
std::vector<boundaryInterface> interfaces = mp.interfaces();
assembler.addWeakC0(interfaces);
assembler.addWeakC1(interfaces);
assembler.initInterfaces();
assembler.options().setReal("PenaltyIfc",1e3);
\end{lstlisting}
Here, \gs{interfaces} is a vector for storing all interfaces of the multi-patch geometry. This vector can be modified to apply penalty coupling to a selection of the patch interfaces. Furthermore, weak boundary conditions as proposed by \cite{Herrema2019} have been implemented in the \gs{gsThinShellAssembler}. Finally, using the code lines provided, the response of the beam can be computed for different penalty parameters, as shown in \cref{fig:gismo_examples_coupling_weak_beam_plot}.
\subsubsection{Strong Coupling}\label{subsubsec:examples_coupling_strong}
Provided the geometry of the car from \cref{example:gismo,example:gsUnstructuredSpline}, with matching patches. On this geometry, an unstructured spline basis can be constructed such that $C^1$ isogeometric analysis can be performed. In particular, a modal analysis of this geometry is performed to infer the eigenfrequencies of the car part. The material used in this example is steel with material parameters $E=210\cdot10^3\:[\text{MPa}]$, $\nu=0.3\:[-]$, $\rho=7850\:[\text{tonnes}/\text{mm}^3]$ and a thickness of $t=10\:[\text{mm}]$. This example is inspired by \cite{Verhelst2023Coupling}, to which the reader is referred for a comparison of multiple methods with respect to a commercial finite element code.\\

\begin{figure}
    \centering
    \begin{subfigure}{\linewidth}
        \centering
        \begin{minipage}{\halfwidth}
            \centering
                \includegraphics[width=\linewidth,trim=150 100 050 00, clip]{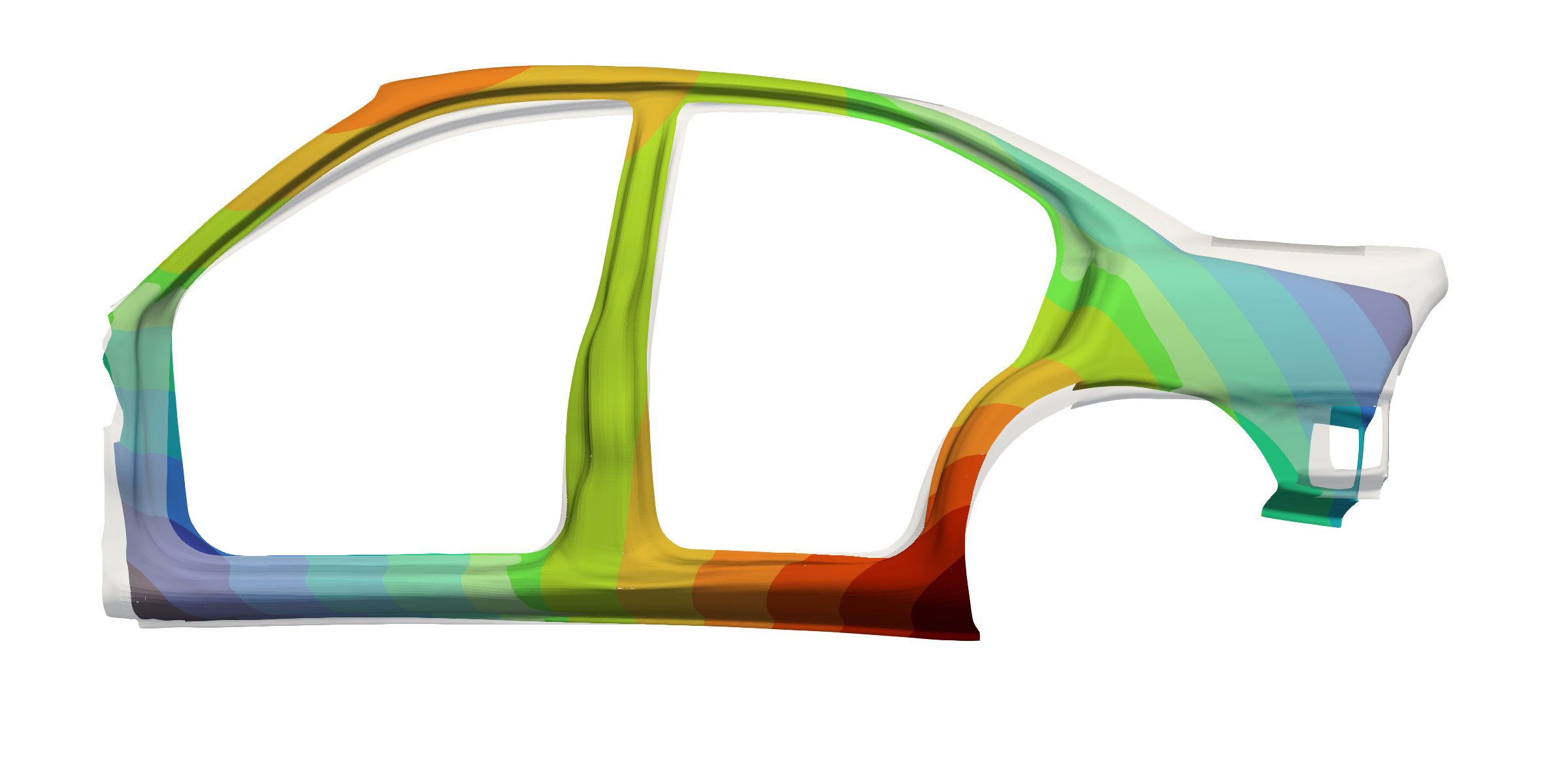}
        \end{minipage}
        \hfill
        \begin{minipage}{\halfwidth}
            \centering
                \includegraphics[width=\linewidth,trim=552 150 352 00, clip]{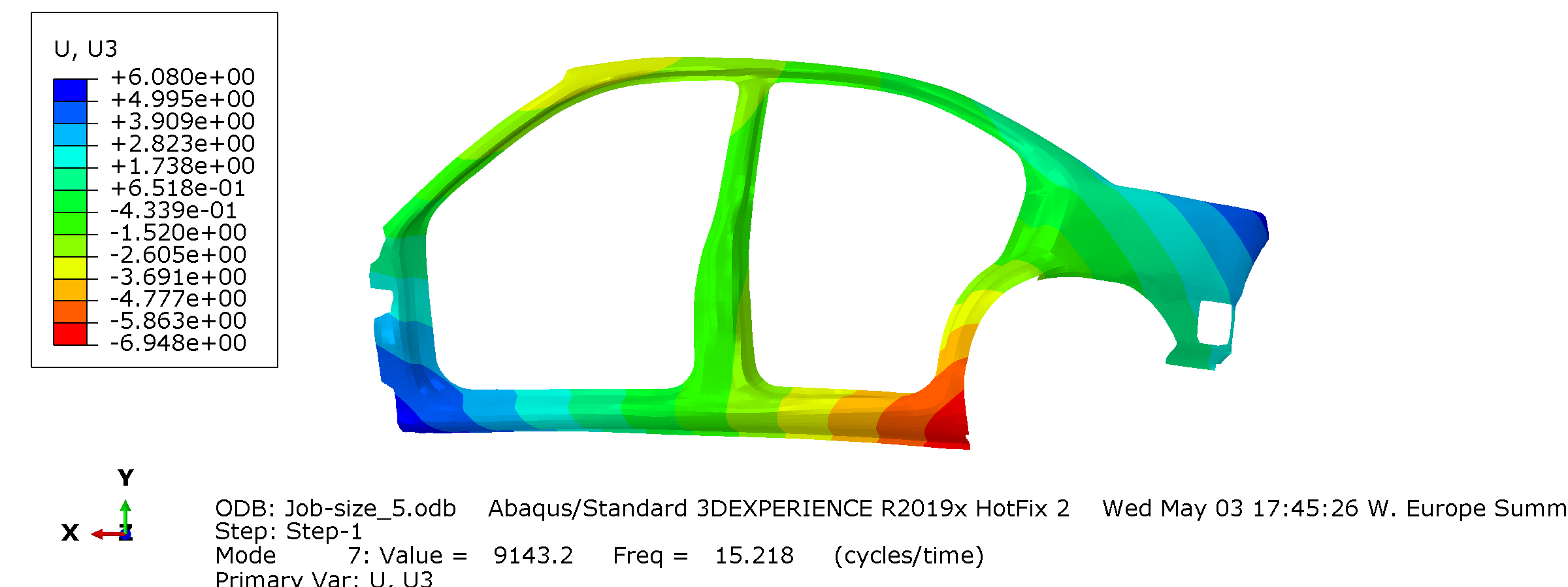}
        \end{minipage}
        \caption{Mode 1}
    \end{subfigure}

    \begin{subfigure}{\linewidth}
        \centering
        \begin{minipage}{\halfwidth}
            \centering
                \includegraphics[width=\linewidth,trim=150 100 050 00, clip]{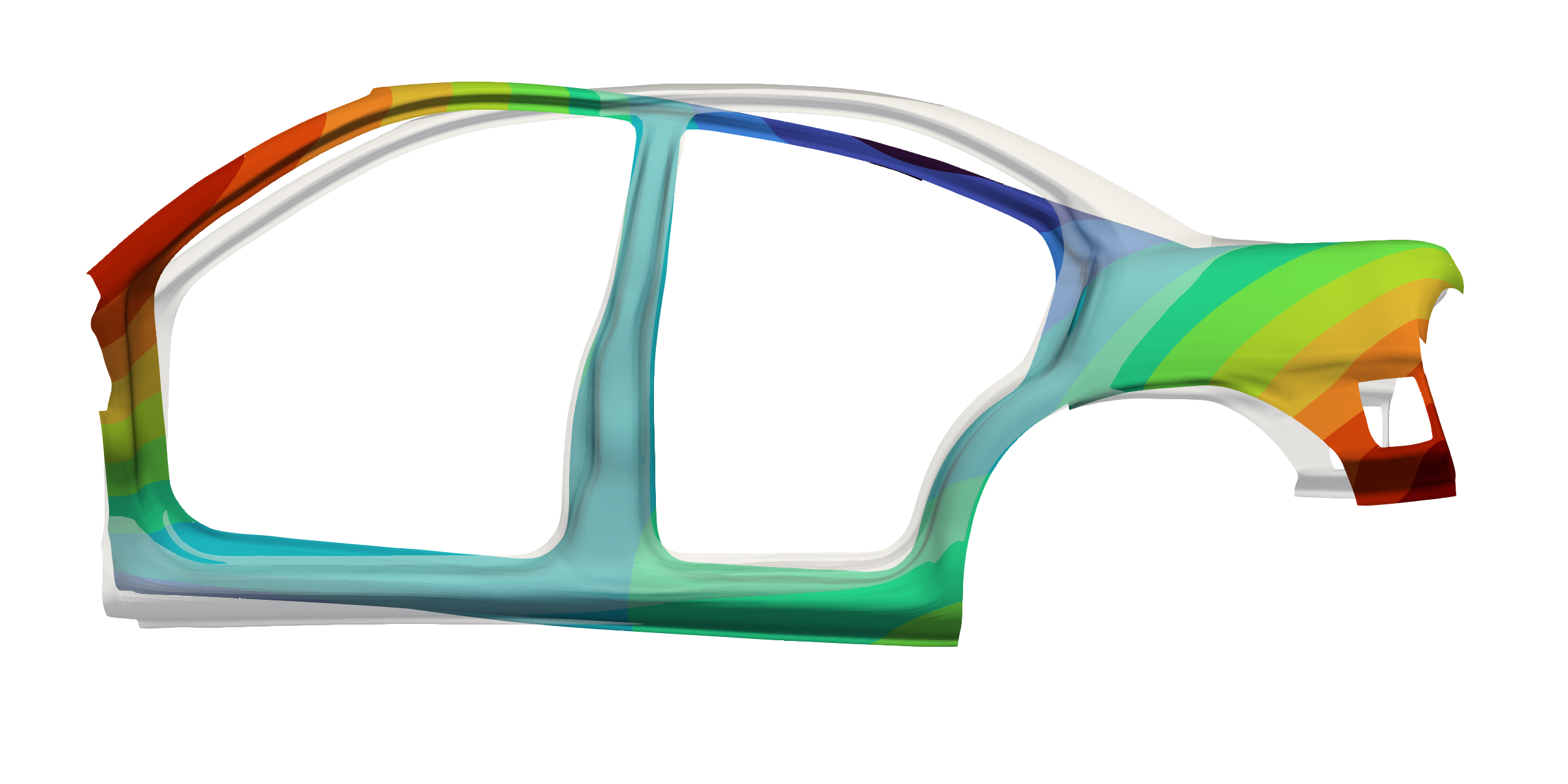}
        \end{minipage}
        \hfill
        \begin{minipage}{\halfwidth}
            \centering
                \includegraphics[width=\linewidth,trim=552 150 352 00, clip]{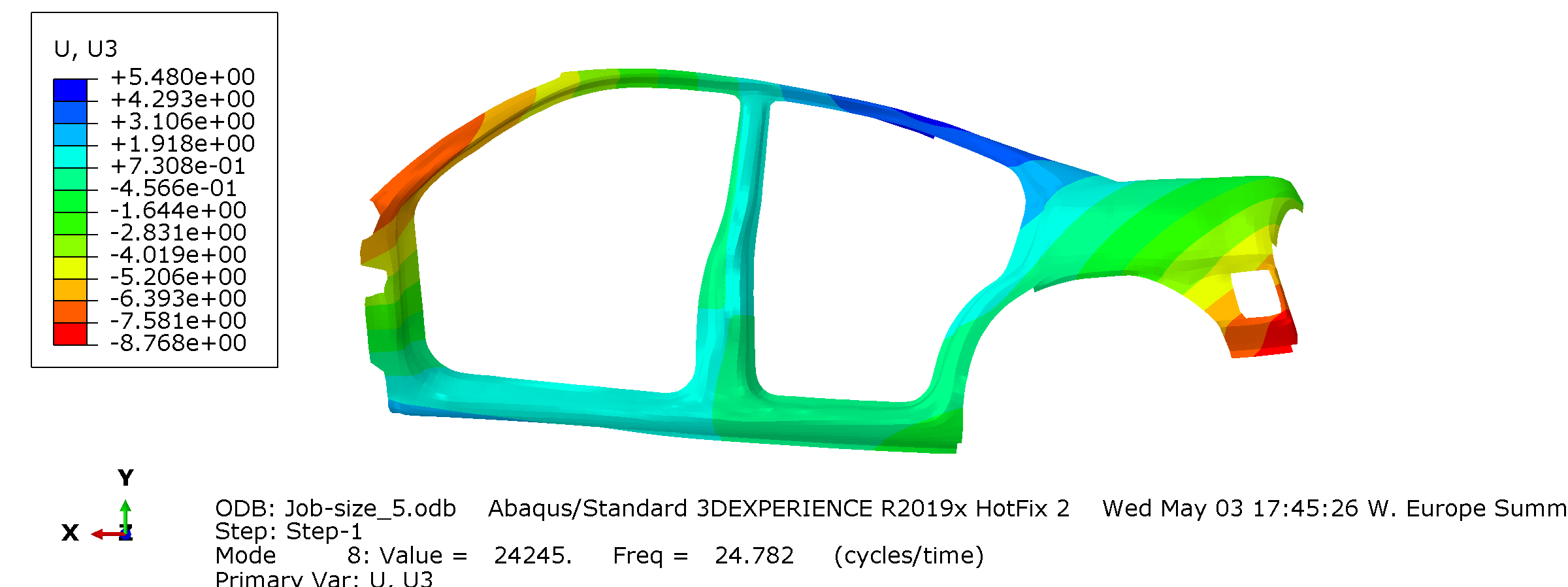}
        \end{minipage}
        \caption{Mode 2}
    \end{subfigure}

    \begin{subfigure}{\linewidth}
        \centering
        \begin{minipage}{\halfwidth}
            \centering
                \includegraphics[width=\linewidth,trim=150 100 050 00, clip]{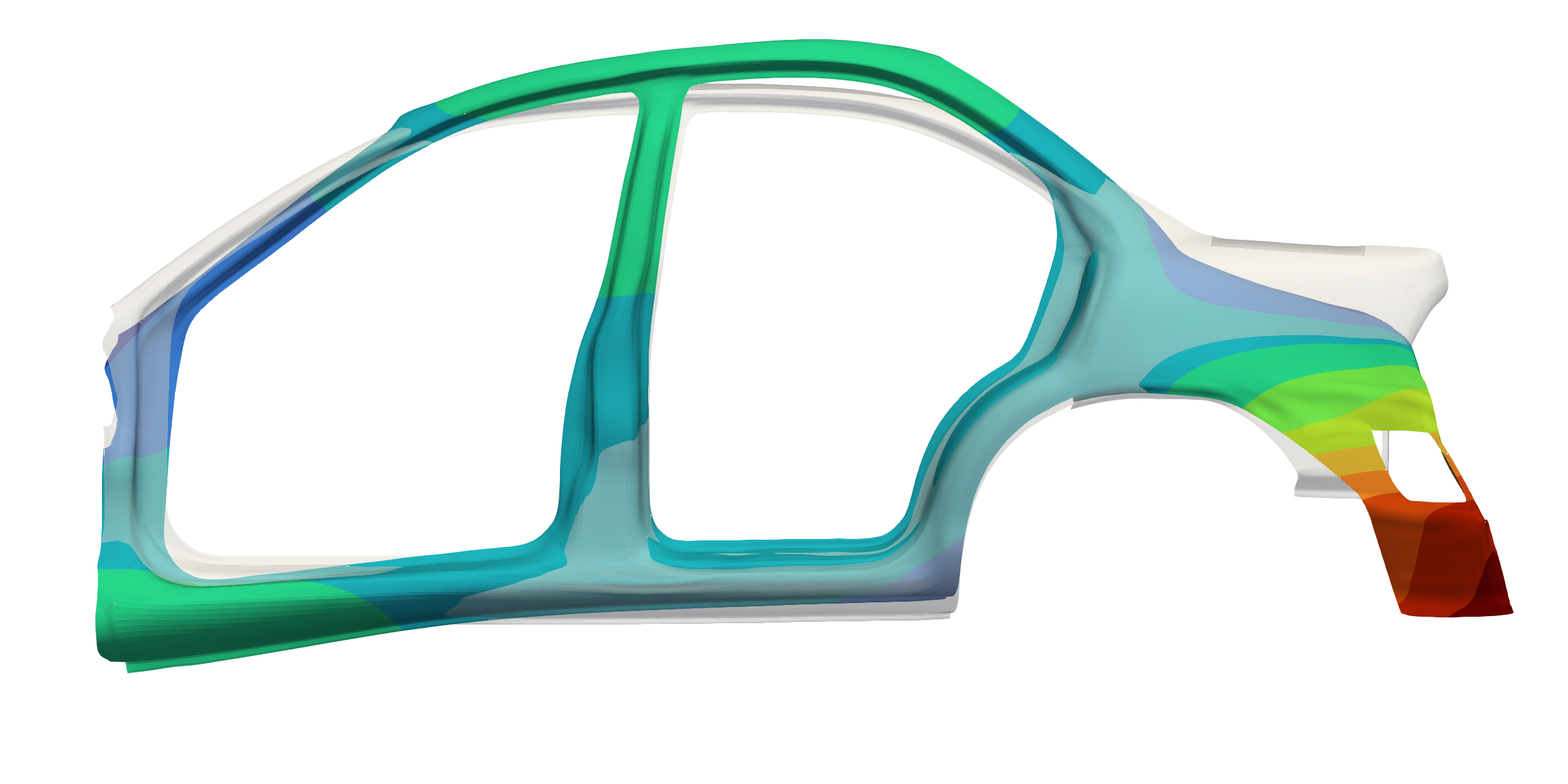}
        \end{minipage}
        \hfill
        \begin{minipage}{\halfwidth}
            \centering
                \includegraphics[width=\linewidth,trim=552 150 352 00, clip]{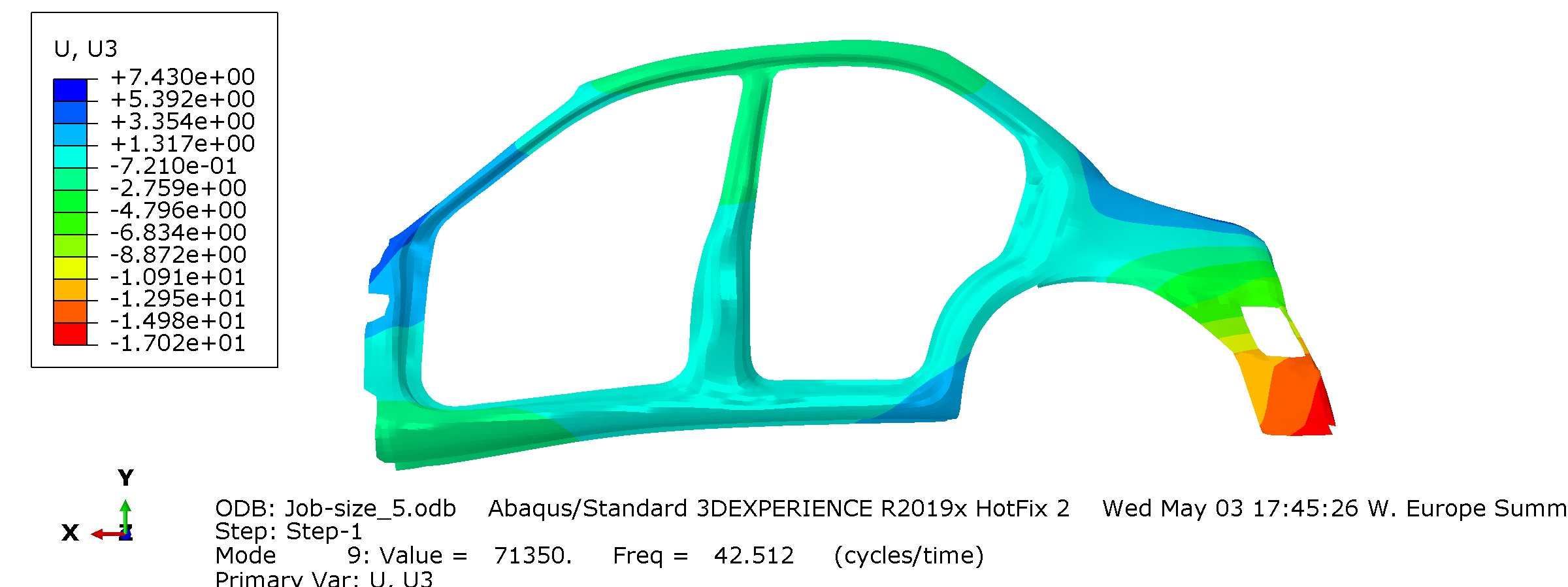}
        \end{minipage}
        \caption{Mode 3}
    \end{subfigure}

    \begin{subfigure}{\linewidth}
        \centering
        \begin{minipage}{\halfwidth}
            \centering
                \includegraphics[width=\linewidth,trim=150 100 050 00, clip]{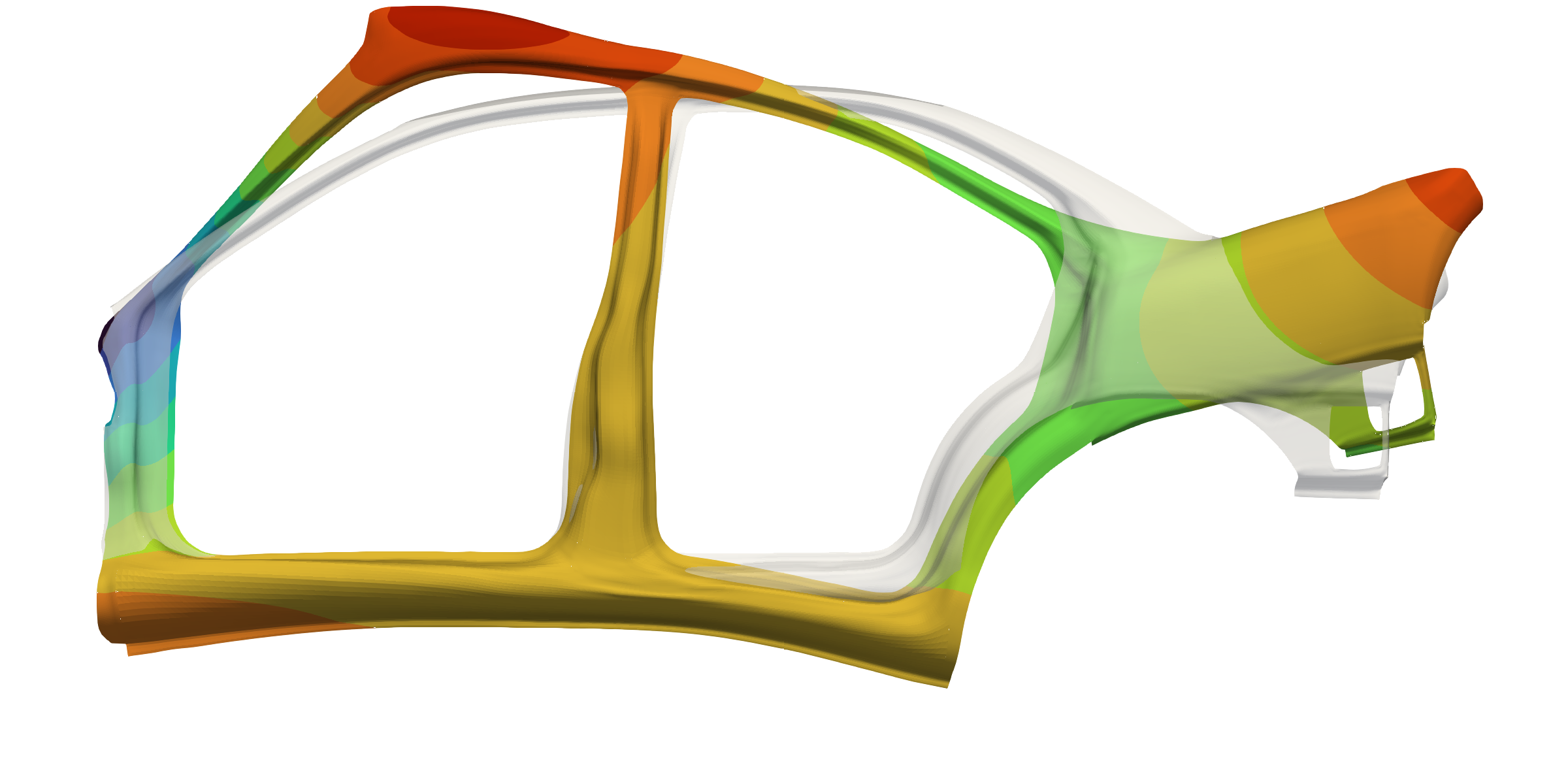}
        \end{minipage}
        \hfill
        \begin{minipage}{\halfwidth}
            \centering
                \includegraphics[width=\linewidth,trim=552 150 352 00, clip]{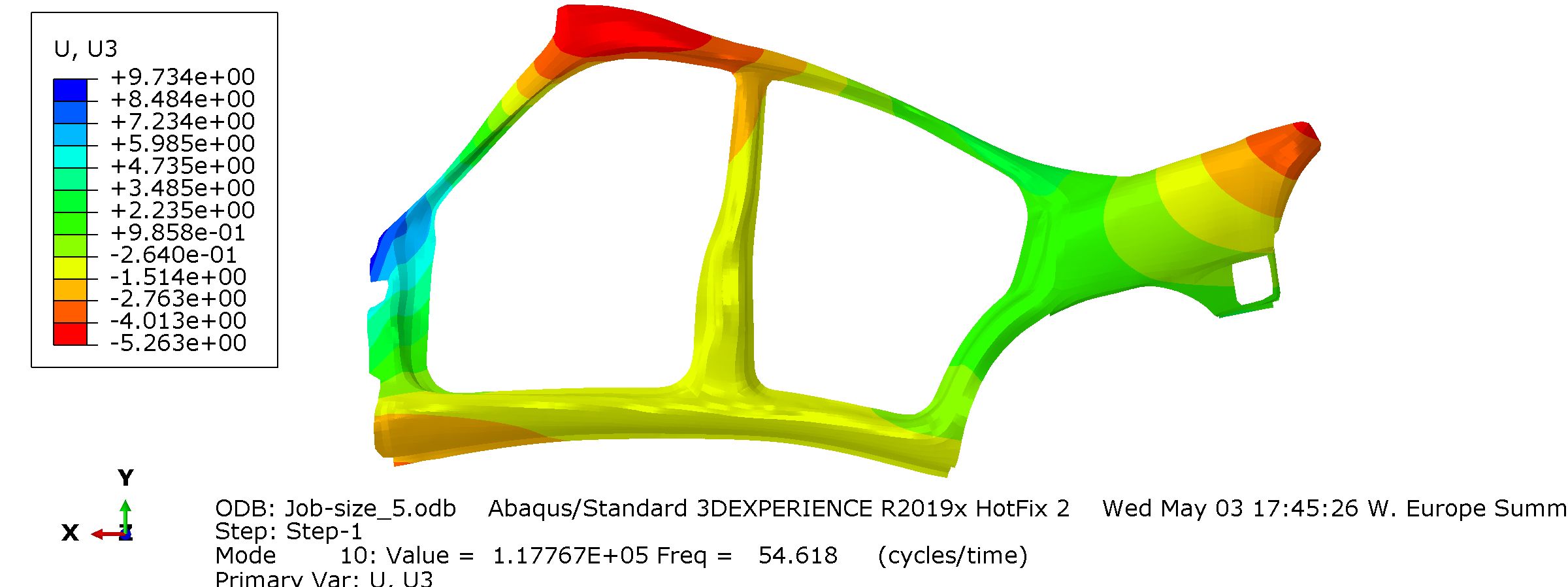}
        \end{minipage}
        \caption{Mode 4}
    \end{subfigure}

        %

        %

        %
    \caption{Out-of-plane deformations of the first four vibration modes of the side of the car from \cref{fig:coupling_car_side}. The results on the left represent the results obtained by the D-Patch construction and the results on the right represent results obtained using ABAQUS (10mm). The mode shapes are all deformation modes warped by the deformation vector and plotted over the undeformed (transparent) geometry.}
    \label{fig:coupling_car_modes}
\end{figure}

Strong coupling of multi-patches is enabled by using unstructured spline constructions, enabled through the \gs{gsUnstructuredSplines} module (see \cref{sec:gsUnstructuredSplines}). To perform shell analysis on \gs{gsMappedBasis} and \gs{gsMappedSpline}, the \gs{gsThinShellAssembler} is generalised to work on geometries defined by arbitrary functions (i.e., inherited from \gs{gsFunctionSet}), which in turn support computation of first and second derivatives, normals, and other geometric data. Moreover, the \gs{gsThinShellAssembler} allows to store a \gs{gsMappedBasis} object with an underlying local basis provided as a \gs{gsMultiBasis} used for quadrature definitions. Using the \gs{mbasis} and \gs{mspline} defined in \cref{example:gsUnstructuredSpline}, the \gs{gsThinShellAssembler} can be configured to use unstructured splines as follows:
\begin{lstlisting}
// Provided bases, mbasis, mspline, and assembler
assembler.setSpaceBasis(bases);
assembler.setSmoothBasis(mbasis);
assembler.setGeometry(mspline);
\end{lstlisting}
Where the objects \gs{bases}, \gs{mbasis} and \gs{mspline} follow, for example, from the Almost-$C^1$ method \cite{Takacs2023}, as demonstrated in \cref{example:gsUnstructuredSpline}.
For the modal analysis, the \gs{gsModalSolver} solver is used, compiled with \gs{Spectra}. Using the shell assembler, this solver is simply set up as:

\begin{lstlisting}
// Provided assembler
gsSparseMatrix<real_t> K, M;
assembler.assemble();        // Asseble eqs. @\ref{eq:gismo_gsKLShell_Linear}@ and @\ref{eq:gismo_externalForce}@
K = assembler.matrix();      // Collect eq. @\ref{eq:gismo_gsKLShell_Linear}@
assembler.assembleMass();    // Assemble eq. @\ref{eq:gismo_massMatrix}@
M = assembler.matrix();      // Collect eq. @\ref{eq:gismo_massMatrix}@
gsModalSolver<> solver(K,M); // Define the solver
solver.compute();            // Compute eq. @\ref{eq:gismo_ModalAnalysisProblem}@
\end{lstlisting}
Using this solver with the almost-$C^1$ basis, the vibration modes of the car geometry are obtained, as presented in \cref{fig:coupling_car_modes} and \cref{tab:ModalShell}. For more details, the reader is referred to \cite{Verhelst2023Coupling}.

\begin{table}[]
    \centering
    \caption{Eigenfrequencies of the Almost-$C^1$ and D-Patch constructions for the car geometry in \cref{fig:gismo_car_patches}. The results of an ABAQUS FEA simulation using the S4R element are provided as a reference. The mode shapes are plotted in \cref{fig:coupling_car_modes}.}
    \label{tab:ModalShell}
    \begin{tabular}{lll|llll}
        \toprule
        Method & & \# DoFs & Mode 1 & Mode 2 & Mode 3 & Mode 4\\
        \midrule
        \multicolumn{2}{l}{\multirow[c]{4}{*}{Almost-$C^1$, $p=2$, $r=1$}}   & 13,731     & 15.740 & 25.567 & 43.829 & 56.654 \\
        \multicolumn{2}{l}{}                                                 & 49,758     & 15.762 & 25.564 & 43.429 & 56.778 \\
        \multicolumn{2}{l}{}                                                 & 189,654    & 15.776 & 25.552 & 43.269 & 56.785 \\
        \multicolumn{2}{l}{}                                                 & 740,814    & 15.774 & 25.531 & 43.177 & 56.746 \\
        \midrule
        \multicolumn{2}{l}{\multirow[c]{3}{*}{D-Patch, $p=2$, $r=1$}}        & 49,437     & 15.785 & 25.607 & 43.641 & 56.902\\
        \multicolumn{2}{l}{}                                                 & 189,333    & 15.780 & 25.561 & 43.323 & 56.807\\
        \multicolumn{2}{l}{}                                                 & 740,493    & 15.775 & 25.533 & 43.191 & 56.748 \\
        \midrule
        \multicolumn{2}{l}{\multirow[c]{2}{*}{D-Patch, $p=3$, $r=1$}}        & 136,839    & 15.749 & 25.593 & 43.348 & 56.786 \\
        \multicolumn{2}{l}{}                                                 & 630,459    & 15.760 & 25.581 & 43.231 & 56.801 \\
        \midrule
        \multicolumn{2}{l}{\multirow[c]{2}{*}{D-Patch, $p=3$, $r=2$}}        & 71,760     & 15.771 & 25.539 & 43.224 & 56.744 \\
        \multicolumn{2}{l}{}                                                 & 226,524    & 15.755 & 25.582 & 43.235 & 56.807 \\
        \midrule
        \multirow{3}*{ABAQUS S4R}& 10mm                                      & 126,966    & 15.303 & 24.881 & 42.629 & 54.887 \\
        & 5mm                                       & 440,076    & 15.224 & 24.780 & 42.516 & 54.627 \\
        & 2.5mm                                     & 1,653,030     & 15.119 & 24.640 & 42.338 & 54.277 \\

        \bottomrule
    \end{tabular}
\end{table}

\subsection{Structural Stability Analysis}\label{subsec:examples_ALM}
In order to demonstrate the \gs{gsALMBase} and \gs{gsAPALM} solvers from the \gs{gsStructuralAnalysis} class, two examples regarding the analysis of structural instabilities are presented. Firstly, the set-up of an arc-length simulation with singular point computation for solving a bifurcation problem without imposing initial perturbations is presented. This method relies on solving the system in \cref{eq:gismo_gsStructuralAnalysis_CrisfieldExtended}. Thereafter, the use of the Adaptive Parallel Arc-Length Method from \cite{Verhelst2023APALM} is presented.\\

\subsubsection{Bifurcation Instabilities}\label{subsubsec:examples_ALM_bifurcation}
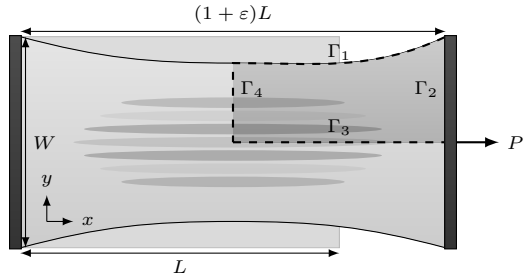
\begin{figure}
    \centering
    \begin{tikzpicture}[scale=0.7]
        \footnotesize
        \def\B{4}
        \def\L{6}
        \def\strain{2};
        \def\overlap{0.02}

        \fill[opacity=0.15, bottom color=black!10, top color=black!05, draw=black] (0,0) -- (\L,0) -- (\L,\B) -- (0,\B) --cycle;

        \filldraw[draw=black,bottom color=black!20, top color=black!10] (0,0) to[out=15,in=180] (\L/2+\strain/2,0.5) to[out=0,in=165] (\L+\strain,0) -- (\L+\strain,\B) to[out=205,in=0] (\L/2+\strain/2,\B-0.5) to[out=180,in=-15] (0,\B) -- cycle;

        \filldraw[dashed,draw=black,bottom color=black!30, top color=black!20,thick] (\L/2+\strain/2,\B/2)  -- (\L+\strain,\B/2) -- (\L+\strain,\B) to[out=205,in=0] (\L/2+\strain/2,\B-0.5) -- cycle;

        \draw[latex-latex] (0,\B+0.1)--node [above]{$(1+\varepsilon)L$}(\L+\strain,\B+0.1);
        \draw[latex-latex] (0,-0.1)--node [below]{$L$}(\L,-0.1);
        \draw[latex-latex] (0.1,0)--node [right]{$W$}(0.1,\B);

        \node at (3*\L/4+3*\strain/4,\B)[rotate=0,below]{$\Gamma_1$};
        \node at (\L+\strain,3*\B/4)[rotate=0,left]{$\Gamma_2$};
        \node at (3*\L/4+3*\strain/4,\B/2)[rotate=0,above]{$\Gamma_3$};
        \node at (\L/2+\strain/2,3*\B/4)[rotate=0,right]{$\Gamma_4$};

        \draw[-latex,thick] (\L+\strain,\B/2) -- (\L+\strain+1,\B/2) node[right]{$P$};

        \filldraw[bottom color=black!80,top color=black!70] (\overlap,-\overlap) -- (\overlap,\B+\overlap) -- (-0.2,\B+\overlap) -- (-0.2,-\overlap)--cycle;
        \filldraw[bottom color=black!80,top color=black!70] (\L-\overlap+\strain,-\overlap) -- (\L-\overlap+\strain,\B+\overlap) -- (\L+0.2+\strain,\B+\overlap) -- (\L+0.2+\strain,-\overlap)--cycle;

        \fill[color=black,opacity=0.1](0.5*\L+0.5*\strain,\B/2) ellipse [x radius=\L/2,y radius=0.1];
        \fill[color=black,opacity=0.2](0.5*\L+0.5*\strain,\B/2+0.25) ellipse [x radius=\L/2-0.2,y radius=0.1];
        \fill[color=black,opacity=0.2](0.5*\L+0.5*\strain,\B/2-0.25) ellipse [x radius=\L/2-0.2,y radius=0.1];
        \fill[color=black,opacity=0.05](0.5*\L+0.5*\strain,\B/2+0.5) ellipse [x radius=\L/2-0.5,y radius=0.1];
        \fill[color=black,opacity=0.05](0.5*\L+0.5*\strain,\B/2-0.5) ellipse [x radius=\L/2-0.5,y radius=0.1];
        \fill[color=black,opacity=0.15](0.5*\L+0.5*\strain,\B/2+0.75) ellipse [x radius=\L/2-0.9,y radius=0.1];
        \fill[color=black,opacity=0.15](0.5*\L+0.5*\strain,\B/2-0.75) ellipse [x radius=\L/2-0.9,y radius=0.1];

        \draw[-latex] (0.5,0.5) -- node[right,inner sep=8pt] {$x$}(1.0,0.5);
        \draw[-latex] (0.5,0.5) -- node[above,inner sep=8pt] {$y$}(0.5,1.0);

    \end{tikzpicture}
    \caption{modeling geometry for the wrinkled sheet. The sheet is modelled as only a quarter using symmetry conditions on $\Gamma_3$ and $\Gamma_4$. The side $\Gamma_1$ is free and the side $\Gamma_2$ is clamped and fixed in all direction except the $x$-direction. The sheet has length $L$ and width $W$. The strain of the sheet due to load $P$ is denoted by $\varepsilon$.}
    \label{fig:wrinkling}
\end{figure}

Provided a rectangular membrane with length $L=280\:[\text{mm}]$, width $W=140\:[\text{mm}]$ and thickness $t=0.14\:[\text{mm}]$, with material parameters $\mu=\mu_1+\mu_2=12.43\cdot10^4 + 31.62\cdot10^4\:[\text{Pa}]$ for an incompressible Mooney-Rivlin material model ($\nu=0.5\:[-]$), which is fixed on the left end and clamped but pulled on the right end, as depicted in \cref{fig:wrinkling}. When the load $P$ is increased, tension is applied, and longitudinal wrinkles will appear and disappear as the strain $\varepsilon$ increases. This benchmark example is based on the work of \cite{Panaitescu2019} and replicated from \cite{Verhelst2021}. This benchmark example can be modelled using arc-length methods with the extended arc-length method to compute the bifurcation point; see \cref{sec:gsStructuralAnalysis_ArcLength}, which will be demonstrated in the following.\\

\begin{figure}[tb!]
    \centering
    \begin{subfigure}{\linewidth}
        \begin{tikzpicture}
            \footnotesize
            \begin{axis}
                [
                xlabel=$\varepsilon$,
                ylabel=$\max(\zeta)/t$,
                legend pos = north east,
                xmin = 0,
                xmax = 0.5,
                ymin = 0,
                ymax = 3.5,
                width=\linewidth,
                height=0.3\textheight,
                grid=major,
                ]
                \addplot+[mark=o,mark size=2,black] table[header=true,x index = {0},y index = {1}, col sep = comma]{Data/Wrinkling/PanaitescuExperiment.csv};\addlegendentry{Exp.}
                \addplot [name path=upper,draw=none,mark=none,forget plot] table[header=true,x index = {2},y index = {3}, col sep = comma]{Data/Wrinkling/PanaitescuExperiment.csv};
                \addplot [name path=lower,draw=none,mark=none,forget plot] table[header=true,x index = {4},y index = {5}, col sep = comma]{Data/Wrinkling/PanaitescuExperiment.csv};
                \addplot [fill=black!10,forget plot] fill between[of=upper and lower];

                \addplot+[style=style0,dashed,no marks,plotcol1!60!black] table[header=true,x index = {0},y index = {1}, col sep = comma]{Data/Wrinkling/LSDynaFull.csv};\addlegendentry{FI}
                \addplot+[style=style1,dotted,no marks,plotcol1!80!black] table[header=true,x index = {0},y index = {1}, col sep = comma]{Data/Wrinkling/LSDynaReduced.csv};\addlegendentry{H-L}
                \addplot+[style=style2,dashdotted,no marks,plotcol1!100!black] table[header=true,x index = {0},y index = {1}, col sep = comma]{Data/Wrinkling/LSDynaS-RReduced.csv};\addlegendentry{H-L S/R}
                \addplot+[style=style0,no marks,plotcol5] table[header=true,x index = {0},y index = {1}, col sep = comma]{Data/Wrinkling/ANSYS.csv};\addlegendentry{SHELL181}

                \pgfplotsset{cycle list shift=-6}

                \addplot+[style=style1,no markers,plotcol2] table[header=true,x index = {0},y index = {1}, col sep = comma]{Data/Wrinkling/IGA_Quart_r5e2_MR.csv};\addlegendentry{MR}
                \addplot+[style=style2,no markers,plotcol3] table[header=true,x index = {0},y index = {1}, col sep = comma]{Data/Wrinkling/IGA_Quart_r5e2_OG.csv};\addlegendentry{OG}
            \end{axis}
        \end{tikzpicture}
        \caption{Strain-amplitude diagram of the tension wrinkling of a thin sheet. The vertical axis represents the maximum amplitude normalised by the shell thickness $t$ and the horizontal axis represents the strain $\varepsilon$ of the sheet. The present model is used to obtain the Mooney-Rivlin (MR) and Ogden (OG) results. The fully integrated (FI), Hughes-Liu (H-L), and Hughes-Liu Selective/Reduced (H-L S/R) results are obtained using LS-DYNA, and the SHELL181 results are obtained using ANSYS. The experimental results (Exp.) from \cite{Panaitescu2019} are plotted as a reference.}
        \label{fig:hyperelasticity_wrinkling_results}
    \end{subfigure}

    \begin{subfigure}{\halfwidth}
        \centering
        \begin{minipage}[b]{\halfwidth}
            \centering
            \includegraphics[trim = 0 100 0 0, clip, width=\linewidth]{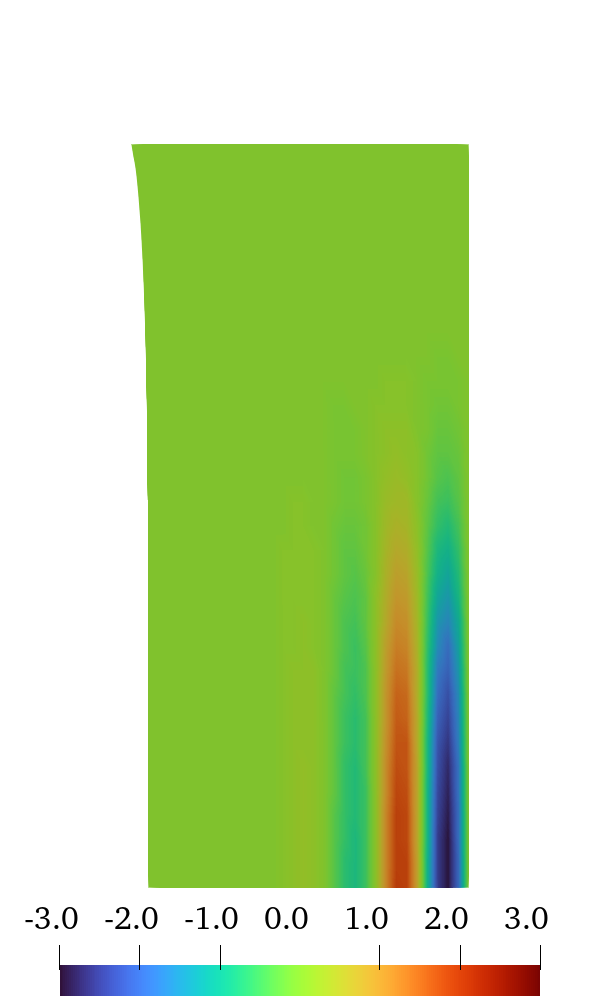}
            \caption*{$\varepsilon=0.1$}
        \end{minipage}
        \begin{minipage}[b]{\halfwidth}
            \centering
            \includegraphics[trim = 0 100 0 0, clip, width=\linewidth]{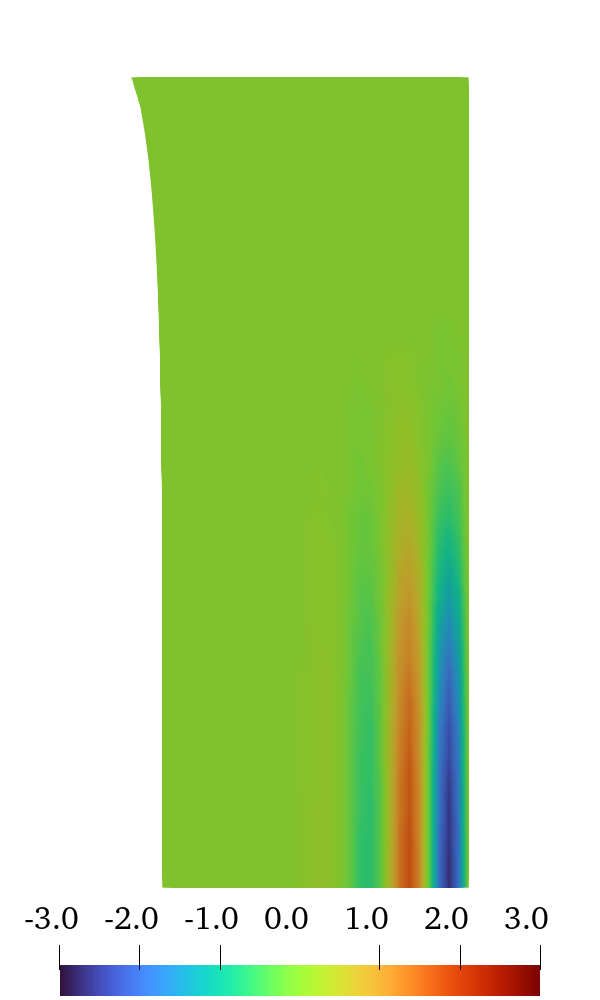}
            \caption*{$\varepsilon=0.2$}
        \end{minipage}

        \begin{minipage}{\fullwidth}
            \centering
            \captionsetup{justification=centering}
            \includegraphics[trim = 0 0 0 900, clip, width=0.8\linewidth]{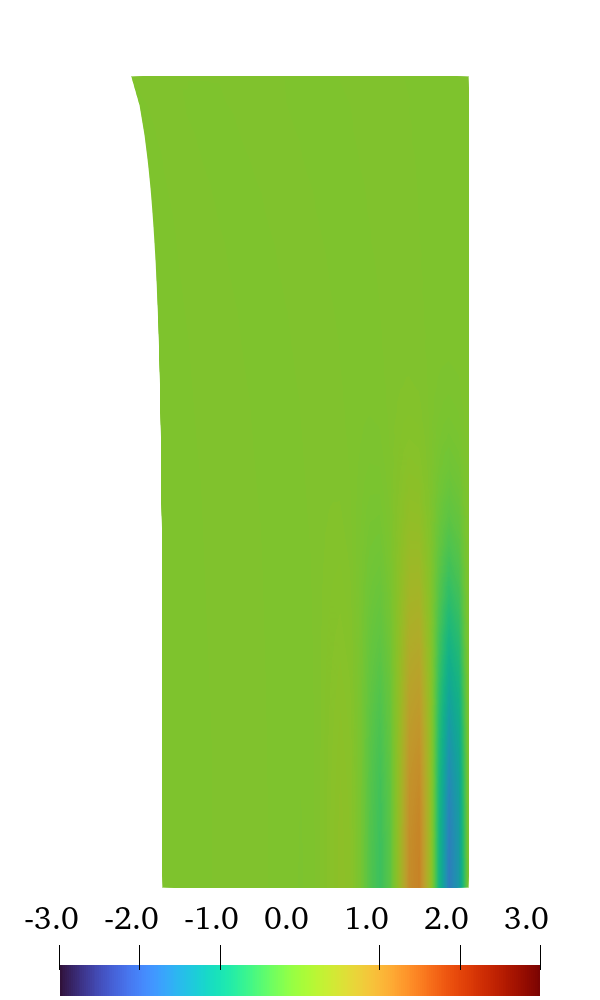}
            \caption*{Normalised wrinkling amplitude $\zeta/t\:[-]$}
        \end{minipage}
        \caption{Contour plot of out-of-plane displacements $w$ for different strains $\varepsilon$ for the MR model.}
    \end{subfigure}
    \hfill
    \begin{subfigure}{\halfwidth}
        \centering
        \begin{minipage}[b]{\halfwidth}
            \centering
            \includegraphics[trim = 0 100 0 0, clip, width=\linewidth]{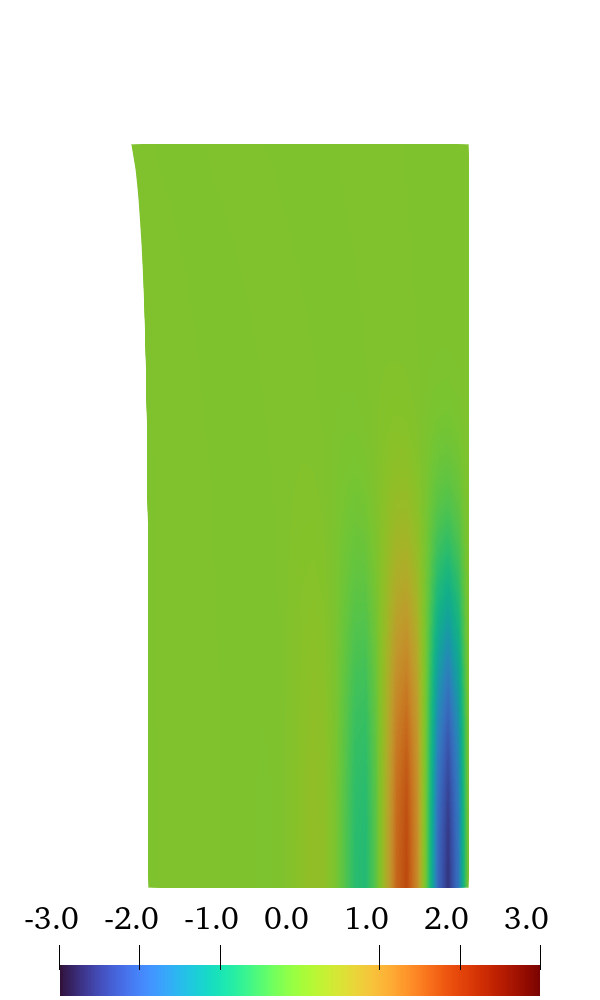}
            \caption*{$\varepsilon=0.1$}
        \end{minipage}
        \begin{minipage}[b]{\halfwidth}
            \centering
            \includegraphics[trim = 0 100 0 0, clip, width=\linewidth]{Figures/WrinklingOG/newstep60}
            \caption*{$\varepsilon=0.2$}
        \end{minipage}

        \begin{minipage}{\fullwidth}
            \centering
            \captionsetup{justification=centering}
            \includegraphics[trim = 0 0 0 900, clip, width=0.8\linewidth]{Figures/WrinklingOG/newstep60}
            \caption*{Normalised wrinkling amplitude $\zeta/t\:[-]$}
        \end{minipage}
        \caption{Contour plot of out-of-plane displacements $w$ for different strains $\varepsilon$ for the OG model.}
    \end{subfigure}

    \caption{Wrinkling formation in a thin sheet subject to tension.}
    \label{fig:hyperelasticity_wrinkling_figs}
\end{figure}

Provided a \gs{gsThinShellAssembler} and a \gs{gsALMCrisfield} arc-length method (see \cref{example:gsStructuralAnalysis}), an arc-length simulation including bifurcation point detection, computation, and branch switching is performed by:

\begin{lstlisting}
// Provided arcLength, Nsteps
index_t k=0;
while (k<Nsteps)
{
    gsStatus status = arcLength.step();
    if (status==gsStatus::NotConverged || status==gsStatus::AssemblyError)
    {
        arcLength.reduceLength(0.5);      // Multiplies the current length by 0.5
        continue;
    }
    arcLength.computeStability();
    if (arcLength.stabilityChange())      // if true, a limit point or bifurcation point
    {                                     // is found
        arcLength.computeSingularPoint(); // Compute eq. @\ref{eq:gismo_gsStructuralAnalysis_CrisfieldExtended}@
        arcLength.switchBranch();
        }
    arcLength.resetLength();
    k++;
}
\end{lstlisting}
As can be seen in the code snippet above, the arc length is reduced when the step is unsuccessful. Furthermore, the singular point computation is activated if the stability of the structure changes, which is checked through internal stability computations based on what is described in \cref{subsec:gsStructuralAnalysis_buckling}. In \cref{fig:hyperelasticity_wrinkling_results}, the resulting solution of the wrinkling benchmark is provided as the uniform meshing results. For more details, the reader is referred to \cite{Verhelst2021}.\\

\subsubsection{Limit-Point Instabilities}
Provided the geometry of a snapping meta-material inspired by \cite{Rafsanjani2015}, see \cref{fig:APALM_snapping_case}. To model the response of this metamaterial, different strategies can be used. One can employ displacement-controlled simulation by incrementally fixing the displacement on the top boundary and computing the corresponding deformations in the domain. Alternatively, one can use arc-length methods to find the snapping behaviour of the metamaterial. In the latter case, the Adaptive Parallel Arc-Length Method can be used to speed up computations, as demonstrated in \cite{Verhelst2023APALM}. In the following, the setup of a displacement-controlled simulation is demonstrated, followed by the setup of an arc-length method, both using the \gs{gsStructuralAnalysis} module.\\

\begin{figure}
    \centering
        \begin{subfigure}{\halfwidth}
            \centering
            \resizebox{\linewidth}{!}
            {
            \def\Nx{3}
            \def\Ny{2}
            \def\al{0.2}
            \def\tw{0.15}
            \def\tb{0.15}
            \def\tg{0.1}
            \def\ts{0.1}
            \def\l{1}
            \def\samples{50}
            \begin{tikzpicture}[scale=1]
            \pgfmathsetmacro\L{\Nx*\l}
            \pgfmathsetmacro\Nym{\Ny-1}
            \pgfmathsetmacro\Nxx{2*\Nx-1}
            \foreach \y in {0,...,\Ny}
            {
                \pgfmathsetmacro\hm{\y*(\ts+2*\tg+\tb)+\al*\l/2-\tg/2}
                \pgfmathsetmacro\hp{\y*(\ts+2*\tg+\tb)+\al*\l/2+1.5*\tg}
                \foreach \x in {1,...,\Nxx}
                {
                    \pgfmathsetmacro\lrm{\x*\l/2-\tw/2}
                    \pgfmathsetmacro\lrp{\x*\l/2+\tw/2}
                    \filldraw[color=gray2] (\lrm,\hm) rectangle (\lrp,\hp);
                }
                \pgfmathsetmacro\lrm{0}
                \pgfmathsetmacro\lrp{\tw/2}
                \filldraw[color=gray2] (\lrm,\hm) rectangle (\lrp,\hp);
                \pgfmathsetmacro\lrm{\L-\tw/2}
                \pgfmathsetmacro\lrp{\L}
                \filldraw[color=gray2] (\lrm,\hm) rectangle (\lrp,\hp);
            }

            \filldraw[color=gray1,domain=0:\L,samples=\samples] plot (\x,{-\al*\l/2*cos(\x/\l*360)})--(\L,-1) -- (0,-1) --cycle;
            \draw[black,thick] (\L,-1) -- (0,-1) node[below,midway]{$\Gamma_1$};
            \draw[black,thick,latex-latex] ($(\L,-1)+(0,0.1)$) -- ($(0,-1)+(0,0.1)$) node[above,midway]{$W$};
            \draw[black,thick,latex-latex] ($(0,-1)-(0.1,0)$) -- ($(0,0)-(0.1,0)$) node[left,midway]{$h_B$};

            \foreach \k in {0,...,\Nym}
            {
                \pgfmathsetmacro\h{\k*(\ts+\tg+\tb)+(\k+1)*\tg}
                \filldraw[color=col2,domain=0:\L,samples=\samples] plot (\x,{-\al*\l/2*cos(\x/\l*360)+\h}) -- plot[domain=\L:0] (\x,{-\al*\l/2*cos(\x/\l*360)+\h+\ts}) --cycle;
                \pgfmathsetmacro\h{\h+\ts+\tg}
                \filldraw[color=col1,domain=0:\L,samples=\samples] plot (\x,{-\al*\l/2*cos(\x/\l*360)+\h}) -- plot[domain=\L:0] (\x,{-\al*\l/2*cos(\x/\l*360)+\h+\tb}) --cycle;
            }

            \pgfmathsetmacro\h{\Ny*(\ts+\tg+\tb)+(\Ny+1)*\tg}
            \filldraw[color=col2,domain=0:\L,samples=\samples] plot (\x,{-\al*\l/2*cos(\x/\l*360)+\h}) -- plot[domain=\L:0] (\x,{-\al*\l/2*cos(\x/\l*360)+\h+\ts}) --cycle;
            \pgfmathsetmacro\H{(\Ny+1)*(2*\tg+\ts+\tb)-\tb}
            \filldraw[color=gray1,domain=0:\L,samples=\samples] plot (\x,{-\al*\l/2*cos(\x/\l*360)+\H})--(\L,\H+1) -- (0,\H+1) --cycle;
            \draw[black,thick] (\L,\H+1) -- (0,\H+1) node[below,midway]{$\Gamma_2$};
            \draw[black,thick,-latex] (\L/2,\H+1) -- (\L/2,\H+1+0.5) node[above]{$\lambda P$};
            \draw[black,thick,latex-latex]($(0,\H)-(0.1,0)$) -- ($(0,\H+1)-(0.1,0)$) node[left,midway]{$h_T$};
            \draw[black,thick,latex-latex] ($(\L,-1)+(0.1,0)$) --  ($(\L,\H+1)+(0.1,0)$)  node[right,midway]{$H$};

            \pgfmathsetmacro\h{\tg}
            \pgfmathsetmacro\hm{\al*\l/2+\tg/2}
            \pgfmathsetmacro\hmg{\al*\l/2+\ts+2*\tg+\tb+\tg/2}
            \draw[color=black,thick] (\l/2,\hm) -- (\l/2+\tw/2,\hm) -- plot[domain=\l/2+\tw/2:3*\l/2-\tw/2,samples=\samples] (\x,{-\al*\l/2*cos(\x/\l*360)+\h}) -- (3*\l/2-\tw/2,\hm) -- (3*\l/2,\hm) -- plot[domain=3*\l/2:\l+\tw/2] (\x,{-\al*\l/2*cos(\x/\l*360)+\h+\ts}) -- plot[domain=\l+\tw/2:3*\l/2] (\x,{-\al*\l/2*cos(\x/\l*360)+\h+\ts+\tg}) -- (3*\l/2,\hmg) -- (3*\l/2-\tw/2,\hmg) -- plot[domain=3*\l/2-\tw/2:\l/2+\tw/2] (\x,{-\al*\l/2*cos(\x/\l*360)+\h+\ts+\tg+\tb}) -- (\l/2+\tw/2,\hmg) -- (\l/2,\hmg)  -- plot[domain=\l/2:\l-\tw/2] (\x,{-\al*\l/2*cos(\x/\l*360)+\h+\ts+\tg})  -- plot[domain=\l-\tw/2:\l/2] (\x,{-\al*\l/2*cos(\x/\l*360)+\h+\ts}) --cycle;

            \end{tikzpicture}
            }
            \caption{A snapping meta-material with $3\times2.5$ building blocks, of which one is outlined. The total multi-patch consists of 132 patches.}
            \label{fig:APALM_snapping_meta-material}
        \end{subfigure}
        \hfill
        \begin{subfigure}{\halfwidth}
            \centering
            \resizebox{\linewidth}{!}
            {
            \def\Nx{3}
            \def\Ny{3}
            \def\al{0.2}
            \def\tw{0.15}
            \def\tb{0.15}
            \def\tg{0.1}
            \def\ts{0.1}
            \def\l{1}
            \def\samples{50}
            \begin{tikzpicture}[scale=4]
                \node (A) at (0,\tg/2){};
                \node (B) at (0,\tg){};
                \node (C) at (\tw/2,\tg/2){};
                \pgfmathsetmacro\x{\tw/2}
                \pgfmathsetmacro\y{\al*\l/2*cos(\x/\l*360)}
                \node (D) at (\x,\y){};
                \node (E) at ($(B)+(0,\ts)$){};
                \node (F) at ($(D)+(0,\ts)$){};
                \pgfmathsetmacro\x{\l/2-\tw/2}
                \pgfmathsetmacro\y{\al*\l/2*cos(\x/\l*360)}
                \node (G) at (\x,\y){};
                \node (H) at ($(G)+(0,\ts)$){};
                \pgfmathsetmacro\x{\l/2+\tw/2}
                \pgfmathsetmacro\y{\al*\l/2*cos(\x/\l*360)}
                \node (I) at (\x,\y){};
                \node (J) at ($(I)+(0,\ts)$){};
                \pgfmathsetmacro\x{\l-\tw/2}
                \node (K) at (\x,\tg/2){};
                \pgfmathsetmacro\x{\l-\tw/2}
                \pgfmathsetmacro\y{\al*\l/2*cos(\x/\l*360)}
                \node (L) at (\x,\y){};
                \node (M) at (\l,\tg/2){};
                \node (N) at (\l,\tg){};
                \node (O) at ($(L)+(0,\ts)$){};
                \node (P) at ($(N)+(0,\ts)$){};

                \node (Q) at ($(E)+(0,\tg)$){};
                \node (R) at ($(F)+(0,\tg)$){};
                \node (S) at ($(H)+(0,\tg)$){};
                \node (T) at ($(J)+(0,\tg)$){};
                \node (U) at ($(O)+(0,\tg)$){};
                \node (V) at ($(P)+(0,\tg)$){};

                \node (W) at ($(Q)+(0,\tb)$){};
                \node (X) at ($(R)+(0,\tb)$){};
                \node (Y) at ($(S)+(0,\tb)$){};
                \node (Z) at ($(T)+(0,\tb)$){};
                \node (AA) at ($(U)+(0,\tb)$){};
                \node (AB) at ($(V)+(0,\tb)$){};

                \node (AC) at (0,\tg+\tb+\ts+\al*\l/2+\tg/2){};
                \node (AD) at (\tw/2,\tg+\tb+\ts+\al*\l/2+\tg/2){};
                \node (AE) at (\x,\tg+\tb+\ts+\al*\l/2+\tg/2){};
                \node (AF) at (\l,\tg+\tb+\ts+\al*\l/2+\tg/2){};

                \draw[fill=gray2]  plot[domain=0:\tw/2,samples=\samples] (\x,{\al*\l/2*cos(\x/\l*360)}) -- (C.center) -- (A.center) -- cycle;
                \draw[fill=gray2]  plot[domain=\l-\tw/2:\l,samples=\samples] (\x,{\al*\l/2*cos(\x/\l*360)}) -- (M.center) -- (K.center) -- cycle;
                \draw[fill=col2]  plot[domain=0:\tw/2,samples=\samples] (\x,{\al*\l/2*cos(\x/\l*360)}) -- plot[domain=\tw/2:0,samples=\samples] (\x,{\al*\l/2*cos(\x/\l*360)+\ts}) -- cycle;
                \draw[fill=col2]  plot[domain=\tw/2:\l/2-\tw/2,samples=\samples] (\x,{\al*\l/2*cos(\x/\l*360)}) -- plot[domain=\l/2-\tw/2:\tw/2,samples=\samples] (\x,{\al*\l/2*cos(\x/\l*360)+\ts}) -- cycle;
                \draw[fill=col2]  plot[domain=\l/2-\tw/2:\l/2+\tw/2,samples=\samples] (\x,{\al*\l/2*cos(\x/\l*360)}) -- plot[domain=\l/2+\tw/2:\l/2-\tw/2,samples=\samples] (\x,{\al*\l/2*cos(\x/\l*360)+\ts}) -- cycle;
                \draw[fill=col2]  plot[domain=\l/2+\tw/2:\l-\tw/2,samples=\samples] (\x,{\al*\l/2*cos(\x/\l*360)}) -- plot[domain=\l-\tw/2:\l/2+\tw/2,samples=\samples] (\x,{\al*\l/2*cos(\x/\l*360)+\ts}) -- cycle;
                \draw[fill=col2]  plot[domain=\l-\tw/2:\l,samples=\samples] (\x,{\al*\l/2*cos(\x/\l*360)}) -- plot[domain=\l:\l-\tw/2,samples=\samples] (\x,{\al*\l/2*cos(\x/\l*360)+\ts}) -- cycle;

                \draw[fill=gray2]  plot[domain=\l/2-\tw/2:\l/2+\tw/2,samples=\samples] (\x,{\al*\l/2*cos(\x/\l*360)+\tg}) -- plot[domain=\l/2+\tw/2:\l/2-\tw/2,samples=\samples] (\x,{\al*\l/2*cos(\x/\l*360)+\ts+\tg}) -- cycle;

                \draw[fill=col1]  plot[domain=0:\tw/2,samples=\samples] (\x,{\al*\l/2*cos(\x/\l*360)+\ts+\tg}) -- plot[domain=\tw/2:0,samples=\samples] (\x,{\al*\l/2*cos(\x/\l*360)+\ts+\tg+\tb}) -- cycle;
                \draw[fill=col1]  plot[domain=\tw/2:\l/2-\tw/2,samples=\samples] (\x,{\al*\l/2*cos(\x/\l*360)+\ts+\tg}) -- plot[domain=\l/2-\tw/2:\tw/2,samples=\samples] (\x,{\al*\l/2*cos(\x/\l*360)+\ts+\tg+\tb}) -- cycle;
                \draw[fill=col1]  plot[domain=\l/2-\tw/2:\l/2+\tw/2,samples=\samples] (\x,{\al*\l/2*cos(\x/\l*360)+\ts+\tg}) -- plot[domain=\l/2+\tw/2:\l/2-\tw/2,samples=\samples] (\x,{\al*\l/2*cos(\x/\l*360)+\ts+\tg+\tb}) -- cycle;
                \draw[fill=col1]  plot[domain=\l/2+\tw/2:\l-\tw/2,samples=\samples] (\x,{\al*\l/2*cos(\x/\l*360)+\ts+\tg}) -- plot[domain=\l-\tw/2:\l/2+\tw/2,samples=\samples] (\x,{\al*\l/2*cos(\x/\l*360)+\ts+\tg+\tb}) -- cycle;
                \draw[fill=col1]  plot[domain=\l-\tw/2:\l,samples=\samples] (\x,{\al*\l/2*cos(\x/\l*360)+\ts+\tg}) -- plot[domain=\l:\l-\tw/2,samples=\samples] (\x,{\al*\l/2*cos(\x/\l*360)+\ts+\tg+\tb}) -- cycle;

                \draw[fill=gray2]  plot[domain=0:\tw/2,samples=\samples] (\x,{\al*\l/2*cos(\x/\l*360)+\ts+\tg+\tb}) -- (AD.center) -- (AC.center) -- cycle;
                \draw[fill=gray2]  plot[domain=\l-\tw/2:\l,samples=\samples] (\x,{\al*\l/2*cos(\x/\l*360)+\ts+\tg+\tb}) -- (AF.center) -- (AE.center) -- cycle;

                \path (H)--(S) node[midway] (HS){};
                \path (J)--(T) node[midway] (JT){};
                \draw[latex-latex] (HS.north) -- (JT.north) node[midway,below]{$t_w$};
                \draw[latex-latex] ($(N)+(\tw/4,0)$) -- ($(M)+(\tw/4,0)$) node[midway,right]{$\nicefrac{t_g}{2}$};
                \draw[latex-latex] ($(P)+(\tw/4,0)$) -- ($(N)+(\tw/4,0)$) node[midway,right]{$t_s$};
                \draw[latex-latex] ($(V)+(\tw/4,0)$) -- ($(P)+(\tw/4,0)$) node[midway,right]{$t_g$};
                \draw[latex-latex] ($(AB)+(\tw/4,0)$) -- ($(V)+(\tw/4,0)$) node[midway,right]{$t_b$};
                \draw[latex-latex] ($(AF)+(\tw/4,0)$) -- ($(AB)+(\tw/4,0)$) node[midway,right]{$\nicefrac{t_g}{2}$};

                \pgfmathsetmacro\y{-\al*\l/2}
                \draw[densely dotted] (-\tw/8,\y) -- (\l/2,\y);
                \draw[densely dotted] ($(B)-(\tw/8,0)$)  -- (B.center);
                \draw[latex-latex] (-\tw/8,\y)  -- ($(B)-(\tw/8,0)$) node[midway,left]{$a$};
                \draw[latex-latex] ($(A)-(\tw/4,0)$) -- ($(AC)-(\tw/4,0)$) node[midway,left]{$h$};

                \draw[latex-latex] (0,-\al*\l/2-\tg/2) -- (\l,-\al*\l/2-\tg/2) node[midway,below]{$l$};
            \end{tikzpicture}
            }
            \caption{The snapping building block, composed of 15 patches outlines in black.}
            \label{fig:APALM_snapping_element}
        \end{subfigure}
    \caption{The problem definition for the snapping meta material using a grid of $3\times2.5$ elements (\subref{fig:APALM_snapping_meta-material}) with the element geometry as defined in (\subref{fig:APALM_snapping_element}). The element dimensions are defined using the thickness of the load-bearing part $t_b=1.5\:[\text{mm}]$ and the thickness of the snapping part $t_s=1.0\:[\text{mm}]$, the thickness of the gap $t_g=1.0\:[\text{mm}]$ and the thickness of the connectors $t_w=1.5\:[\text{mm}]$, such that the height $h=t_b+t_s+2t_g$. The length of the element is $\ell=10\:[\text{mm}]$, and the amplitude of the cosine wave defining the element shape is given by $a=0.3l$. Since the meta-material has $3\times2.5$ elements, the total width is $W=3\ell$. The height of the total metamaterial is given by $H=3h+2t_g+t_s+h_B+h_T$, where $h_B=h_T=5t_g$ are the buffer zones on the top and the bottom. The thickness of the specimen (in out-of-plane direction) is $b=3\:[\text{mm}]$. The material is defined using a compressible Neo-Hookean material model with Young's modulus $E=78\:[\text{N}/\text{mm}^2]$ and Poisson ratio $\nu=0.4\:[-]$. The bottom boundary $\Gamma_1$ is fixed using $u_x=u_y=0$, and the top boundary $\Gamma_2$ is fixed in the horizontal direction ($u_x=0$) and coupled in the vertical direction $u_y$. The load applied on the top boundary is a variable defined by $\lambda P$.}
    \label{fig:APALM_snapping_case}
\end{figure}

Firstly, a displacement-controlled simulation can be done by incrementally increasing the Dirichlet boundary condition on $\Gamma_2$ and performing a static solve in each load step. This routine is simplified by the \gs{gsControlDisplacement} class, taking a static solver as an input:

\ResetLineNumber
\begin{lstlisting}
// Provided K, F, Jacobian, ALResidual, assembler, topBdr, b, W, H,
// probePatch, probePoint, Emax
gsStaticNewton<real_t> staticNR(K,F,Jacobian,ALResidual);
gsControlDisplacement<real_t> control(&staticNR);
real_t sigma, eps;
gsMultiPatch<real_t> deformed, displacement;
while (eps<=Emax)
{
    // Perform a step
    control.step(dy);
    // Construct the deformed geometry
    assembler.constructSolution(U,deformed);
    // Compute the equivalent stress on the top boundary
    sigma = assembler.boundaryForce(deformed,topBdr)[1] / (b*W);
    // Construct the displacement field
    assembler.constructDisplacement(solVector,displacement);
    // Compute epsilon based on the displacement at probePatch and probePoint;
    eps = displacement.patch(probePatch).eval(probePoint) / H;
}
\end{lstlisting}
Here, the \gs{probePoint} defined on a \gs{probePatch} is used to evaluate the displacement of the top boundary for the computation of the total strain $\varepsilon$ of the metamaterial.\\

Secondly, the recently developed Adaptive Parallel Arc-Length Method (APALM) \cite{Verhelst2023APALM} allows to run quasi-static simulations in parallel and is implemented in \gs{gsAPALM}. Since the APALM can be implemented in a non-intrusive way and therefore can be simply wrapped around a \gs{gsALMBase} class. To keep track of the solutions, the \gs{gsAPALMData} accompanies the \gs{gsAPALM} class. In \gismo, the APALM is used as follows:

\begin{lstlisting}
// Provided arcLength
gsAPALMData<real_t,solution_t> apalmData;
const gsMpi & MPI               = gsMpi::init(argc, argv);
gsMpiComm              MPI_COMM = MPI.worldComm();
gsAPALM<real_t> apalm(arcLength,apalmData,MPI_COMM);
\end{lstlisting}
Here, the second template argument of the class \gs{gsAPALMData} is the way solutions are provided. Here, \gs{solution_t} is \gs{std::pair<gsVector<real_t>,real_t>}, where the first vector is the discrete vector of solutions and the second is a scalar representing the load. Furthermore, it can be seen that the APALM algorithms are decoupled from the APALM data structure. The reason for this is that the data structure in \gs{gsAPALMData} could be used for other applications beyond arc-length methods. The maximum refinement level of the APALM routine and the tolerance are specified in the options of the data structure \gs{gsAPALMData}, since the data structure decides whether an interval should be refined or not. The number of sub-intervals per refinement is specified in the \gs{gsAPALM} class, since this is the class submitting solutions into the data structure. When the APALM is defined, the function \gs{solve(N)} can be called to perform fully parallel load-stepping using the APALM using \gs{N} steps, or the functions \gs{serialSolve(N)} and \gs{parallelSolve(N)} can be called to perform serial initialisation and parallel correction, respectively, in a segregated way (see \cite{Verhelst2023APALM}). The \gs{gsAPALM} class is compatible with singular point detection and approach, as described above.\\

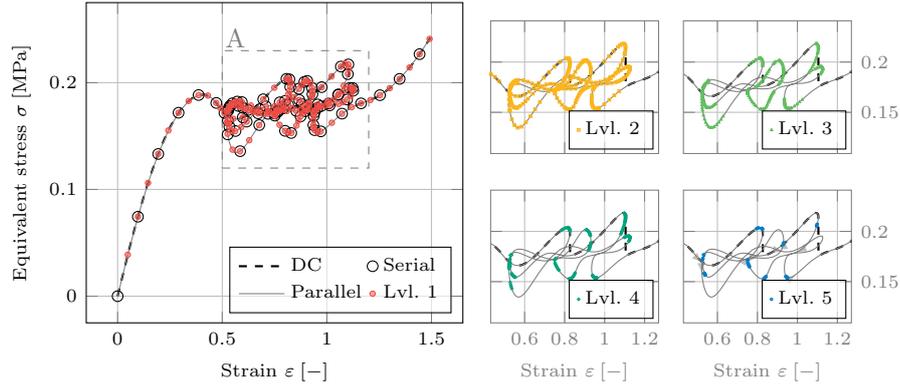
\begin{figure}
    \def\minXi{0.5}
    \def\maxXi{1.2}
    \def\minYi{0.12}
    \def\maxYi{0.23}
    \centering
    \begin{subfigure}[t]{\linewidth}
        \centering
        \begin{tikzpicture}
            \begin{axis}
                [
                xlabel={Strain $\varepsilon\:[-]$},
                ylabel={Equivalent stress $\sigma\:[\text{MPa}]$},
                legend pos = south east,
                xmin = 0,
                xmax = 1.5,
                ymin = 0,
                ymax = 0.25,
                enlarge x limits = true,
                enlarge y limits = true,
                width=0.5\linewidth,
                height=0.3\textheight,
                grid=major,
                legend columns=2,
                restrict x to domain = 0:1.5
                ]
                \addplot+[dashed,no markers,black,thick] table[header=true,x expr = \thisrowno{3},y expr = -\thisrowno{4}/1e6, col sep = comma]{Data/Snapping_DC.csv};
                \addlegendentry{DC}
                \def\level{0}
                \addplot+[only marks,mark=o,mark size=2.0,black] table[header=true,x expr =\thisrowno{3},y expr = \thisrowno{4}/1e6, col sep = comma,restrict expr to domain={\thisrowno{7}}{\level:\level}]{Data/Snapping_parallel.csv};
                \addlegendentry{Serial}
                \addplot+[no markers,gray,thin,solid,
                ] table[header=true,x expr = \thisrowno{3},y expr = \thisrowno{4}/1e6, col sep = comma]{Data/Snapping_parallel.csv};
                \addlegendentry{Parallel}

                \def\level{1}
                \pgfplotsset{cycle list shift=-3}
                    \addplot+[only marks,mark size=1.0] table[header=true,x expr =\thisrowno{3},y expr = \thisrowno{4}/1e6, col sep = comma,restrict expr to domain={\thisrowno{7}}{\level:\level}]{Data/Snapping_parallel.csv};
                    \addlegendentry{Lvl. 1}
                \draw[gray,dashed] (axis cs: \minXi,\minYi) -- (axis cs: \maxXi,\minYi) -- (axis cs: \maxXi,\maxYi) -- (axis cs: \minXi,\maxYi) node[above right, inner sep=1pt] {A} --cycle;
            \end{axis}
        \end{tikzpicture}
        \begin{tikzpicture}
            \begin{groupplot}
                [
                group style={
                    group size=2 by 2,
                    xlabels at=edge bottom,
                    xticklabels at=edge bottom,
                    ylabels at=edge right,
                    yticklabels at=edge right,
                    vertical sep=0.025\textheight,
                    horizontal sep=0.025\linewidth,
                },
                xlabel style={color=gray},
                x tick label style={color = gray},
                ylabel style={color=gray},
                y tick label style={color = gray},
                legend pos = south east,
                width=0.29\linewidth,
                height=0.172\textheight,
                grid=major,
                legend columns=2,
                restrict x to domain = 0:1.5,
                ylabel={\phantom{.}},
                xlabel={Strain $\varepsilon\:[-]$}
                ]
                \nextgroupplot[
                xmin = \minXi,
                xmax = \maxXi,
                ymin = \minYi,
                ymax = \maxYi,
                draw = gray,
                ]
                \addplot+[dashed,no markers,black,thick,forget plot] table[header=true,x expr = \thisrowno{3},y expr = -\thisrowno{4}/1e6, col sep = comma]{Data/Snapping_DC.csv};
                \addplot+[no markers,gray,thin,solid,forget plot] table[header=true,x expr = \thisrowno{3},y expr = \thisrowno{4}/1e6, col sep = comma]{Data/Snapping_parallel.csv};
                \pgfplotsset{cycle list shift=-1}
                \def\level{2}
                    \addplot+[only marks,mark size=0.5] table[header=true,x expr =\thisrowno{3},y expr = \thisrowno{4}/1e6, col sep = comma,restrict expr to domain={\thisrowno{7}}{\level:\level}]{Data/Snapping_parallel.csv};
                    \addlegendentry{Lvl. \level}

                \nextgroupplot[
                xmin = \minXi,
                xmax = \maxXi,
                ymin = \minYi,
                ymax = \maxYi,
                draw = gray,
                ]
                \addplot+[dashed,no markers,black,thick,forget plot] table[header=true,x expr = \thisrowno{3},y expr = -\thisrowno{4}/1e6, col sep = comma]{Data/Snapping_DC.csv};
                \addplot+[no markers,gray,thin,solid,forget plot] table[header=true,x expr = \thisrowno{3},y expr = \thisrowno{4}/1e6, col sep = comma]{Data/Snapping_parallel.csv};
                \pgfplotsset{cycle list shift=2}
                \def\level{3}
                    \addplot+[only marks,mark size=0.5] table[header=true,x expr =\thisrowno{3},y expr = \thisrowno{4}/1e6, col sep = comma,restrict expr to domain={\thisrowno{7}}{\level:\level}]{Data/Snapping_parallel.csv};
                    \addlegendentry{Lvl. \level}

                \nextgroupplot[
                xmin = \minXi,
                xmax = \maxXi,
                ymin = \minYi,
                ymax = \maxYi,
                draw = gray,
                ]
                \addplot+[dashed,no markers,black,thick,forget plot] table[header=true,x expr = \thisrowno{3},y expr = -\thisrowno{4}/1e6, col sep = comma]{Data/Snapping_DC.csv};
                \addplot+[no markers,gray,thin,solid,forget plot] table[header=true,x expr = \thisrowno{3},y expr = \thisrowno{4}/1e6, col sep = comma]{Data/Snapping_parallel.csv};
                \pgfplotsset{cycle list shift=3}
                \def\level{4}
                    \addplot+[only marks,mark size=0.5] table[header=true,x expr =\thisrowno{3},y expr = \thisrowno{4}/1e6, col sep = comma,restrict expr to domain={\thisrowno{7}}{\level:\level}]{Data/Snapping_parallel.csv};
                    \addlegendentry{Lvl. \level}

                \nextgroupplot[
                xmin = \minXi,
                xmax = \maxXi,
                ymin = \minYi,
                ymax = \maxYi,
                draw = gray,
                ]
                \addplot+[dashed,no markers,black,thick,forget plot] table[header=true,x expr = \thisrowno{3},y expr = -\thisrowno{4}/1e6, col sep = comma]{Data/Snapping_DC.csv};
                \addplot+[no markers,gray,thin,solid,forget plot,
                postaction={decorate, decoration={markings,
                        mark=at position 0.2 with {\arrow{latex};},
                        mark=at position 0.3 with {\arrow{latex};},
                        mark=at position 0.4 with {\arrow{latex};},
                        mark=at position 0.5 with {\arrow{latex};},
                        mark=at position 0.6 with {\arrow{latex};},
                        mark=at position 0.7 with {\arrow{latex};},
                        mark=at position 0.8 with {\arrow{latex};},
                        mark=at position 0.9 with {\arrow{latex};},
                }}] table[header=true,x expr = \thisrowno{3},y expr = \thisrowno{4}/1e6, col sep = comma ]{Data/Snapping_parallel.csv};
                \def\level{5}
                \pgfplotsset{cycle list shift=4}
                    \addplot+[only marks,mark size=0.5] table[header=true,x expr =\thisrowno{3},y expr = \thisrowno{4}/1e6, col sep = comma,restrict expr to domain={\thisrowno{7}}{\level:\level}]{Data/Snapping_parallel.csv};
                    \addlegendentry{Lvl. \level}

            \end{groupplot}
        \end{tikzpicture}
    \end{subfigure}
    \caption{Stress-strain diagram for the snapping meta-material from \cref{fig:APALM_snapping_case}. The vertical axis depicts the equivalent stress $\sigma=\lambda P / (b W)$, and the horizontal axis represents the strain $\varepsilon=u_y/H$, where $u_y$ is the displacement of the top boundary $\Gamma_2$. The complete curve with the displacement-controlled (DC) results, the points obtained in serial initialisation, and the line obtained by parallel corrections are presented on the left. The figures on the right present the points from different hierarchical levels at the inset depicted in the left diagram. The simulation is performed with a tolerance of $\varepsilon_l=\varepsilon_u=10^{-3}$ and an increment length of $\Delta L=0.05$.}
    \label{fig:APALM_snapping}
\end{figure}

The equilibrium paths following from the APALM and the displacement-controlled simulation are provided in \cref{fig:APALM_snapping}. For a discussion of the results, the reader is referred to \cite{Verhelst2023APALM}.

\subsection{Error Analysis and Adaptivity}\label{subsec:examples_adaptivity}
Provided the wrinkling problem from \cref{subsubsec:examples_ALM_bifurcation}, an adaptive meshing procedure can be set up. Here, the setup of the dual-weighted residual (DWR) error estimator for Kirchhoff-Love shells is presented, which allows for goal-oriented error estimation, hence goal-adaptive meshing. For more details, the reader is referred to the paper \cite{Verhelst2023Adaptive}.\\

\begin{figure}
    \centering
    \begin{tikzpicture}
        \begin{groupplot}[
            group style={
                group name=my plots,
                group size=1 by 3,
                xlabels at=edge bottom,
                xticklabels at=edge bottom,
                ylabels at=edge left,
                yticklabels at=edge left,
                vertical sep=5pt,
                horizontal sep=5pt
            },
            width=\linewidth,
            height=0.2\textheight,
            tickpos=left,
            grid = major,
            xlabel=$\epsilon$,
            xmin = 0,
            xmax = 0.5,
        ]
        \nextgroupplot[
            ylabel=$\max(\zeta)/t$,
            ymin = 0,
            ymax = 3.5,
        ]
        \addplot+[no marks,black,dotted,mark=diamond ,mark size=1] table[header=true,x index = {0},y index = {1}, col sep = comma]{Data/ArcLengthWrinkling/LSDynaFull.csv};
        \addplot+[no marks,black,solid,mark=+,mark size=1] table[header=true,x index = {0},y index = {1}, col sep = comma]{Data/ArcLengthWrinkling/ANSYS.csv};
        \addplot+[mark=o,mark size=2,black,mark options={black,solid}] table[header=true,x index = {0},y index = {1}, col sep = comma]{Data/ArcLengthWrinkling/PanaitescuExperiment.csv};
        \addplot [name path=upper,draw=none,mark=none,forget plot] table[header=true,x index = {2},y index = {3}, col sep = comma]{Data/ArcLengthWrinkling/PanaitescuExperiment.csv};
        \addplot [name path=lower,draw=none,mark=none,forget plot] table[header=true,x index = {4},y index = {5}, col sep = comma]{Data/ArcLengthWrinkling/PanaitescuExperiment.csv};
        \addplot [fill=black!10,forget plot] fill between[of=upper and lower];

        \pgfplotsset{cycle list shift=-3}
        \addplot+[] table[header=true,x index = {1},y index = {2}, col sep = comma]{Data/ArcLengthWrinkling/Wrinkling_r5_adaptive_R050_C0005.csv};
        \addplot [forget plot,only marks,black,mark indices={9,10,11,12,40,41,42,43},mark=o] table[header=true,x index = {1},y index = {2}, col sep = comma]{Data/ArcLengthWrinkling/Wrinkling_r5_adaptive_R050_C0005.csv};
        \addplot+[no marks,solid] table[header=true,x index = {1},y index = {2}, col sep = comma]{Data/ArcLengthWrinkling/Wrinkling_r5.csv};
        \addplot+[no marks,solid] table[header=true,x index = {1},y index = {2}, col sep = comma]{Data/ArcLengthWrinkling/Wrinkling_r6.csv};

        \nextgroupplot[ylabel={$\Delta\mathbcal{L}$},ymode=log,height=0.4\textheight]

        \addlegendimage{no marks,black,dotted,mark=diamond ,mark size=1}\addlegendentry{LS-DYNA}
        \addlegendimage{no marks,black,solid,mark=+,mark size=1}\addlegendentry{SHELL181}
        \addlegendimage{mark=o,mark size=2,black,mark options={black,solid}}\addlegendentry{Exp. }

        \fill[black!70,opacity=0.5] (axis cs:0,1e-14)-- (axis cs:0.5,1e-14)--(axis cs:0.5,1e-10)--(axis cs:0,1e-10)--cycle;

        \addplot+[] table[header=true,x index = {1},y index = {4}, col sep = comma]{Data/ArcLengthWrinkling/Wrinkling_r5_adaptive_R050_C0005.csv};\addlegendentry{adaptive, $32\times32$}
        \addplot [forget plot,only marks,black,mark indices={9,10,11,12,40,41,42,43},mark=o] table[header=true,x index = {1},y index = {4}, col sep = comma]{Data/ArcLengthWrinkling/Wrinkling_r5_adaptive_R050_C0005.csv};
        \addplot+[,solid] table[header=true,x index = {1},y index = {4}, col sep = comma]{Data/ArcLengthWrinkling/Wrinkling_r5.csv};\addlegendentry{uniform, $32\times 32$}
        \addplot+[,solid] table[header=true,x index = {1},y index = {4}, col sep = comma]{Data/ArcLengthWrinkling/Wrinkling_r6.csv};\addlegendentry{uniform, $64\times 64$}

        \nextgroupplot[ylabel={\# DoFs},ymode=log,ymin=5e2, ymax=5e4]
        \addplot+[] table[header=true,x index = {1},y index = {3}, col sep = comma]{Data/ArcLengthWrinkling/Wrinkling_r5_adaptive_R050_C0005.csv};
        \addplot [forget plot,only marks,black,mark indices={9,10,11,12,40,41,42,43},mark=o,
                    nodes near coords={\coordindex},
                    every node near coord/.append style={anchor=south,font=\tiny},
                    ] table[header=true,x index = {1},y index = {3}, col sep = comma]{Data/ArcLengthWrinkling/Wrinkling_r5_adaptive_R050_C0005.csv};
        \addplot+[no markers,solid] table[header=true,x index = {1},y index = {3}, col sep = comma]{Data/ArcLengthWrinkling/Wrinkling_r5.csv};
        \addplot+[no markers,solid] table[header=true,x index = {1},y index = {3}, col sep = comma]{Data/ArcLengthWrinkling/Wrinkling_r6.csv};
        \end{groupplot}
        \end{tikzpicture}
    \caption{The non-dimensional maximal wrinkling amplitude $\max(\zeta)/t$ (top), the goal functional error $\Delta\mathcal{L}$ (mid), and the number of degrees of freedom of the computational mesh (bottom) with respect to the strain of the sheet $\epsilon$. The markers with a black border represent the points at which the meshes are provided in \cref{fig:adaptiveMeshing_wrinkling_results_meshes}. The coloured lines are the solutions obtained by the present model with a Mooney-Rivlin material model with uniform or adaptive meshes. For the adaptive simulation, the parameters $(\rho_r,\rho_c)=(0.5,0.005)$ are used. The solid line in the top figure is a SHELL181 result obtained using ANSYS; the dotted line in the top figure is a result obtained using the fully integrated shell in LS-DYNA; and the dashed line with markers represents experimental data obtained by \cite{Panaitescu2019} (Ref.).}
    \label{fig:adaptiveMeshing_wrinkling_results}
\end{figure}
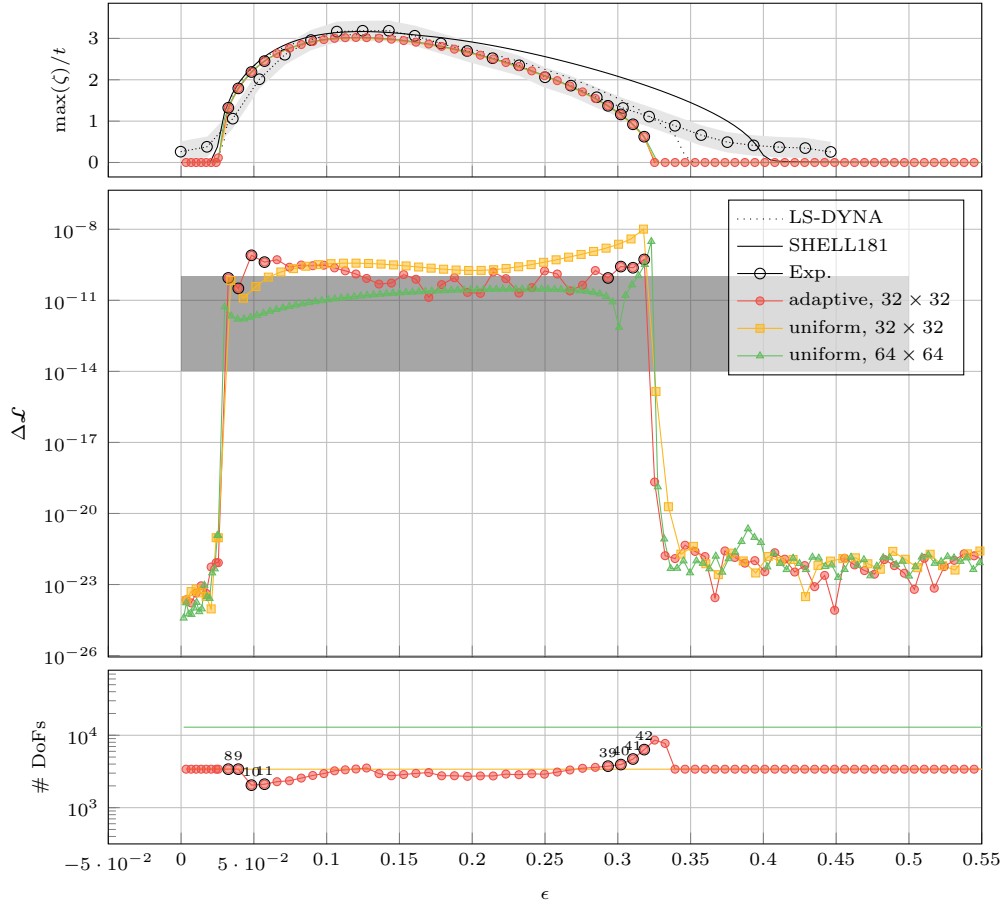

\begin{figure}
    \centering
    \begin{subfigure}{\quartwidth}
    \centering
    \captionsetup{justification=centering}
    \begin{subfigure}{\halfwidth}
    \centering
    \captionsetup{justification=centering}
    \includegraphics[width=\linewidth]{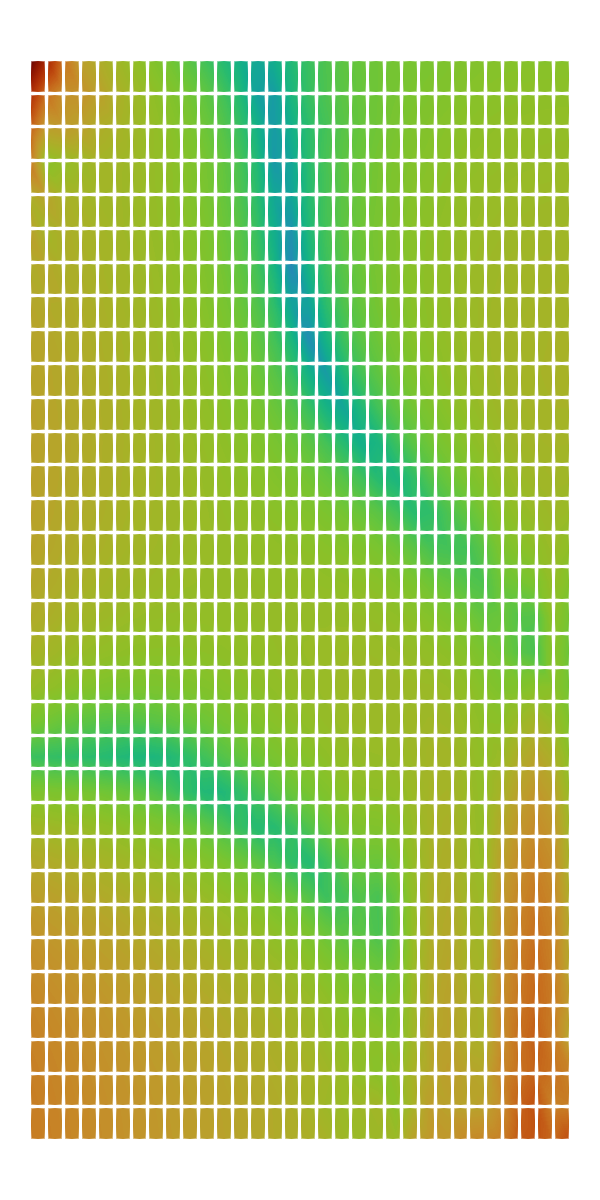}
    \end{subfigure}
    \hfill
    \begin{subfigure}{\halfwidth}
    \centering
    \captionsetup{justification=centering}
    \includegraphics[width=\linewidth]{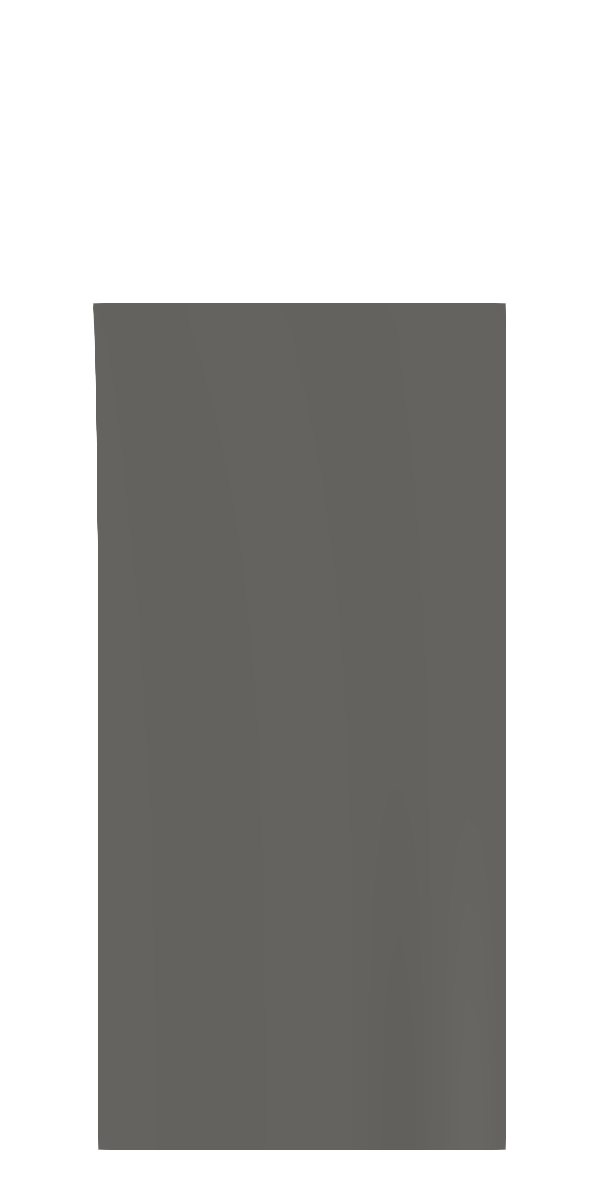}
    \end{subfigure}
    \caption*{Point 8}
    \end{subfigure}
    \hfill
    \begin{subfigure}{\quartwidth}
    \centering
    \captionsetup{justification=centering}
    \begin{subfigure}{\halfwidth}
    \centering
    \captionsetup{justification=centering}
    \includegraphics[width=\linewidth]{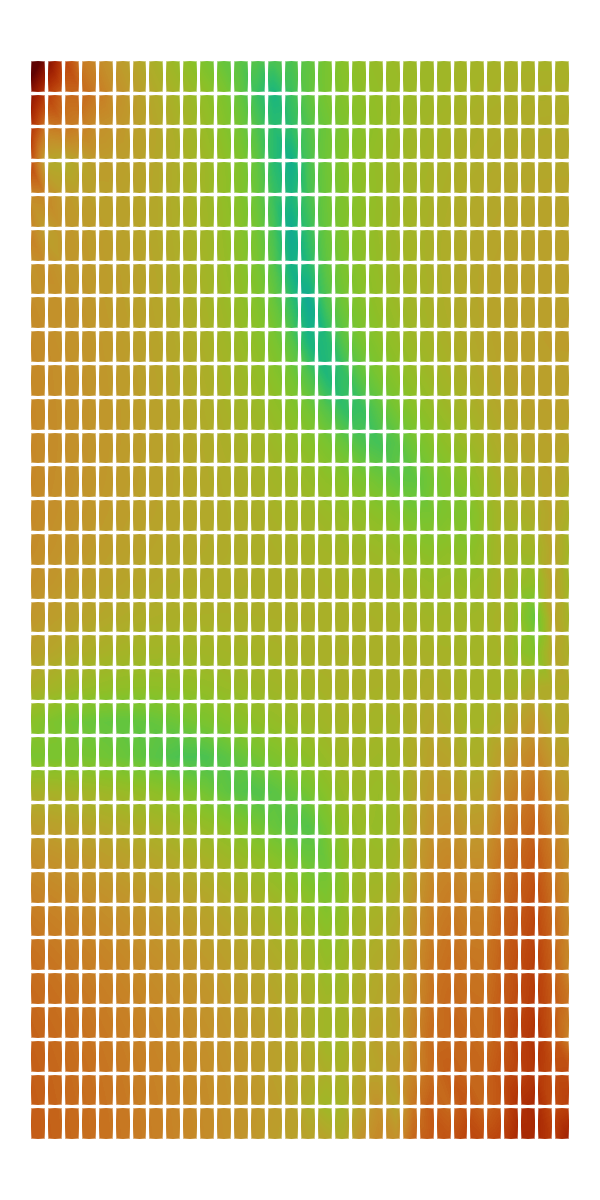}
    \end{subfigure}
    \hfill
    \begin{subfigure}{\halfwidth}
    \centering
    \captionsetup{justification=centering}
    \includegraphics[width=\linewidth]{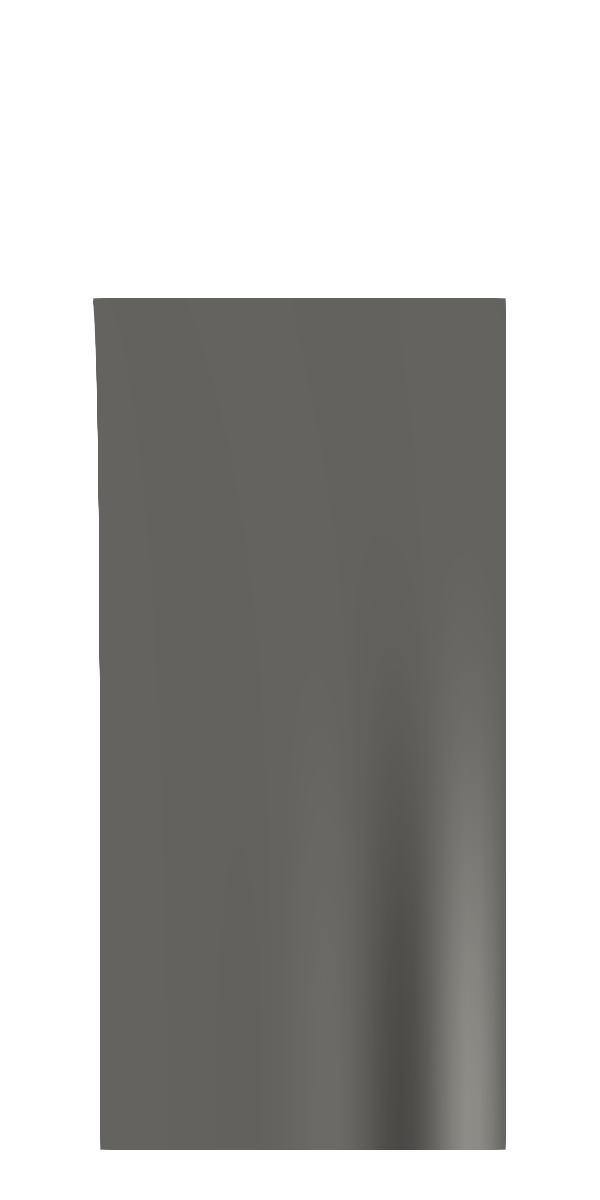}
    \end{subfigure}
    \caption*{Point 9}
    \end{subfigure}
    \hfill
    \begin{subfigure}{\quartwidth}
    \centering
    \captionsetup{justification=centering}
    \begin{subfigure}{\halfwidth}
    \centering
    \captionsetup{justification=centering}
    \includegraphics[width=\linewidth]{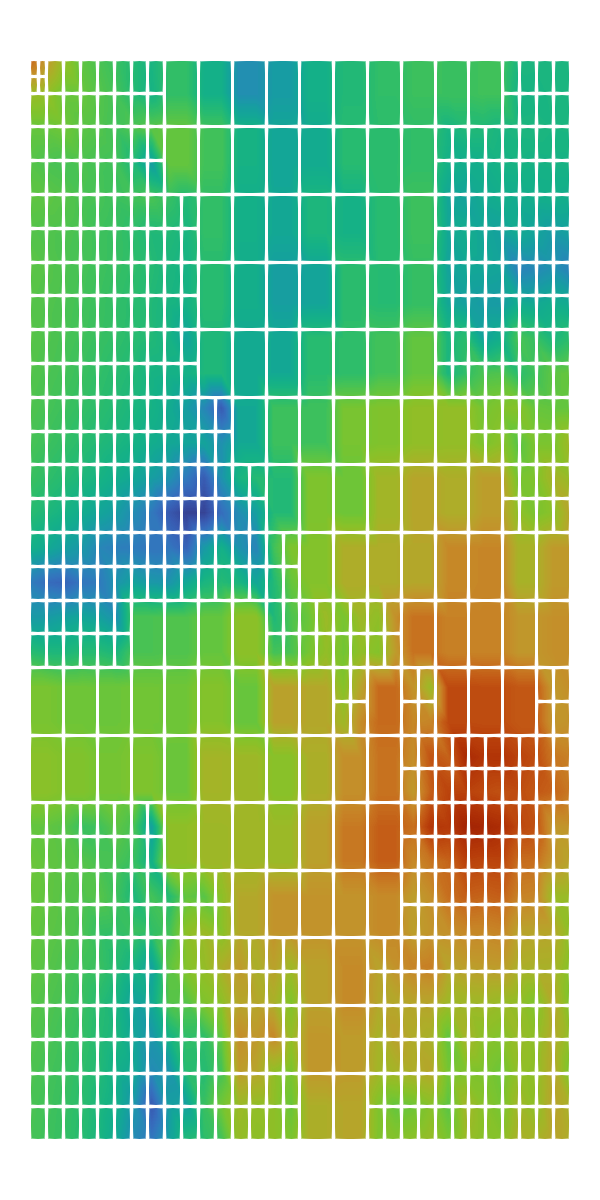}
    \end{subfigure}
    \hfill
    \begin{subfigure}{\halfwidth}
    \centering
    \captionsetup{justification=centering}
    \includegraphics[width=\linewidth]{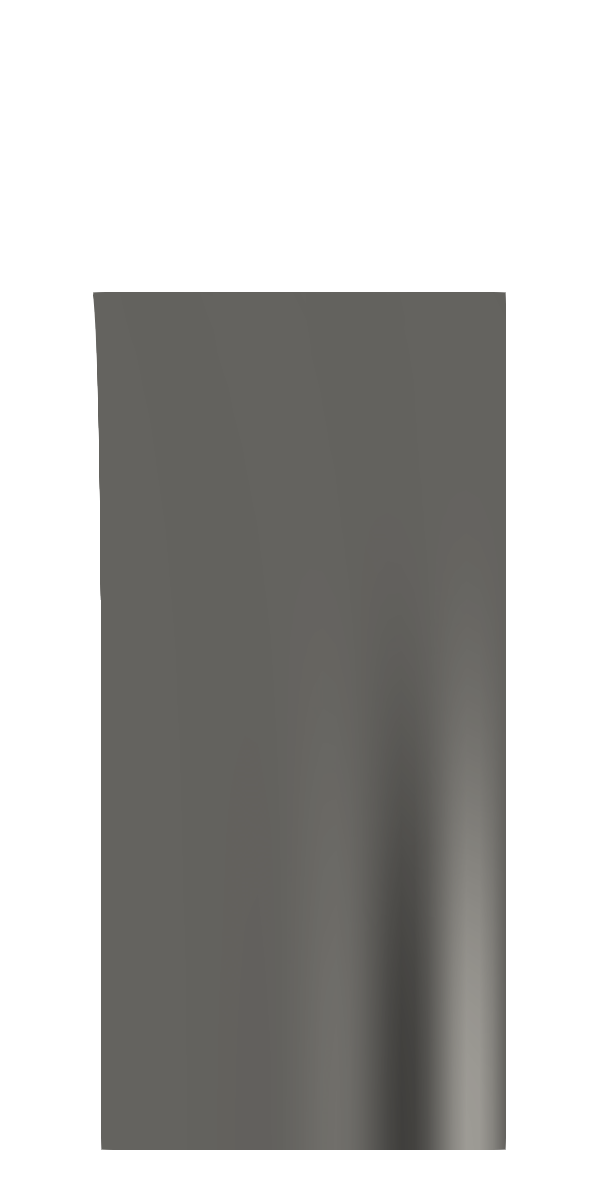}
    \end{subfigure}
    \caption*{Point 10}
    \end{subfigure}
    \hfill
    \begin{subfigure}{\quartwidth}
    \centering
    \captionsetup{justification=centering}
    \begin{subfigure}{\halfwidth}
    \centering
    \captionsetup{justification=centering}
    \includegraphics[width=\linewidth]{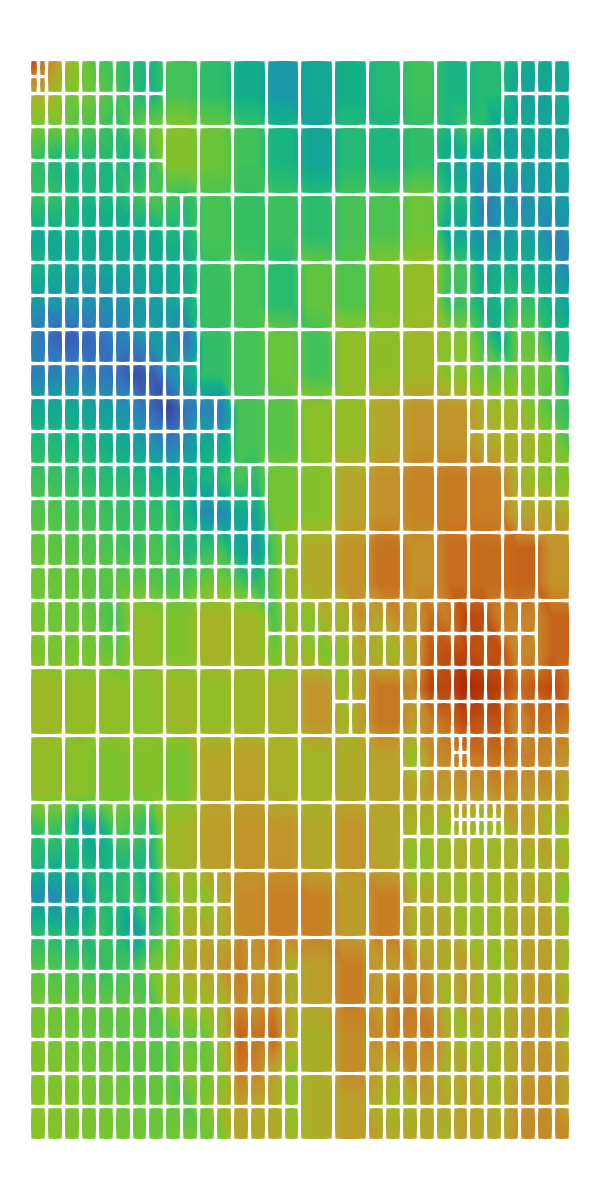}
    \end{subfigure}
    \hfill
    \begin{subfigure}{\halfwidth}
    \centering
    \captionsetup{justification=centering}
    \includegraphics[width=\linewidth]{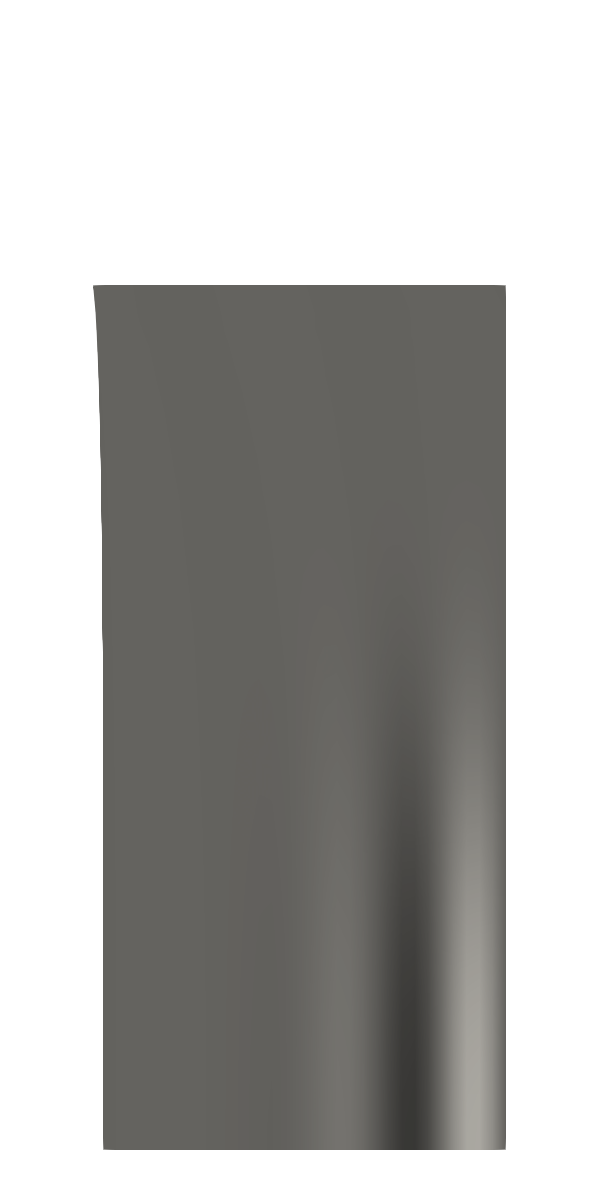}
    \end{subfigure}
    \caption*{Point 11}
    \end{subfigure}

    \begin{subfigure}{\quartwidth}
    \centering
    \captionsetup{justification=centering}
    \begin{subfigure}{\halfwidth}
    \centering
    \captionsetup{justification=centering}
    \includegraphics[width=\linewidth]{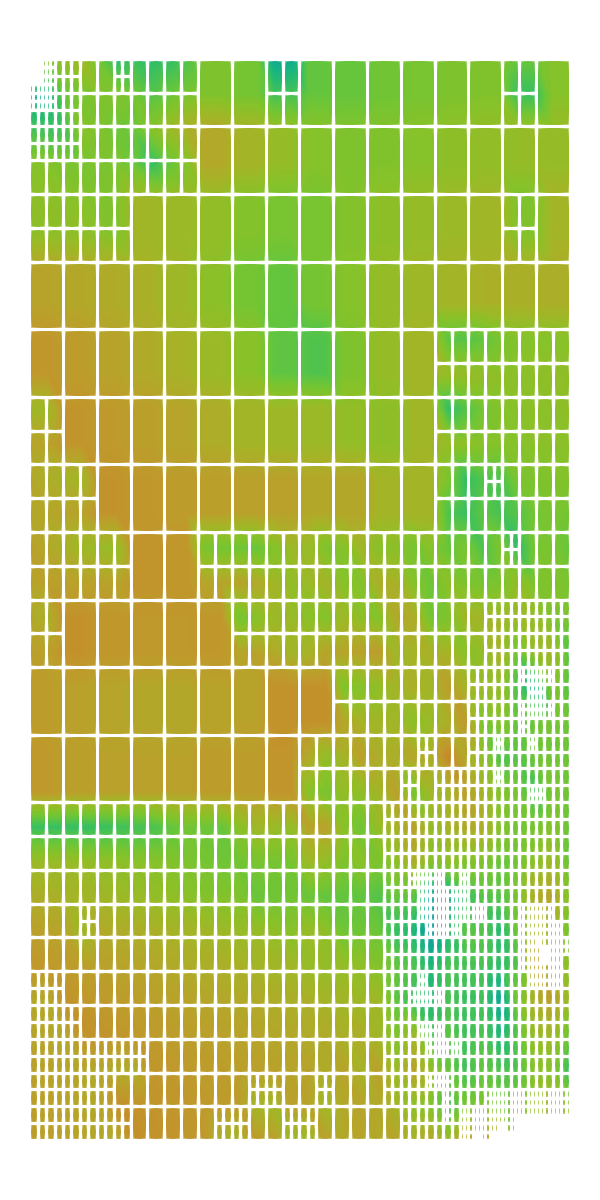}
    \end{subfigure}
    \hfill
    \begin{subfigure}{\halfwidth}
    \centering
    \captionsetup{justification=centering}
    \includegraphics[width=\linewidth]{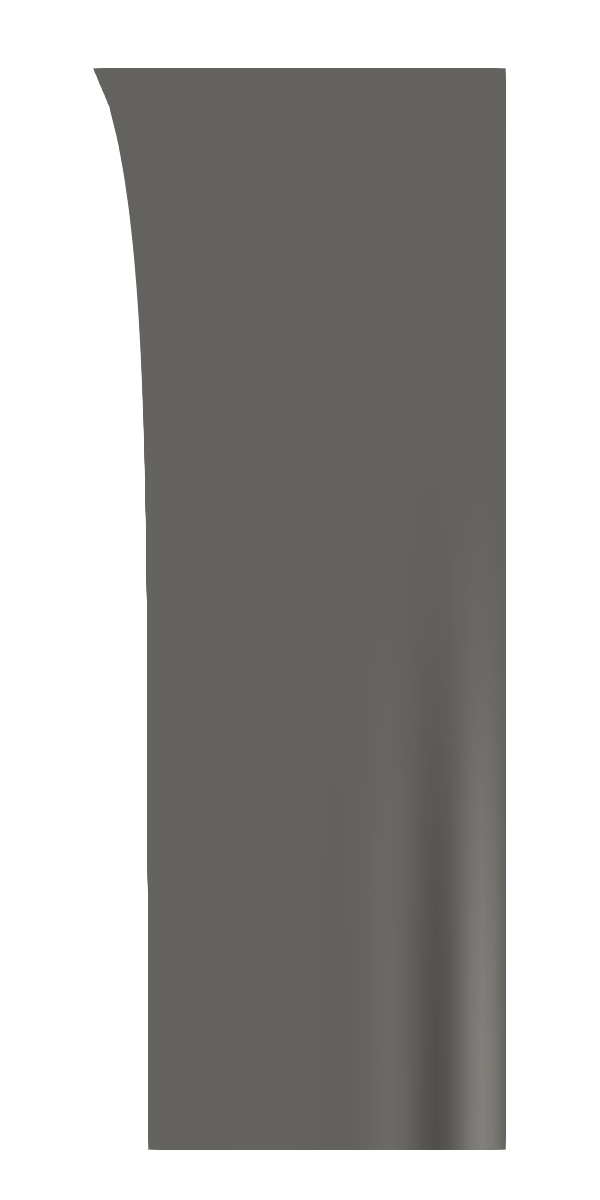}
    \end{subfigure}
    \caption*{Point 39}
    \end{subfigure}
    \hfill
    \begin{subfigure}{\quartwidth}
    \centering
    \captionsetup{justification=centering}
    \begin{subfigure}{\halfwidth}
    \centering
    \captionsetup{justification=centering}
    \includegraphics[width=\linewidth]{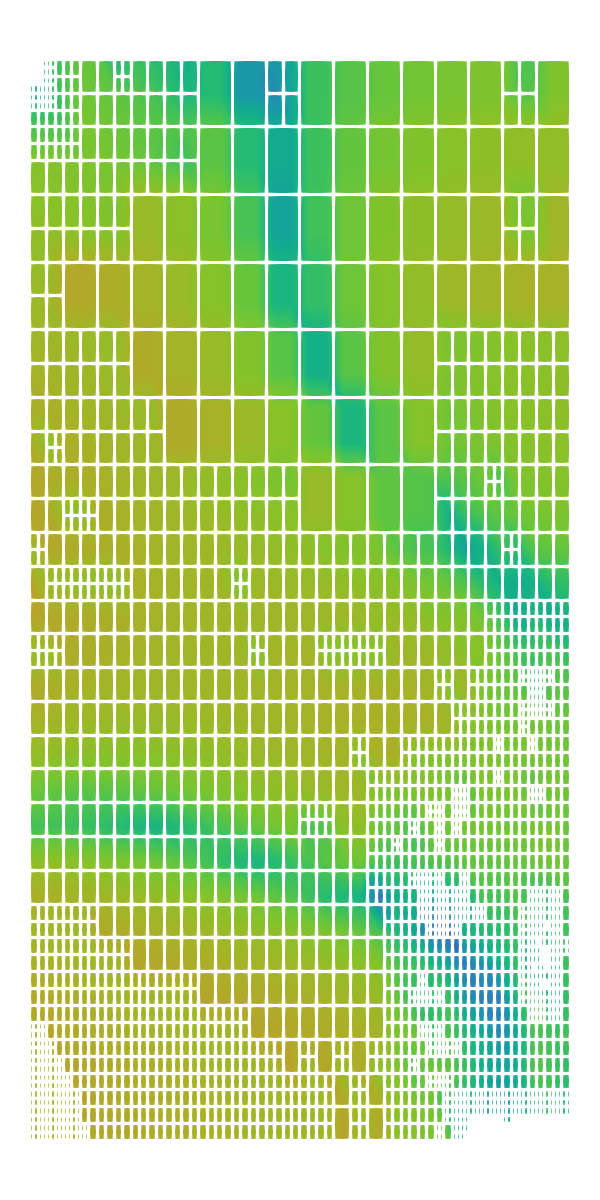}
    \end{subfigure}
    \hfill
    \begin{subfigure}{\halfwidth}
    \centering
    \captionsetup{justification=centering}
    \includegraphics[width=\linewidth]{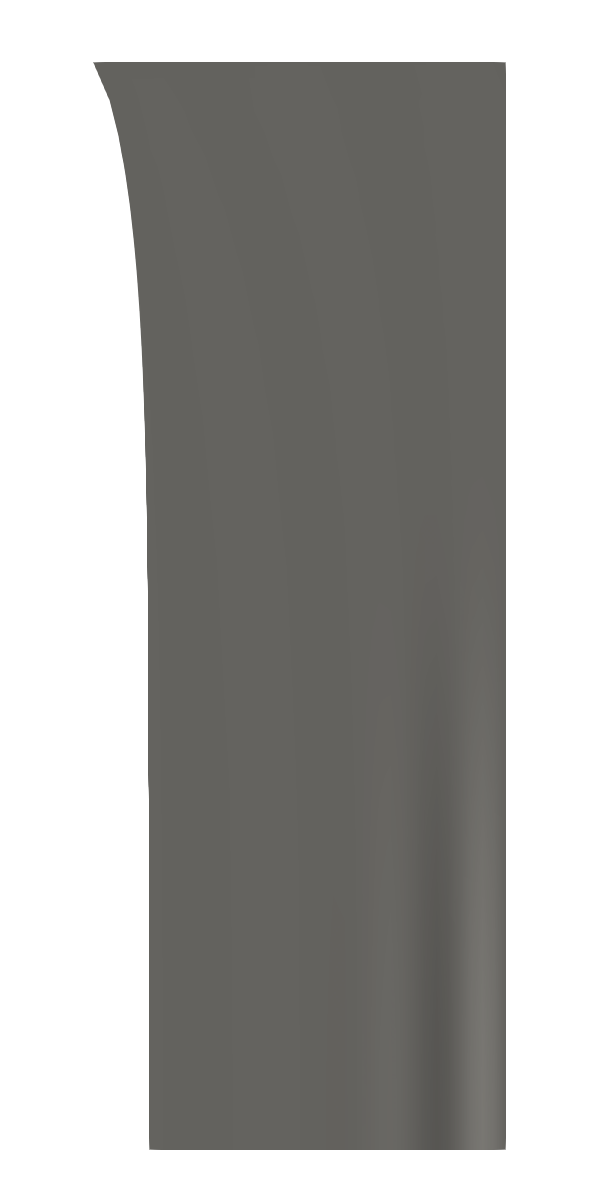}
    \end{subfigure}
    \caption*{Point 40}
    \end{subfigure}
    \hfill
    \begin{subfigure}{\quartwidth}
    \centering
    \captionsetup{justification=centering}
    \begin{subfigure}{\halfwidth}
    \centering
    \captionsetup{justification=centering}
    \includegraphics[width=\linewidth]{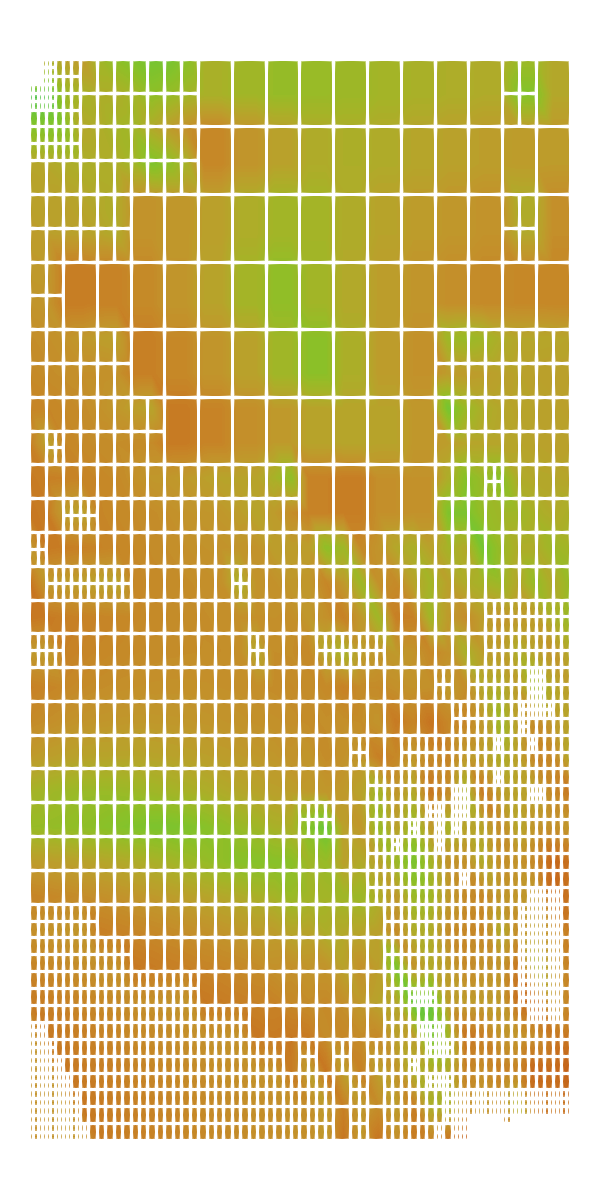}
    \end{subfigure}
    \hfill
    \begin{subfigure}{\halfwidth}
    \centering
    \captionsetup{justification=centering}
    \includegraphics[width=\linewidth]{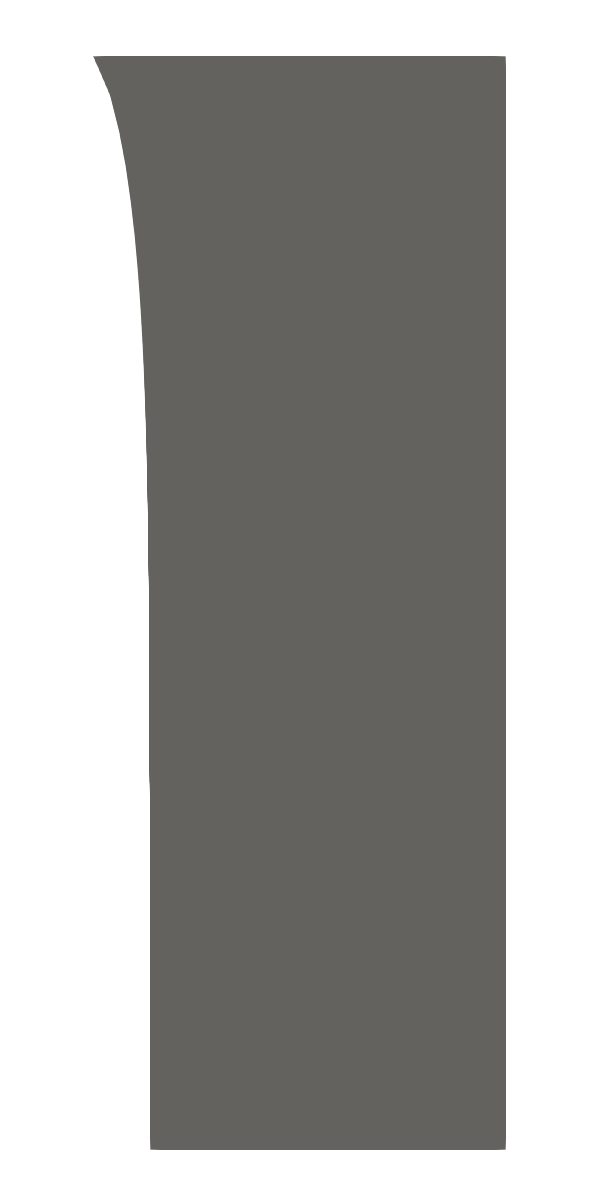}
    \end{subfigure}
    \caption*{Point 41}
    \end{subfigure}
    \hfill
    \begin{subfigure}{\quartwidth}
    \centering
    \captionsetup{justification=centering}
    \begin{subfigure}{\halfwidth}
    \centering
    \captionsetup{justification=centering}
    \includegraphics[width=\linewidth]{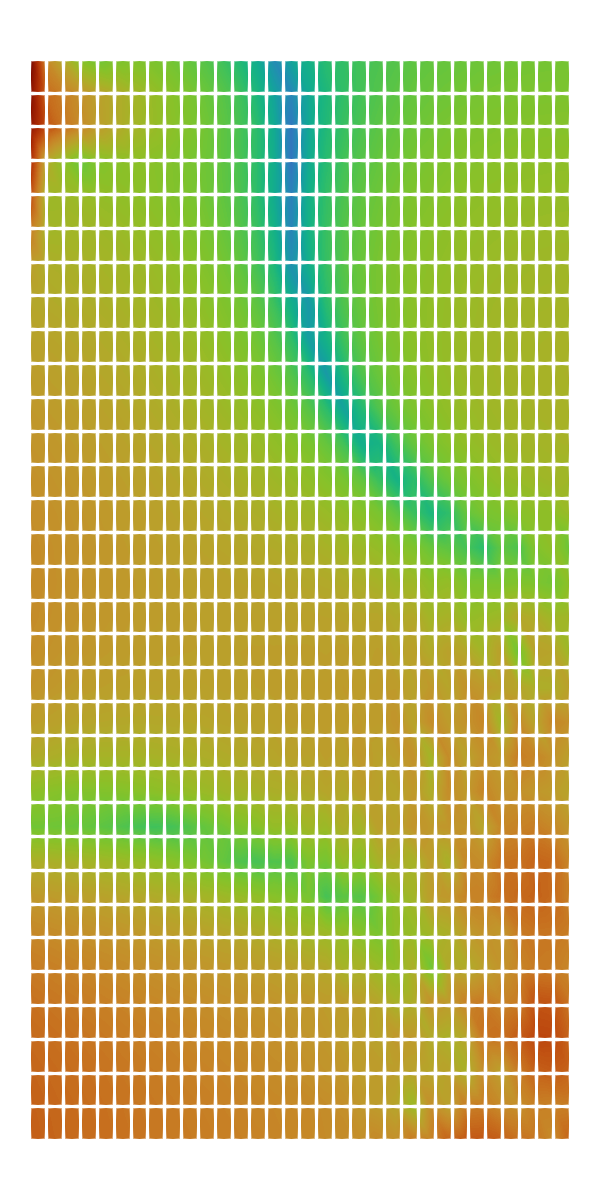}
    \end{subfigure}
    \hfill
    \begin{subfigure}{\halfwidth}
    \centering
    \captionsetup{justification=centering}
    \includegraphics[width=\linewidth]{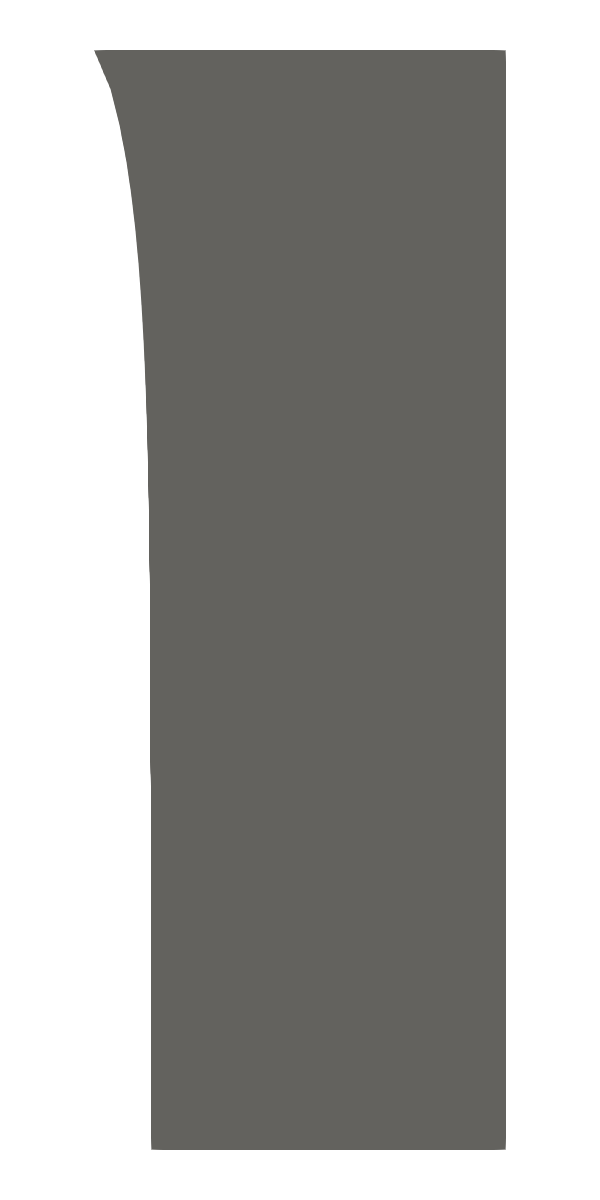}
    \end{subfigure}
    \caption*{Point 42}
    \end{subfigure}

    \begin{subfigure}{\halfwidth}
    \centering
    \captionsetup{justification=centering}
    \includegraphics[width=\twothirdwidth,trim= 0 0 0 1100,clip]{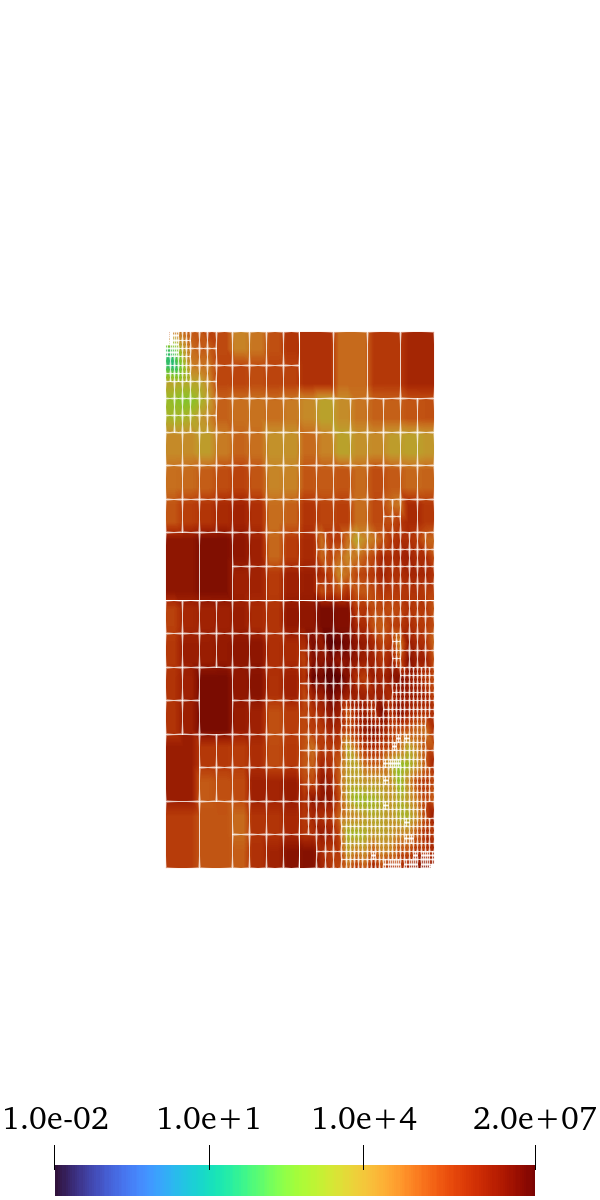}
    \caption*{$e_k/\Delta\mathcal{L}^2$}
    \end{subfigure}
    \begin{subfigure}{\halfwidth}
    \centering
    \captionsetup{justification=centering}
    \includegraphics[width=\twothirdwidth,trim= 0 0 0 1100,clip]{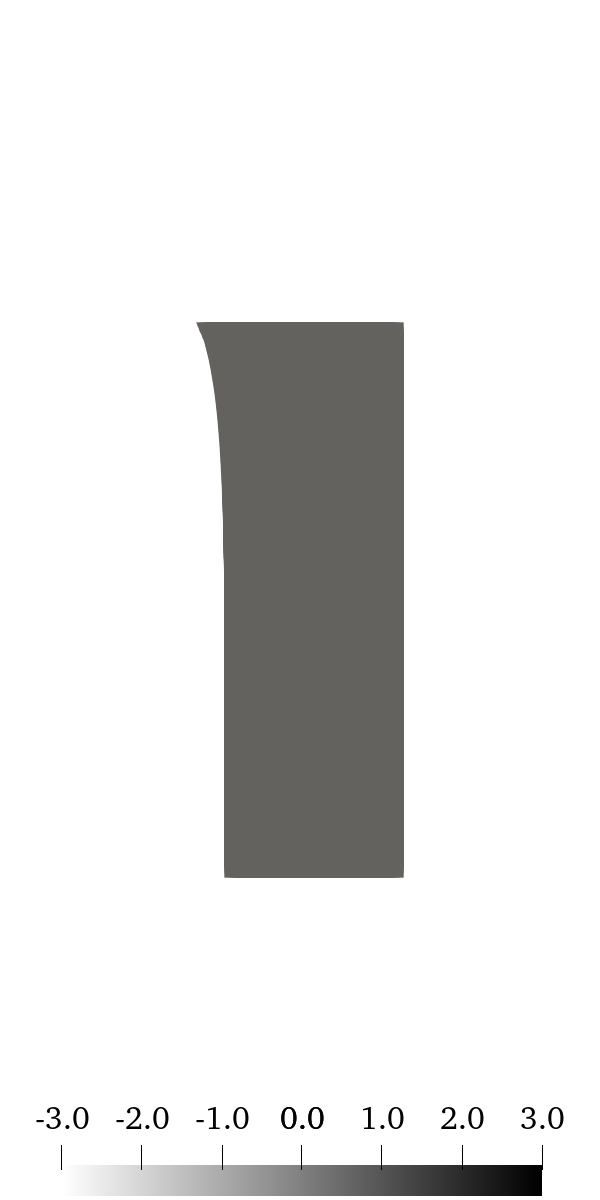}
    \caption*{$u_z/t\:[-]$}
    \end{subfigure}
    \caption{Normalised element errors (left) and the normalised wrinkling amplitude (right) for the wrinkling benchmark plotted on the undeformed geometry with the corresponding elements.}
    \label{fig:adaptiveMeshing_wrinkling_results_meshes}
\end{figure}

The DWR method provided in the \gs{gsKLShell} module provides error estimates in terms of goal functionals such that adaptive meshing can be used. The use of the error estimation routines is implemented in a straight-forward way. The DWR method requires a primal and dual problem to be solved, defined on the basis used for analysis, $\mathbb{S}^p_h$, and an enriched basis, $\tilde{\mathbb{S}}^{p}_h$. Since the error estimator is defined on two spaces, the \gs{gsThinShellAssemblerDWR} class is called the regular \gs{gsThinShellAssemblerDWR}, but with two bases instead of one. The configuration of the \gs{gsThinShellAssemblerDWR} works exactly the same as the configuration of a \gs{gsThinShellAssembler}, but it adds the configuration of a goal functional via the \gs{setGoal}:

\ResetLineNumber
\begin{lstlisting}
// Provided mp, basisL, basisH, bcs, force, materialMatrix
gsThinShellAssemblerDWR<d,T,bool> DWRassembler(mp,basisL,basisH,bcs,force,materialMatrix);
DWRassembler.setGoal(GoalFunction::MembraneStress)
\end{lstlisting}
Here, the \gs{basisL} and \gs{basisH} are the bases $\mathbb{S}^p_h$ and $\tilde{\mathbb{S}}^{p}_h$, respectively, and the other variables are as in \cref{example:gsKLShell}. In order to estimate the error using the DWR, the (non-linear) primal problem has to be solved on $\mathbb{S}^p_h$, and two linear dual problems have to be solved on $\mathbb{S}^p_h$ and $\tilde{\mathbb{S}}^p_h$, respectively. The matrices and vectors for the dual problems are assembled for the provided goal functional using \gs{assembler.assembleDualL()} and can be accessed by \gs{assembler.matrixL()} and \gs{assembler.dualL()} for the matrix and vector, respectively. Since this routine is the same for all problems involving DWR error estimation, the error estimation part has been integrated in the \gs{gsDWRHelper}. Mesh adaptivity can be performed using the \gs{gsAdaptiveMeshing} class, for which the reader is referred to the documentation. Given an arc-length method (see \cref{example:gsStructuralAnalysis}), a \gs{gsDWRHelper}, and an adaptive meshing class, the following loop can be used to adaptively refine a snapping example:

\ContinueLineNumber
\begin{lstlisting}
// Provided DWRassembler, arcLength, mp
gsAdaptiveMeshing<real_t>    mesher(mp);
gsThinShellDWRHelper<real_t> DWRhelper(DWRassembler);
std::vector<real_t>          elErrors;             // element-wise errors
gsVector<real_t>             solutionVector;
gsHBoxContainer<2,real_t>    markRef, markCrs;     // Containers for marked elements
for (index_t k=0; k!=Nsteps; k++)
{
    arcLength.step();                              // Perform an arc-length step
    solutionVector = arcLength.solutionU();        // Obtain the solution vector
    DWRhelper.computeError(solutionVector);        // Compute the errors (element-wise)
    elErrors = DWRhelper.sqErrors(true);           // Obtain squared errors, normalised
                                                   // w.r.t. the global squared error:
    mesher.markRef_into(elErrors,markRef);         // Mark elements for refinement
    mesher.markCrs_into(elErrors,markRef,markCrs); // Mark elements for coarsening,
                                                   // provided markRef
    mesher.refine(markRef);                        // Refine the geometry
    mesher.unrefine(markCrs);                      // Unrefine the geometry
}
\end{lstlisting}
This code, which is a simplified version of the actual code used for the wrinkling example (see \gs{benchmark_Wrinkling_DWR.cpp}), can be used for any goal functional chosen in \gs{gsThinShellAssemblerDWR} and for refinement using hierarchical and truncated hierarchical B-splines, which are refined using the suitable grading algorithms from \cite{Bracco2018} by default within the \gs{gsAdaptiveMeshing} and \gs{gsHBoxContainer} classes. For more information, the reader is referred to the documentation of these classes and the corresponding examples. In \cref{fig:adaptiveMeshing_wrinkling_results,fig:adaptiveMeshing_wrinkling_results_meshes}, the results of the goal-adaptive meshing simulation of the tension wrinkling benchmark problem are provided. For more details, the reader is referred to \cite{Verhelst2023Adaptive}.

\section{Summary}\label{sec:Conclusion}

This paper introduces open-source software for structural stability analysis using isogeometric membranes and Kirchhoff--Love shell elements. The software is developed as modules within the Geometry + Simulation modules library (\gismo); this C++ library contains various routines for geometric processing and the analysis of systems of equations using isogeometric analysis. The novel modules in this paper contain new methods from recent publications of the authors, including hyperelastic material models for shells and membranes \cite{Verhelst2021,Verhelst2025}, goal-oriented error estimators for Kirchhoff--Love shells \cite{Verhelst2023Adaptive}, an adaptive parallel arc-length method \cite{Verhelst2023APALM} and various unstructured spline constructions for isogeometric analysis \cite{Verhelst2023Coupling}. The modules provide off-the-shelf solvers for isogeometric analysis, aimed at future extension with novel material models, fast solvers for structural analysis, and new unstructured spline constructions. Through brief code snippets, examples, and brief class details, this paper aims to inspire the reader for the design of similar software packages and to stimulate the reader to use the presented modules for verification purposes or to use them in new developments. Overall, the experiments presented in this paper demonstrate that isogeometric analysis is able to provide efficient modeling tools and robust solutions to shell problems and beyond.\\

\backmatter

\bmhead{Acknowledgements}
The authors want to acknowledge all contributors to the Geometry + Simulation Modules. Specifically, the authors want to thank A. Farahat and P. Weinmüller for their contributions to the Unstructured Splines module and J. Li to their contributions to the Kirchhoff-Love shell module.

\section*{Declarations}

\paragraph{Data availability} No data are associated with this article.
\paragraph{Material availability} No physical materials are associated with this article.
\paragraph{Code availability} The code corresponding to this article is made publicly available through the Zenodo DOIs in \cref{tab:gismo_versions}, for completeness: \url{https://doi.org/10.5281/zenodo.15874172} (\code{gismo v25.07.0}), \url{https://doi.org/10.5281/zenodo.15875018} (\code{gsKLShell v25.07.0}), \url{https://doi.org/10.5281/zenodo.15875033} (\code{gsStructuralAnalysis v25.07.0}) and \url{https://doi.org/10.5281/zenodo.15875020} (\code{gsUnstructuredSplines v25.07.0}).
\paragraph{Competing interests} The authors declare that they have no conflicts of interest.
\paragraph{Funding} HMV is grateful for the financial support from the Italian Ministry of University and Research (MUR) through the PRIN projects COSMIC (no. 2022A79M75) and ASTICE (no. 202292JW3F), with the contribution of the European union -- Next Generation EU.






\begin{appendices}
\clearpage
\section{Installation}\label{sec:installation}
The software in this paper is part of the Geometry + Simulation Modules (\gismo). Therefore, the installation instructions for the software presented in this paper are derived from \gismo. The following instructions are for MacOS and Linux systems. For Windows, the reader is referred to the \gismo documentation.\\

\subsection{Downloading G+Smo}
The examples provided in this paper are reproduced using \gismo and its submodules with the versions as listed in \cref{tab:gismo_versions}. To obtain the \gismo repository from GitHub, one can use the following command:
\ResetLineNumber
\begin{lstlisting}[language=bash]
git clone --branch v25.07.0 https://github.com/gismo/gismo
\end{lstlisting}
The submodules \gs{gsKLShell}, \gs{gsStructuralAnalysis} and \gs{gsUnstructuredSplines} are optional submodules within \gismo. If enabled during build, they will be automatically cloned into their latest commit in the \code{optional/} directory of \gismo. However, for the reproduction of the work in this paper, we recommend manual cloning of the submodules with their correct version into the \gs{optional/} directory:
\ContinueLineNumber
\begin{lstlisting}[language=bash]
# Navigate to the optional/ folder in gismo/
cd optional

# Clone the gsKLShell, gsStructuralAnalysis and gsUnstructuredSplines modules
git clone --branch v25.07.0 https://github.com/gismo/gsKLShell
git clone --branch v25.07.0 https://github.com/gismo/gsStructuralAnalysis
git clone --branch v25.07.0 https://github.com/gismo/gsUnstructuredSplines
\end{lstlisting}

\begin{table}[h]
    \centering
    \caption{Versions for \gismo and the submodules \gs{gsKLShell}, \gs{gsStructuralAnalysis}, \gs{gsUnstructuredSplines} used in this paper.}
    \label{tab:gismo_versions}
    \footnotesize
    \begin{tabular}{p{0.2\linewidth}p{0.075\linewidth}p{0.425\linewidth}p{0.15\linewidth}}
        \toprule
        \textbf{Module} & \textbf{Version} & \textbf{GitHub Repository}\textsuperscript{\textdagger} & \textbf{Zenodo DOI}\textsuperscript{\textdaggerdbl}\\
        \midrule
        \gs{gismo} & \texttt{v25.07.0} & \href{https://github.com/gismo/gismo/releases/tag/v25.07.0}{\texttt{gismo/releases/tag/v25.07.0}} & \href{https://doi.org/10.5281/zenodo.15874172}{\texttt{zenodo.15874172}}\\
        \gs{gsKLShell} & \texttt{v25.07.0} & \href{https://github.com/gismo/gsKLShell/releases/tag/v25.07.0}{\texttt{gsKLShell/releases/tag/v25.07.0}} & \href{https://doi.org/10.5281/zenodo.15875018}{\texttt{zenodo.15875018}}\\
        \gs{gsStructuralAnalysis} & \texttt{v25.07.0} & \href{https://github.com/gismo/gsStructuralAnalysis/releases/tag/v25.07.0}{\texttt{gsStructuralAnalysis/releases/tag/v25.07.0}} & \href{https://doi.org/10.5281/zenodo.15875033}{\texttt{zenodo.15875033}}\\
        \gs{gsUnstructuredSplines} & \texttt{v25.07.0} & \href{https://github.com/gismo/gsUnstructuredSplines/releases/tag/v25.07.0}{\texttt{gsUnstructuredSplines/releases/tag/v25.07.0}} & \href{https://doi.org/10.5281/zenodo.15875020}{\texttt{zenodo.15875020}}\\
        \bottomrule
    \end{tabular}
{}\textsuperscript{\textdagger}: Preceded by \texttt{https://github.com/gismo/}\\
{}\textsuperscript{\textdaggerdbl}: Preceded by \texttt{https://doi.org/10.5281/}
\end{table}

\subsection{Installing G+Smo}
Installation of \gismo is done using CMake configuration. In the following\footnote{The procedure is tested on a notebook with Ubuntu 22.04 LTS, an Intel i7 13-700H CPU and 32GB RAM. For other installation instructions, the reader is referred to the \code{README.md} file in the \gismo repository.}, \gismo is built with \code{OpenMP} parallelization for assembly, MPI parallelization for the APALM method, using the Spectra \cite{Qiu2023} eigenvalue solver for sparse matrices and using the \code{optim} optimization library \cite{OHara2025} optionally used in the file to reproduce \cref{fig:cylinder}
\ResetLineNumber
\begin{lstlisting}[language=bash]
# Navigate to the directory of gismo
cd <path/to/gismo>

# Make a build folder and go inside
mkdir build
cd build

# Initialize the build with the correct CMake flags
cmake . -DGISMO_WITH_OPENMP=ON
        -DGISMO_WITH_MPI=ON
        -DGISMO_SUBMODULES="gsSpectra;gsOptim;gsKLShell;gsStructuralAnalysis;gsUnstructuredSplines"

# Build gismo
make gsKLShell-examples gsStructuralAnalysis-all gsUnstructuredSplines-examples
\end{lstlisting}

\subsection{Verifying the Installation}
The installation of \gismo can be verified using the unit-tests provided within the library. The unit-tests in \gismo can be compiled using the following commands:
\begin{lstlisting}[language=bash]
    # Set the path to the Pybind11 directory
    cmake . -DGISMO_BUILD_UNITTESTS=ON

    # Build pygismo
    make unittests
\end{lstlisting}
Consequently, the unit-tests related to the modules presented in this paper can be run as follows:
\begin{lstlisting}[language=bash]
# From the build folder
./bin/unittests
\end{lstlisting}

\clearpage
\section{Result Reproduction}
\label{sec:reproduction}
Most of the examples provided in this paper are based on results from earlier papers. For the reproducibility of the examples corresponding to earlier papers the reader is therefore directed to the respective paper for reproduction instructions. For the newly provided examples in this paper, \cref{tab:gismo_reproducibility} provides reproduction instructions.

\begin{table}[!h]
    \centering
    \caption{File name and run arguments required for the reproducibility of the figures in this paper. Arguments with a single dash (\code{-}) require an argument.}
    \footnotesize
    \begin{tabular}{p{0.10\textwidth}p{0.05\textwidth}p{0.30\textwidth}p{0.4\textwidth}}
        \toprule
        \textbf{Figure} & \multicolumn{3}{l}{\texttt{Run File}}\\
        & \textbf{Arg.} & \textbf{Description} & \textbf{Values}\\
        \midrule
        \Cref{fig:hyperelasticity_uniaxial} & \multicolumn{3}{l}{\code{benchmark_UniaxialTension}}\\
        & \code{-M} & Material model & \code{1}: NH, \code{2}: MR, \code{4}: OG (with \code{I 3})\\
        & \code{-I} & Implementation & \code{1}: Analytical, \code{2}: Generalised, \code{3}: Spectral\\
        & \code{-c} & Compressibility & \code{0}: Incompressible, \code{1}: Compressible\\
        \midrule
        \Cref{fig:cylinder} & \multicolumn{3}{l}{\code{benchmark_cylinder_DC}}\\
        & \code{-e} & Number of degree elevation steps & \code{1}: $p=2$,\newline \code{2}: $p=3$\\
        & \code{-r} & Number of uniform refinement steps & \code{4}: $256$ elements,\newline \code{5}: $1024$ elements\\
        & \code{-t} & Test case & \code{0}: Annulus,\newline \code{1}: Cylinder\\
        & \code{-L} & Displacement step length & \code{1e-1}\\
        & \code{-N} & Number of displacement step & \code{10}\\
        & \code{--NR} & Use Newton--Raphson solver & \\
        & \code{--DR} & Use Dynamic--Relaxation solver & \\
        & \code{-a} & Mass scaling parameter for DR & \code{1e14}\\
        & \code{-a} & Damping parameter for DR & \code{0} (use kinetic damping)\\
        & \code{--TFT} & Trigger TFT model & \\
        \multicolumn{4}{p{1.0\textwidth}}{\textbf{Note:} This benchmark can take several hours for both the shell model and the TFT-membrane model.}\\
        \midrule
        \Cref{fig:gismo_examples_coupling_weak_beam_plot} & \multicolumn{3}{l}{\code{static_shell_XML}}\\
        & \code{-i} & Input file & \code{gsStructuralAnalysis/filedata/pde/beam.xml}\\
        & \code{-e} & Number of degree elevation steps & \code{1}: $p=2$\\
        & \code{-r} & Number of uniform refinement steps & \code{2}: $4\times4$ elements per patch\\
        \multicolumn{4}{p{1.0\textwidth}}{\textbf{Post-processing:} The parameter $\alpha$ can be set via the \code{IfcPenalty} option in the XML file and the free end twist is obtained by $\theta=\tan^{-1}\qty(\frac{z_L-z_R}{\Delta x})$ with $z_L$ and $z_R$ the $z$-displacements on the left and right side of the flange, respectively, and $\Delta x$ the flange width, equal to $2$.}\\
        \bottomrule
    \end{tabular}
\end{table}

\begin{table}[!h]
    \centering
    \ContinuedFloat
    \caption{(Continued)}
    \label{tab:gismo_reproducibility}
    \footnotesize
    \begin{tabular}{p{0.10\textwidth}p{0.05\textwidth}p{0.30\textwidth}p{0.4\textwidth}}
        \toprule
        \textbf{Figure} & \multicolumn{3}{l}{\texttt{Run File}}\\
        & \textbf{Arg.} & \textbf{Description} & \textbf{Values}\\
        \midrule
        \Cref{fig:coupling_car_modes} & \multicolumn{3}{l}{\code{kirchhoff-Love_multipatch_vibration_XML_example}}\\
        \Cref{tab:ModalShell}
        & \code{-G} & Geometry file & \\
        & & \multicolumn{2}{l}{\code{surfaces/neon/dpatch_p2_s1_r1_geom.xml}}\\
        & & \multicolumn{2}{l}{\code{surfaces/neon/dpatch_p2_s1_r2_geom.xml} \code{...}}
        \\
        & & \multicolumn{2}{l}{\code{surfaces/neon/almostC1_p2_s1_r1_geom.xml}}\\
        & & \multicolumn{2}{l}{\code{surfaces/neon/almostC1_p2_s1_r2_geom.xml} \code{...}}
        \\
        & \code{-b} & Basis file & \\
        & & \multicolumn{2}{l}{\code{surfaces/neon/dpatch_p2_s1_r1_basis.xml}}\\
        & & \multicolumn{2}{l}{\code{surfaces/neon/dpatch_p2_s1_r2_basis.xml} \code{...}}\\
        & & \multicolumn{2}{l}{\code{surfaces/neon/almostC1_p2_s1_r1_basis.xml}}\\
        & & \multicolumn{2}{l}{\code{surfaces/neon/almostC1_p2_s1_r2_basis.xml} \code{...}}\\
        & \code{-B} & Problem definition file & \code{pde/shell/car_bvp.xml}\\
        \midrule
        \Cref{fig:hyperelasticity_wrinkling_results} & \multicolumn{3}{l}{\code{benchmark_TensionWrinkling}}\\
        & \code{-M} & Material model & \code{1}: NH, \code{2}: MR, \code{4}: OG (with \code{I 3})\\
        & \code{-I} & Implementation & \code{1}: Analytical, \code{2}: Generalised, \code{3}: Spectral\\
        \multicolumn{4}{p{1.0\textwidth}}{\textbf{Post-processing:} The displacements in the end-point is provided in $[\text{m}]$, hence $x$-displacement has to be divided by the length of the sheet to obtain the strain and the $z$-displacement has to be normalized with the thickness.}\\
        \midrule
        \Cref{fig:APALM_snapping} & \multicolumn{3}{l}{\code{snapping_example_shell_APALM}}\\
        &
        \multicolumn{3}{l}{\code{snapping_example_shell_DC}}\\
        & \code{-e} & Number of degree elevations & \code{1}\\
        & \code{-L} & Load-step size & \code{2.5e-4}: DC,\newline  \code{5e-2}: APALM\\
        & \code{-X},\code{-Y} & Number of unit-cells in $x$- and $y$-direction & \code{3,3}\\
        \multicolumn{4}{p{1.0\textwidth}}{\textbf{Note:} This benchmark can be ran with a combination of OpenMP and MPI commands, distribution the number of ranks and the number of threads per rank over the available processor. On SLURM-compatible system, the reader is referred to the SLURM documentation. Due to the large number of load steps, this example might take hours depending on the available resources.}\\
        \midrule
        \Cref{fig:adaptiveMeshing_wrinkling_results} & \multicolumn{3}{l}{\code{benchmark_Wrinkling_DWR}}\\
        \Cref{fig:adaptiveMeshing_wrinkling_results_meshes} & \code{-r} & Number of initial uniform refinements & \code{5}: $32\times32$, \code{6}: $64\times64$,\\
        & \code{-T} & Target error & \code{1e-12}\\
        & \code{-B} & Error lower-bound multiplier & \code{0.1}\\
        & \code{-D} & Coasening error limit & \code{1e-19}\\
        & \code{-g} & Goal functional & \code{1}: Displacement\\
        & \code{-C} & ...and its component & \code{2}: $x$-component\\
        & \code{-O} & Adaptive meshing options & \code{shell_mesher_options_Wrinkling.xml}\\
        & \code{--adaptMesh} & Perform refinement loops & \\
        \multicolumn{4}{p{1.0\textwidth}}{\textbf{Post-processing:} The same as \code{benchmark_TensionWrinkling}}\\
        \multicolumn{4}{p{1.0\textwidth}}{\textbf{Note:} This benchmark can take several hours.}\\
        \bottomrule
    \end{tabular}
\end{table}

\clearpage
\section{Class overview per module}
\label{sec:overview}

\begin{table}[h]
          \centering
          \caption{Overview of the main classes in the \gs{gsKLShell} module. The template argument \gs{d} denotes the parametric dimension, and the argument \gs{T} denotes the type for real numbers, e.g., \gs{double} or \gs{long double}. The symbol \bottombar denotes inheritance.}
          \label{tab:gsKLShell}
          \footnotesize
          \begin{tabular}{p{0.35\linewidth}p{0.55\linewidth}}
                    \toprule
                    \textbf{Class} & \textbf{Description} \\
                    \midrule
                    \gs{gsMaterialMatrixBase} & Base class providing the constitutive implementation for the \gs{gsThinShellAssembler}.\\
                    \bottombar \gs{gsMaterialMatrixBaseDim} & Dimension-dependent base class inheriting from \gs{gsMaterialMatrixBase} and providing all necessary basis computations for the constitutive relations.\\
                    \quad\bottombar \gs{gsMaterialMatrixLinear} & Implementation of the Saint-Venant Kirchoff constitutive relation.\\
                    \quad\bottombar \gs{gsMaterialMatrixComposite} & Uses the Saint-Venant Kirchoff constitutive relation for laminates with different orientations, described in \cite{Herrema2019}\\
                    \quad\bottombar \gs{gsMaterialMatrixNonlinear} & Implementation of the Neo-Hookean, Mooney-Rivlin, and Ogden material models, following the works of \cite{Kiendl2015} and \cite{Verhelst2021}. The class is templated over the material model \gs{mat} and over the implementation: (i) $\Sten$ and $\Ccten$ are given analytically; (ii) $\Psi$ and its derivatives to $\Cten$ are given analytically; or (iii) $\Psi$ and its derivatives to $\lambda$ are given analytically.\\
                    \quad\bottombar \gs{gsMaterialMatrixTFT} & A tension field theory-based material matrix for wrinkling modeling. This material matrix uses another material matrix as an input and applies the modification scheme from \cite{Nakashino2005} for linear elastic materials or the scheme from \cite{Verhelst2025} for hyperelastic materials to implicitly model wrinkling.\\
                    \midrule
                    \gs{gsThinShellAssembler} & The actual assembler for the Kirchhoff--Love shell model. The assembler takes a geometry (\gs{gsFunctionSet}), a basis (\gs{gsMultiBasis} or \gs{gsMappedBasis}), a set of boundary conditions (\gs{gsBoundaryConditions}), a body force (\gs{gsFunctionSet}), and a material matrix (\gs{gsMaterialMatrixBase}) as input. Optionally, point loads, follower pressures, or an elastic foundation can be provided. The class provides the linear stiffness matrix $\MAT{K}$, the external force vector $\VEC{P}$, the residual vector $\VEC{R}(\uvec)$, the Jacobian $\MAT{K}(\uvec)$, and the mass matrix $\MAT{M}$.\\
                    \gs{gsThinShellAssemblerDWR} & The same as \gs{gsThinShellAssembler}, but takes a primal and a dual space (both \gs{gsMultiBasis}) and a goal functional. The class outputs operators for the Dual-Weighted Residual (DWR) method.\\
                    \bottomrule
          \end{tabular}
\end{table}

\begin{table}[t]
\centering
\caption{Most important solvers in the \gs{gsStructuralAnalysis} module. The symbol \bottombar denotes inheritance.}
\label{tab:gsStructuralAnalysis}
\footnotesize
\begin{tabular}{p{0.25\linewidth}p{0.65\linewidth}}
\toprule
\textbf{Class} & \textbf{Description} \\
\midrule
\gs{gsStaticBase} & Provides a base class for static analysis. The class provides common functions for the derived classes, such as functions accessing solutions or stability computations.\\
\bottombar \gs{gsStaticDR} & Class that implements the Dynamic Relaxation (DR) method. The class takes the residual vector $\VEC{R}(\uvec)$ and the (lumped) mass matrix $\MAT{M}$.\\
\bottombar \gs{gsStaticNewton} & Class that implements Newton-Raphson (NR) iterations for solving non-linear static problems. The class takes the force vector $\VEC{P}$, the residual vector $\VEC{R}(\uvec)$, the stiffness matrix $\MAT{K}$, and the Jacobian operator $\MAT{K}(\uvec)$ or the MIP Jacobian $\MAT{K}(\uvec,\Delta\uvec)$. \\
\bottombar \gs{gsStaticComposite} & Class that implements a sequential static solver. Multiple static solvers with different tolerances are provided, and this class sequentially calls the static solvers and provides the updates.\\
\midrule
\gs{gsEigenProblemBase} & Base class for eigenvalue problems for structural analysis. The class provides the solution routines, whereas the derived classes provide the set-up of the eigenvalue problem. The eigenproblems are solved using the generalised eigensolver of \gs{Eigen} \cite{eigenweb} or using \gs{Spectra} \cite{Qiu2023} for large problems, optionally with shifts.\\
\bottombar \gs{gsModalSolver} & A solver that performs modal analysis, i.e., that solves the eigenvalue problem \cref{eq:gismo_ModalAnalysisProblem}, hence requires the mass and linear stiffness matrices $\MAT{M}$ and $\MAT{K}$, respectively.\\
\bottombar \gs{gsBucklingSolver} & A solver that performs linear buckling analysis, i.e., that solves the eigenvalue problem \cref{eq:gismo_gsStructuralAnalysis_BucklingAnalysisProblem}. The solver can be enabled by providing the matrices $\MAT{K}$ and $\MAT{K}(\uvec_0)$ or by providing the matrix $\MAT{K}$, the load vector $\VEC{P}$, and the operator for the Jacobian $\MAT{K}(\uvec)$.\\
\midrule
\gs{gsALMBase} & Base class for Arc-Length Methods (ALMs). This class also includes the extended ALM \cite{Wriggers1988}, used for computing singular points. All solvers require the external load vector $\VEC{P}$, the arc-length residual $\VEC{R}(\uvec,\lambda)$, and the Jacobian matrix $\MAT{K}(\uvec)$.\\
\bottombar \gs{gsALMLoadControl} & A class performing load-controlled simulations.\\
\bottombar \gs{gsALMRiks} & An implementation of Riks' arc-length method, see \cite{Riks1979}.\\
\bottombar \gs{gsALMCrisfield} & A class using Crisfield's method, see \cite{Crisfield1981} and \cref{eq:gismo_gsStructuralAnalysis_Crisfield}.\\
\bottomrule
\end{tabular}
\end{table}

\begin{table}[t]
\centering
\caption{Most important classes in the \gs{gsUnstructuredSplines} module. The symbol \bottombar denotes inheritance.}
\label{tab:gsUnstructuredSplines}
\footnotesize
\begin{tabular}{p{0.25\linewidth}p{0.65\linewidth}}
\toprule
\textbf{Class} & \textbf{Description} \\
\midrule
\gs{gsDPatchBase} & The base class for smoothing methods that rely on local refinements. The class contains basic routines for interface and vertex smoothing and requires specifications to handle extraordinary vertices for derived classes. \\
\bottombar \gs{gsSmoothInterfaces} & The most basic class for smoothing interfaces. It provides $C^1$ smoothing on interfaces and ordinary vertices but does not enforce exact or almost $C^1$ smoothing at extraordinary vertices. The class can be initialised for tensor B-spline bases and geometries or THB-spline geometries with degree $p\geq2$ and regularity $0<r\leq p-1$\\
\bottombar \gs{gsDPatch} & Class implementing the degenerate patch (D-Patch) approach from \cite{Toshniwal2017}. It relies on local THB refinement around the extraordinary vertices, building a space with reduced local regularity and refinement. The class can be initialised for tensor B-spline bases and geometries with degree $p\geq2$ and regularity $0<r\leq p-1$\\
\bottombar \gs{gsAlmostC1} & Class implementing the almost-$C^1$ method from \cite{Takacs2023}. As the D-Patch, it relies on local THB refinement. The class can be initialised for tensor B-spline bases and geometries or THB-spline geometries with degree $p=2$ and regularity $r=1$\\
\gs{gsC1SurfSpline} & Class implementing the analysis-suitable $G^1$ (AS-$G^1$) spline construction from \cite{Farahat2023,Farahat2023a}. This class relies on different local bases for the vertex, interface, and interior spaces and requires a degree of $p\geq3$ and a regularity of $r\leq p-2$\\
\gs{gsApproxC1Spline} & Class implementing the approximate $C^1$ spline construction from \cite{Weinmuller2021,Weinmuller2022}. Its implementation is similar to the AS-$G^1$ construction, but it can be constructed for large regularity, namely, it requires a degree of $p\geq3$ and a regularity of $r\leq p-1$\\
\gs{gsMPBESSpline} & Class implementing the multi-patch B-splines with enhanced smoothness (MPBES) construction from \cite{Buchegger2016}. The construction requires $p\geq 2$ and $r \leq p-1$. Note that this construction is $C^{p-1}$ everywhere, but only $C^0$ in the vicinity of the EV.\\
\bottomrule
\end{tabular}
\end{table}






\end{appendices}


\clearpage
\bibliography{sn-bibliography}

\end{document}

%% file: Figures/Illustration_patches_from_mesh1.tikz
\begin{tikzpicture}
[
        ,boundary/.style={thick, draw=black},
        ,interior/.style={thin, draw=gray},
		,green_thin/.style={thin, draw=gray},
		,red_thin/.style={thin, draw=gray},
		,blue_thin/.style={thin, draw=gray},
		,orange_thin/.style={thin, draw=gray},
		,yellow_thin/.style={thin, draw=gray},
		,pink_thin/.style={thin, draw=gray},
		,purple_thin/.style={thin, draw=gray},
		,brown_thin/.style={thin, draw=gray},
		,lime_thin/.style={thin, draw=gray},
		,cyan_thin/.style={thin, draw=gray},
		,magenta_thin/.style={thin, draw=gray},
		,EV/.style={fill=gray},
		,bEV/.style={fill=black},
]
		\node [] (0) at (0, 3) {};
		\node [] (1) at (0, 3.5) {};
		\node [] (2) at (0, 3.95) {};
		\node [] (3) at (0, 2) {};
		\node [] (4) at (0, 1.5) {};
		\node [] (5) at (0, 1) {};
		\node [] (6) at (0, 0.5) {};
		\node [] (7) at (0, 4.5) {};
		\node [] (8) at (0.5, 2.5) {};
		\node [] (9) at (1, 2.5) {};
		\node [] (10) at (1.5, 2.5) {};
		\node [] (11) at (2, 2.5) {};
		\node [] (12) at (-0.5, 2.5) {};
		\node [] (13) at (-1, 2.5) {};
		\node [] (14) at (-1.45, 2.5) {};
		\node [] (15) at (-2, 2.5) {};
		\node [] (16) at (-0.675, 3.175) {};
		\node [] (17) at (-1, 3.5) {};
		\node [] (18) at (-1.425, 3.925) {};
		\node [] (19) at (0.675, 3.175) {};
		\node [] (20) at (1, 3.5) {};
		\node [] (21) at (1.425, 3.925) {};
		\node [] (22) at (-1.425, 1.075) {};
		\node [] (23) at (-1, 1.5) {};
		\node [] (24) at (-0.675, 1.8) {};
		\node [] (25) at (0.35, 2.15) {};
		\node [] (26) at (1, 1.5) {};
		\node [] (27) at (0.5, 1.5) {};
		\node [] (28) at (1, 2) {};
		\node [] (29) at (1.5, 2) {};
		\node [] (30) at (1.5, 1.15) {};
		\node [] (31) at (0.5, 1) {};
		\node [] (32) at (1, 1) {};
		\node [] (36) at (2, 2) {};
		\node [] (37) at (-0.325, 2.825) {};
		\node [] (38) at (-0.325, 2.175) {};
		\node [] (39) at (0.325, 2.825) {};
		\node [] (40) at (0, 0) {};
		\node [] (41) at (0.5, 0) {};
		\node [] (42) at (1, 0) {};
		\node [] (47) at (2, 1) {};
		\node [] (54) at (1.5, 0.5) {};
		\node [] (55) at (0.825, 0.5) {};
		\node [] (56) at (0.5, 0.5) {};
		\node [] (57) at (2, 1.5) {};
		\node [] (58) at (1.5, 1.5) {};

		\draw [style=boundary] (11.center) to (21.center);
		\draw [style=boundary] (21.center) to (7.center);
		\draw [style=boundary] (7.center) to (18.center);
		\draw [style=boundary] (18.center) to (15.center);
		\draw [style=boundary] (15.center) to (22.center);
		\draw [style=boundary] (22.center) to (6.center);
		\draw [style=boundary] (8.center) to (25.center);
		\draw [style=boundary] (25.center) to (3.center);
		\draw [style=interior] (19.center) to (20.center);
		\draw [style=boundary] (3.center) to (38.center);
		\draw [style=boundary] (38.center) to (12.center);
		\draw [style=boundary] (12.center) to (37.center);
		\draw [style=boundary] (37.center) to (0.center);
		\draw [style=boundary] (0.center) to (39.center);
		\draw [style=boundary] (39.center) to (8.center);
		\draw [style=interior] (24.center) to (38.center);
		\draw [style=interior] (12.center) to (13.center);
		\draw [style=interior] (37.center) to (16.center);
		\draw [style=interior] (0.center) to (1.center);
		\draw [style=interior] (39.center) to (19.center);
		\draw [style=interior] (16.center) to (17.center);
		\draw [style=interior] (20.center) to (21.center);
		\draw [style=interior] (1.center) to (2.center);
		\draw [style=interior] (2.center) to (7.center);
		\draw [style=interior] (17.center) to (18.center);
		\draw [style=interior] (14.center) to (15.center);
		\draw [style=interior] (23.center) to (22.center);
		\draw [style=interior] (23.center) to (24.center);
		\draw [style=interior] (13.center) to (14.center);
		\draw [style=boundary] (40.center) to (41.center);
		\draw [style=boundary] (41.center) to (42.center);
		\draw [style=boundary] (6.center) to (40.center);
		\draw [style={green_thin}] (32.center) to (26.center);
		\draw [style={green_thin}] (26.center) to (28.center);
		\draw [style={green_thin}] (28.center) to (9.center);
		\draw [style={green_thin}] (9.center) to (19.center);
		\draw [style={green_thin}] (19.center) to (1.center);
		\draw [style={green_thin}] (1.center) to (16.center);
		\draw [style={green_thin}] (16.center) to (13.center);
		\draw [style={green_thin}] (13.center) to (24.center);
		\draw [style={green_thin}] (24.center) to (4.center);
		\draw [style={green_thin}] (4.center) to (27.center);
		\draw [style={green_thin}] (27.center) to (26.center);
		\draw [style={red_thin}] (32.center) to (31.center);
		\draw [style={red_thin}] (31.center) to (5.center);
		\draw [style={red_thin}] (5.center) to (23.center);
		\draw [style={red_thin}] (23.center) to (14.center);
		\draw [style={red_thin}] (14.center) to (17.center);
		\draw [style={red_thin}] (17.center) to (2.center);
		\draw [style={red_thin}] (2.center) to (20.center);
		\draw [style={red_thin}] (20.center) to (10.center);
		\draw [style={red_thin}] (10.center) to (29.center);
		\draw [style={blue_thin}] (6.center) to (5.center);
		\draw [style={blue_thin}] (5.center) to (4.center);
		\draw [style={blue_thin}] (4.center) to (3.center);
		\draw [style={pink_thin}] (25.center) to (28.center);
		\draw [style={pink_thin}] (28.center) to (29.center);
		\draw [style={pink_thin}] (29.center) to (36.center);
		\draw [style={purple_thin}] (25.center) to (27.center);
		\draw [style={purple_thin}] (27.center) to (31.center);
		\draw [style={purple_thin}] (56.center) to (41.center);
		\draw [style={purple_thin}] (31.center) to (56.center);
		\draw [style={interior}] (11.center) to (10.center);
		\draw [style={interior}] (10.center) to (9.center);
		\draw [style={interior}] (9.center) to (8.center);
		\draw [style={boundary}](42.center) to (54.center);
		\draw [style={boundary}](54.center) to (47.center);
		\draw [style={orange_thin}] (32.center) to (54.center);
		\draw [style={lime_thin}] (30.center) to (47.center);
		\draw [style={brown_thin}] (55.center) to (42.center);
		\draw [style={cyan_thin}] (32.center) to (55.center);
		\draw [style={magenta_thin}] (32.center) to (30.center);
		\draw [style={yellow_thin}] (56.center) to (6.center);
		\draw [style={yellow_thin}] (55.center) to (56.center);
		\draw [style={green_thin}] (26.center) to (58.center);
		\draw [style={red_thin}] (58.center) to (30.center);
		\draw [style={red_thin}] (29.center) to (58.center);
		\draw [style={green_thin}] (58.center) to (57.center);
		\draw [style={boundary}](11.center) to (36.center);
		\draw [style={boundary}](36.center) to (57.center);
		\draw [style={boundary}](57.center) to (47.center);

		\fill[bEV] (6) circle [radius=0.05];
		\fill[EV] (55) circle [radius=0.05];
		\fill[EV] (30) circle [radius=0.05];
		\fill[EV] (32) circle [radius=0.05];
		\fill[bEV] (25) circle [radius=0.05];
\end{tikzpicture}

%% file: Figures/Illustration_patches_from_mesh2.tikz
\begin{tikzpicture}
[
        ,boundary/.style={thick, draw=black},
        ,interior/.style={thin, draw=gray},
		,green_thin/.style={thin, draw=TUcol4,-latex},
		,red_thin/.style={thin, draw=TUcol5,-latex},
		,blue_thin/.style={thin, draw=gray},
		,orange_thin/.style={thin, draw=TUcol7,-latex},
		,yellow_thin/.style={thin, draw=gray},
		,pink_thin/.style={thin, draw=gray},
		,purple_thin/.style={thin, draw=gray},
		,brown_thin/.style={thin, draw=gray},
		,lime_thin/.style={thin, draw=gray},
		,cyan_thin/.style={thin, draw=TUcol13,-latex},
		,magenta_thin/.style={thin, draw=TUcol14,-latex},
		,EV/.style={fill=gray},
		,bEV/.style={fill=black},
]
		\node [] (0) at (0, 3) {};
		\node [] (1) at (0, 3.5) {};
		\node [] (2) at (0, 3.95) {};
		\node [] (3) at (0, 2) {};
		\node [] (4) at (0, 1.5) {};
		\node [] (5) at (0, 1) {};
		\node [] (6) at (0, 0.5) {};
		\node [] (7) at (0, 4.5) {};
		\node [] (8) at (0.5, 2.5) {};
		\node [] (9) at (1, 2.5) {};
		\node [] (10) at (1.5, 2.5) {};
		\node [] (11) at (2, 2.5) {};
		\node [] (12) at (-0.5, 2.5) {};
		\node [] (13) at (-1, 2.5) {};
		\node [] (14) at (-1.45, 2.5) {};
		\node [] (15) at (-2, 2.5) {};
		\node [] (16) at (-0.675, 3.175) {};
		\node [] (17) at (-1, 3.5) {};
		\node [] (18) at (-1.425, 3.925) {};
		\node [] (19) at (0.675, 3.175) {};
		\node [] (20) at (1, 3.5) {};
		\node [] (21) at (1.425, 3.925) {};
		\node [] (22) at (-1.425, 1.075) {};
		\node [] (23) at (-1, 1.5) {};
		\node [] (24) at (-0.675, 1.8) {};
		\node [] (25) at (0.35, 2.15) {};
		\node [] (26) at (1, 1.5) {};
		\node [] (27) at (0.5, 1.5) {};
		\node [] (28) at (1, 2) {};
		\node [] (29) at (1.5, 2) {};
		\node [] (30) at (1.5, 1.15) {};
		\node [] (31) at (0.5, 1) {};
		\node [] (32) at (1, 1) {};
		\node [] (36) at (2, 2) {};
		\node [] (37) at (-0.325, 2.825) {};
		\node [] (38) at (-0.325, 2.175) {};
		\node [] (39) at (0.325, 2.825) {};
		\node [] (40) at (0, 0) {};
		\node [] (41) at (0.5, 0) {};
		\node [] (42) at (1, 0) {};
		\node [] (47) at (2, 1) {};
		\node [] (54) at (1.5, 0.5) {};
		\node [] (55) at (0.825, 0.5) {};
		\node [] (56) at (0.5, 0.5) {};
		\node [] (57) at (2, 1.5) {};
		\node [] (58) at (1.5, 1.5) {};

		\draw [style=boundary] (11.center) to (21.center);
		\draw [style=boundary] (21.center) to (7.center);
		\draw [style=boundary] (7.center) to (18.center);
		\draw [style=boundary] (18.center) to (15.center);
		\draw [style=boundary] (15.center) to (22.center);
		\draw [style=boundary] (22.center) to (6.center);
		\draw [style=boundary] (8.center) to (25.center);
		\draw [style=boundary] (25.center) to (3.center);
		\draw [style=interior] (19.center) to (20.center);
		\draw [style=boundary] (3.center) to (38.center);
		\draw [style=boundary] (38.center) to (12.center);
		\draw [style=boundary] (12.center) to (37.center);
		\draw [style=boundary] (37.center) to (0.center);
		\draw [style=boundary] (0.center) to (39.center);
		\draw [style=boundary] (39.center) to (8.center);
		\draw [style=interior] (24.center) to (38.center);
		\draw [style=interior] (12.center) to (13.center);
		\draw [style=interior] (37.center) to (16.center);
		\draw [style=interior] (0.center) to (1.center);
		\draw [style=interior] (39.center) to (19.center);
		\draw [style=interior] (16.center) to (17.center);
		\draw [style=interior] (20.center) to (21.center);
		\draw [style=interior] (1.center) to (2.center);
		\draw [style=interior] (2.center) to (7.center);
		\draw [style=interior] (17.center) to (18.center);
		\draw [style=interior] (14.center) to (15.center);
		\draw [style=interior] (23.center) to (22.center);
		\draw [style=interior] (23.center) to (24.center);
		\draw [style=interior] (13.center) to (14.center);
		\draw [style=boundary] (40.center) to (41.center);
		\draw [style=boundary] (41.center) to (42.center);
		\draw [style=boundary] (6.center) to (40.center);
		\draw [style={green_thin}] (32.center) to (26.center);
		\draw [style={green_thin}] (26.center) to (28.center);
		\draw [style={green_thin}] (28.center) to (9.center);
		\draw [style={green_thin}] (9.center) to (19.center);
		\draw [style={green_thin}] (19.center) to (1.center);
		\draw [style={green_thin}] (1.center) to (16.center);
		\draw [style={green_thin}] (16.center) to (13.center);
		\draw [style={green_thin}] (13.center) to (24.center);
		\draw [style={green_thin}] (24.center) to (4.center);
		\draw [style={green_thin}] (4.center) to (27.center);
		\draw [style={green_thin}] (27.center) to (26.center);
		\draw [style={red_thin}] (32.center) to (31.center);
		\draw [style={red_thin}] (31.center) to (5.center);
		\draw [style={red_thin}] (5.center) to (23.center);
		\draw [style={red_thin}] (23.center) to (14.center);
		\draw [style={red_thin}] (14.center) to (17.center);
		\draw [style={red_thin}] (17.center) to (2.center);
		\draw [style={red_thin}] (2.center) to (20.center);
		\draw [style={red_thin}] (20.center) to (10.center);
		\draw [style={red_thin}] (10.center) to (29.center);
		\draw [style={blue_thin}] (6.center) to (5.center);
		\draw [style={blue_thin}] (5.center) to (4.center);
		\draw [style={blue_thin}] (4.center) to (3.center);
		\draw [style={pink_thin}] (25.center) to (28.center);
		\draw [style={pink_thin}] (28.center) to (29.center);
		\draw [style={pink_thin}] (29.center) to (36.center);
		\draw [style={purple_thin}] (25.center) to (27.center);
		\draw [style={purple_thin}] (27.center) to (31.center);
		\draw [style={purple_thin}] (56.center) to (41.center);
		\draw [style={purple_thin}] (31.center) to (56.center);
		\draw [style={interior}] (11.center) to (10.center);
		\draw [style={interior}] (10.center) to (9.center);
		\draw [style={interior}] (9.center) to (8.center);
		\draw [style={boundary}](42.center) to (54.center);
		\draw [style={boundary}](54.center) to (47.center);
		\draw [style={orange_thin}] (32.center) to (54.center);
		\draw [style={lime_thin}] (30.center) to (47.center);
		\draw [style={brown_thin}] (55.center) to (42.center);
		\draw [style={cyan_thin}] (32.center) to (55.center);
		\draw [style={magenta_thin}] (32.center) to (30.center);
		\draw [style={yellow_thin}] (56.center) to (6.center);
		\draw [style={yellow_thin}] (55.center) to (56.center);
		\draw [style={green_thin}] (26.center) to (58.center);
		\draw [style={red_thin}] (58.center) to (30.center);
		\draw [style={red_thin}] (29.center) to (58.center);
		\draw [style={green_thin}] (58.center) to (57.center);
		\draw [style={boundary}](11.center) to (36.center);
		\draw [style={boundary}](36.center) to (57.center);
		\draw [style={boundary}](57.center) to (47.center);

		\fill[bEV] (6) circle [radius=0.05];
		\fill[EV] (55) circle [radius=0.05];
		\fill[EV] (30) circle [radius=0.05];
		\fill[EV] (32) circle [radius=0.05];
		\fill[bEV] (25) circle [radius=0.05];
\end{tikzpicture}

%% file: Figures/Illustration_patches_from_mesh3.tikz
\begin{tikzpicture}
[
        ,boundary/.style={thick, draw=black},
        ,interior/.style={thin, draw=gray},
		,green_thin/.style={thin, draw=TUcol4},
		,red_thin/.style={thin, draw=TUcol5},
		,blue_thin/.style={thin, draw=TUcol6},
		,orange_thin/.style={thin, draw=TUcol7},
		,yellow_thin/.style={thin, draw=TUcol8},
		,pink_thin/.style={thin, draw=TUcol9},
		,purple_thin/.style={thin, draw=TUcol10},
		,brown_thin/.style={thin, draw=TUcol11},
		,lime_thin/.style={thin, draw=TUcol12},
		,cyan_thin/.style={thin, draw=TUcol13},
		,magenta_thin/.style={thin, draw=TUcol14},
		,EV/.style={fill=gray},
		,bEV/.style={fill=black},
]
		\node [] (0) at (0, 3) {};
		\node [] (1) at (0, 3.5) {};
		\node [] (2) at (0, 3.95) {};
		\node [] (3) at (0, 2) {};
		\node [] (4) at (0, 1.5) {};
		\node [] (5) at (0, 1) {};
		\node [] (6) at (0, 0.5) {};
		\node [] (7) at (0, 4.5) {};
		\node [] (8) at (0.5, 2.5) {};
		\node [] (9) at (1, 2.5) {};
		\node [] (10) at (1.5, 2.5) {};
		\node [] (11) at (2, 2.5) {};
		\node [] (12) at (-0.5, 2.5) {};
		\node [] (13) at (-1, 2.5) {};
		\node [] (14) at (-1.45, 2.5) {};
		\node [] (15) at (-2, 2.5) {};
		\node [] (16) at (-0.675, 3.175) {};
		\node [] (17) at (-1, 3.5) {};
		\node [] (18) at (-1.425, 3.925) {};
		\node [] (19) at (0.675, 3.175) {};
		\node [] (20) at (1, 3.5) {};
		\node [] (21) at (1.425, 3.925) {};
		\node [] (22) at (-1.425, 1.075) {};
		\node [] (23) at (-1, 1.5) {};
		\node [] (24) at (-0.675, 1.8) {};
		\node [] (25) at (0.35, 2.15) {};
		\node [] (26) at (1, 1.5) {};
		\node [] (27) at (0.5, 1.5) {};
		\node [] (28) at (1, 2) {};
		\node [] (29) at (1.5, 2) {};
		\node [] (30) at (1.5, 1.15) {};
		\node [] (31) at (0.5, 1) {};
		\node [] (32) at (1, 1) {};
		\node [] (36) at (2, 2) {};
		\node [] (37) at (-0.325, 2.825) {};
		\node [] (38) at (-0.325, 2.175) {};
		\node [] (39) at (0.325, 2.825) {};
		\node [] (40) at (0, 0) {};
		\node [] (41) at (0.5, 0) {};
		\node [] (42) at (1, 0) {};
		\node [] (47) at (2, 1) {};
		\node [] (54) at (1.5, 0.5) {};
		\node [] (55) at (0.825, 0.5) {};
		\node [] (56) at (0.5, 0.5) {};
		\node [] (57) at (2, 1.5) {};
		\node [] (58) at (1.5, 1.5) {};

		\draw [style=boundary] (11.center) to (21.center);
		\draw [style=boundary] (21.center) to (7.center);
		\draw [style=boundary] (7.center) to (18.center);
		\draw [style=boundary] (18.center) to (15.center);
		\draw [style=boundary] (15.center) to (22.center);
		\draw [style=boundary] (22.center) to (6.center);
		\draw [style=boundary] (8.center) to (25.center);
		\draw [style=boundary] (25.center) to (3.center);
		\draw [style=interior] (19.center) to (20.center);
		\draw [style=boundary] (3.center) to (38.center);
		\draw [style=boundary] (38.center) to (12.center);
		\draw [style=boundary] (12.center) to (37.center);
		\draw [style=boundary] (37.center) to (0.center);
		\draw [style=boundary] (0.center) to (39.center);
		\draw [style=boundary] (39.center) to (8.center);
		\draw [style=interior] (24.center) to (38.center);
		\draw [style=interior] (12.center) to (13.center);
		\draw [style=interior] (37.center) to (16.center);
		\draw [style=interior] (0.center) to (1.center);
		\draw [style=interior] (39.center) to (19.center);
		\draw [style=interior] (16.center) to (17.center);
		\draw [style=interior] (20.center) to (21.center);
		\draw [style=interior] (1.center) to (2.center);
		\draw [style=interior] (2.center) to (7.center);
		\draw [style=interior] (17.center) to (18.center);
		\draw [style=interior] (14.center) to (15.center);
		\draw [style=interior] (23.center) to (22.center);
		\draw [style=interior] (23.center) to (24.center);
		\draw [style=interior] (13.center) to (14.center);
		\draw [style=boundary] (40.center) to (41.center);
		\draw [style=boundary] (41.center) to (42.center);
		\draw [style=boundary] (6.center) to (40.center);
		\draw [style={green_thin}] (32.center) to (26.center);
		\draw [style={green_thin}] (26.center) to (28.center);
		\draw [style={green_thin}] (28.center) to (9.center);
		\draw [style={green_thin}] (9.center) to (19.center);
		\draw [style={green_thin}] (19.center) to (1.center);
		\draw [style={green_thin}] (1.center) to (16.center);
		\draw [style={green_thin}] (16.center) to (13.center);
		\draw [style={green_thin}] (13.center) to (24.center);
		\draw [style={green_thin}] (24.center) to (4.center);
		\draw [style={green_thin}] (4.center) to (27.center);
		\draw [style={green_thin}] (27.center) to (26.center);
		\draw [style={red_thin}] (32.center) to (31.center);
		\draw [style={red_thin}] (31.center) to (5.center);
		\draw [style={red_thin}] (5.center) to (23.center);
		\draw [style={red_thin}] (23.center) to (14.center);
		\draw [style={red_thin}] (14.center) to (17.center);
		\draw [style={red_thin}] (17.center) to (2.center);
		\draw [style={red_thin}] (2.center) to (20.center);
		\draw [style={red_thin}] (20.center) to (10.center);
		\draw [style={red_thin}] (10.center) to (29.center);
		\draw [style={blue_thin}] (6.center) to (5.center);
		\draw [style={blue_thin}] (5.center) to (4.center);
		\draw [style={blue_thin}] (4.center) to (3.center);
		\draw [style={pink_thin}] (25.center) to (28.center);
		\draw [style={pink_thin}] (28.center) to (29.center);
		\draw [style={pink_thin}] (29.center) to (36.center);
		\draw [style={purple_thin}] (25.center) to (27.center);
		\draw [style={purple_thin}] (27.center) to (31.center);
		\draw [style={purple_thin}] (56.center) to (41.center);
		\draw [style={purple_thin}] (31.center) to (56.center);
		\draw [style={interior}] (11.center) to (10.center);
		\draw [style={interior}] (10.center) to (9.center);
		\draw [style={interior}] (9.center) to (8.center);
		\draw [style={boundary}](42.center) to (54.center);
		\draw [style={boundary}](54.center) to (47.center);
		\draw [style={orange_thin}] (32.center) to (54.center);
		\draw [style={lime_thin}] (30.center) to (47.center);
		\draw [style={brown_thin}] (55.center) to (42.center);
		\draw [style={cyan_thin}] (32.center) to (55.center);
		\draw [style={magenta_thin}] (32.center) to (30.center);
		\draw [style={yellow_thin}] (56.center) to (6.center);
		\draw [style={yellow_thin}] (55.center) to (56.center);
		\draw [style={green_thin}] (26.center) to (58.center);
		\draw [style={red_thin}] (58.center) to (30.center);
		\draw [style={red_thin}] (29.center) to (58.center);
		\draw [style={green_thin}] (58.center) to (57.center);
		\draw [style={boundary}](11.center) to (36.center);
		\draw [style={boundary}](36.center) to (57.center);
		\draw [style={boundary}](57.center) to (47.center);

		\fill[bEV] (6) circle [radius=0.05];
		\fill[EV] (55) circle [radius=0.05];
		\fill[EV] (30) circle [radius=0.05];
		\fill[EV] (32) circle [radius=0.05];
		\fill[bEV] (25) circle [radius=0.05];

		\fill[gray,opacity=0.1] (40.center)--(41.center)--(56.center)--(6.center)--cycle;
		\fill[gray,opacity=0.1] (41.center)--(42.center)--(55.center)--(56.center)--cycle;
		\fill[gray,opacity=0.1] (42.center)--(54.center)--(32.center)--(55.center)--cycle;
		\fill[gray,opacity=0.1] (54.center)--(47.center)--(30.center)--(32.center)--cycle;
		\fill[gray,opacity=0.1] (47.center)--(57.center)--(58.center)--(30.center)--cycle;
		\fill[gray,opacity=0.1] (58.center)--(57.center)--(36.center)--(29.center)--cycle;
		\fill[gray,opacity=0.1] (29.center)--(36.center)--(11.center)--(10.center)--cycle;

		\fill[gray,opacity=0.1] (6.center)--(56.center)--(31.center)--(5.center)--cycle;
		\fill[gray,opacity=0.1] (56.center)--(55.center)--(32.center)--(31.center)--cycle;

		\fill[gray,opacity=0.1] (5.center)--(31.center)--(27.center)--(4.center)--cycle;
		\fill[gray,opacity=0.1] (31.center)--(32.center)--(26.center)--(27.center)--cycle;
		\fill[gray,opacity=0.1] (32.center)--(30.center)--(58.center)--(26.center)--cycle;
		\fill[gray,opacity=0.1] (26.center)--(58.center)--(29.center)--(28.center)--cycle;
		\fill[gray,opacity=0.1] (28.center)--(29.center)--(10.center)--(9.center)--cycle;

		\fill[gray,opacity=0.1] (4.center)--(27.center)--(25.center)--(3.center)--cycle;
		\fill[gray,opacity=0.1] (27.center)--(26.center)--(28.center)--(25.center)--cycle;
		\fill[gray,opacity=0.1] (25.center)--(28.center)--(9.center)--(8.center)--cycle;

		\fill[gray,opacity=0.2] (8.center)--(9.center)--(19.center)--(1.center)--(16.center)--(13.center)--(24.center)--(4.center)--(3.center)--(38.center)--(12.center)--(37.center)--(0.center)--(39.center)--cycle;
		\fill[gray,opacity=0.3] (9.center)--(10.center)--(20.center)--(2.center)--(17.center)--(14.center)--(23.center)--(5.center)--(4.center)--(24.center)--(13.center)--(16.center)--(1.center)--(19.center)--cycle;
		\fill[gray,opacity=0.4] (10.center)--(11.center)--(21.center)--(7.center)--(18.center)--(15.center)--(22.center)--(6.center)--(5.center)--(23.center)--(14.center)--(17.center)--(2.center)--(20.center)--cycle;
\end{tikzpicture}